\documentclass[a4paper]{amsart}

\usepackage{amsmath, amssymb, amsthm}
\usepackage[hidelinks]{hyperref} 
\usepackage[capitalise, noabbrev]{cleveref}
\usepackage{bbm}
\usepackage{comment}
\usepackage[shortlabels]{enumitem}  
\usepackage{esint, bigints}  
\usepackage{floatrow}
\usepackage{fullpage}
\usepackage{mathrsfs}
\usepackage{mathtools}
\usepackage{setspace}
\usepackage{stmaryrd}
\SetSymbolFont{stmry}{bold}{U}{stmry}{m}{n}
\usepackage{bm}  
\usepackage{tikz, tikz-cd, tikz-3dplot}
\usepackage{xcolor}

\newcommand*{\id}{\mathrm{id}}

\newcommand*{\hor}{\mathrm{hor}}

\newcommand*{\vol}{\mathrm{vol}}

\newcommand*{\C}{\mathbb{C}}

\newcommand*{\Heis}{\mathbb{H}}

\newcommand*{\R}{\mathbb{R}}

\renewcommand*{\S}{\mathbb{S}}
\newcommand*{\T}{\mathbb{T}}
\newcommand*{\Z}{\mathbb{Z}}

\newcommand*{\Cyl}{\mathbf{C}}
\newcommand*{\Exc}{\mathbf{E}}
\newcommand*{\IC}{\mathbf{I}}
\newcommand*{\ICmet}{\mathbf{I}^\mathrm{met}}
\newcommand*{\Mass}{\mathbf{M}}
\newcommand*{\NC}{\mathbf{N}}

\newcommand*{\Haus}{\mathcal{H}}

\newcommand*{\Leb}{\mathcal{L}}

\newcommand*{\Ord}{\mathcal{O}}

\newcommand*{\dd}{\mathrm{d}}
\newcommand*{\Dd}{\mathrm{D}}
\newcommand*{\del}{\partial}
\newcommand*{\eps}{\varepsilon}

\newcommand*{\Jac}{\mathcal{J}}
\newcommand*{\into}{\hookrightarrow}
\newcommand*{\Sym}{\mathrm{Sym}}
\newcommand*{\weakto}{\rightharpoonup}

\newcommand*{\der}[1]{\frac{\del}{\del #1}}

\newcommand{\ddt}{\left.\frac{\dd}{\dd t}\right|_{t=0}}

\DeclareMathOperator{\diam}{diam}
\DeclareMathOperator{\dist}{dist}

\DeclareMathOperator{\Lip}{Lip}

\DeclareMathOperator{\res}{\mathbin{\mbox{\Large$\llcorner$}}}

\DeclareMathOperator{\slot}{\mathbin{\mbox{\Large$\lrcorner$}}}
\DeclareMathOperator{\spt}{spt}
\DeclareMathOperator{\Span}{span}

\DeclarePairedDelimiter\DBrack{\llbracket}{\rrbracket}

\newtheorem{thm}{Theorem}[section]
\newtheorem{cor}[thm]{Corollary}
\newtheorem{lemma}[thm]{Lemma}

\newtheorem{prop}[thm]{Proposition}

\theoremstyle{definition}
\newtheorem{defn}[thm]{Definition}

\newtheorem{rmk}[thm]{Remark}
\newtheorem{assumption}[thm]{Assumption}

\numberwithin{equation}{section}

\newlist{steps}{enumerate}{1}
\setlist[steps, 1]{label={\underline{{Step} \arabic*.}}, leftmargin=*, wide, labelindent=0pt, itemsep=5pt}

\usepackage[backend=biber,style=alphabetic,sorting=nyt,doi=true,isbn=false,url=true,eprint=true,maxbibnames=99]{biblatex}
\renewbibmacro{in:}{}

\AtEveryBibitem{%
\ifentrytype{article}{
    \clearfield{url}%
    \clearfield{urldate}%
}{}
\ifentrytype{book}{
    \clearfield{url}%
    \clearfield{urldate}%
}{}
\ifentrytype{collection}{
    \clearfield{url}%
    \clearfield{urldate}%
}{}
\ifentrytype{incollection}{
    \clearfield{url}%
    \clearfield{urldate}%
}{}
}

\DeclareFieldFormat{edition}{#1 edition}

\def\BHeis{\mathcal{B}} 
\newcommand{\FKnorm}[1]{\| #1 \|_{\mathrm{FK}}}
\renewcommand{\Vec}[1]{\mathbf{#1}}
\newcommand{\FKdist}{\dist_\mathrm{FK}}

\def\idxiso{1}
\def\idxiterateiso{2}
\def\idximp{*}
\def\idxed{3}
\def\idxtilt{4}

\def\epsreg{\bar{\eps}}
\def\epshb{\eps_0}
\def\epsed{\eps_*}

\author{Gerard Orriols}
\address{ETH, Rämistrasse 101, 8092 Zürich, Switzerland}
\email{gerard.orriols@math.ethz.ch}
\date{June 13, 2024}

\title{Existence and partial regularity of Legendrian area-minimizing currents}
\begin{document}
\maketitle
\begin{abstract}
  We show that Legendrian integral currents in a contact manifold that locally minimize the mass among Legendrian competitors have a regular set which is open and dense in their support. We apply this to show existence and partial regularity of solutions of the Legendrian Plateau problem in the $n$th Heisenberg group for an arbitrary horizontal $(n-1)$-cycle as prescribed boundary, and of mass-minimizing Legendrian integral currents in any $n$-dimensional homology class of a closed contact $(2n+1)$-manifold.  In the case of the Heisenberg group, our result applies to Ambrosio--Kirchheim metric currents with respect to the Carnot--Carath\'eodory distance. Our results do not assume any compatibility between the subriemannian metric and the symplectic form.
\end{abstract}

\tableofcontents

\section{Introduction}
\label{sec:introduction}

The goal of this paper is to lay the foundations of an existence and regularity theory of area-minimizing currents among Legendrian integral currents in contact manifolds of arbitrary dimension. The main motivation for applying the methods of classical Geometric Measure Theory to Legendrian submanifolds comes from the project initiated by Schoen and Wolfson to develop a variational theory for the area of Lagrangian submanifolds in symplectic manifolds \cite{schoen-wolfson-lapietra, schoen-wolfson-jdg, schoen-wolfson-other}.

There are many geometric reasons to study minimizers of the area among Lagrangian submanifolds. It is natural to ask whether a homology class in a symplectic manifold which admits a Lagrangian representative contains one of minimal area which enjoys some partial regularity, as in the unconstrained case. Moreover, a good understanding of the homology minimization problem for Lagrangian currents seems to be fundamental before exploring the much more delicate topic of regularity of Hamiltonian-stationary Lagrangian surfaces and of existence of global area-minimizers in Hamiltonian isotopy classes, of which very little is known. See \cite{oh-hamiltonian-stationary} for an introduction to the problem and some motivating conjectures and \cite{viterbo-jams} for a lower bound on the area for Hamiltonian isotopy classes in some ambient spaces.

A deeper motivation to study this minimization problem comes from mirror symmetry and the SYZ conjecture: when the symplectic manifold is Calabi--Yau, hypothetically there should exist many special Lagrangian submanifolds, which are calibrated and hence absolute homology minimizers. Therefore one expects to be able to find them by minimizing area among Lagrangian integral currents in homology classes which admit them. A very short and elegant computation of Schoen and Wolfson \cite{schoen-wolfson-jdg} shows that this strategy can be realized, not only in the Calabi--Yau setting but more generally in K\"ahler--Einstein manifolds, and indeed Lagrangian-stationary closed surfaces are automatically minimal (if the ambient manifold is Calabi--Yau then they must be special Lagrangian). However, for their argument to work, it is fundamental to know that this object is regular enough, and in fact Micallef and Wolfson later found examples of homology classes where singularities cannot be avoided \cite{micallef-wolfson-singular, wolfson-singular}.

It was observed by Schoen and Wolfson, building on an example of Minicozzi \cite{minicozzi-willmore-lagrangian}, that Lagrangian area-minimizing surfaces can be very badly behaved unless they satisfy an exactness condition. More precisely, Minicozzi showed that cylinders of the form $\S^1_\eps \times \R \subset \C \times \C = \R^4$ are area-minimizing among Lagrangian surfaces and, unlike classical area minimizers, they do not satisfy a monotonicity formula for the area (not even a universal lower bound for the mass ratio in a ball). Unless one can somehow control the presence of regions of this type, developing a regularity theory in this generality seems far out of reach\footnote{Note however the a-priori estimates from \cite{schoen-wolfson-lapietra} which do not assume exactness of the surface but rely on the a-priori square integrability of the mean curvature, which already excludes many interesting singularities like the Schoen--Wolfson cones.}.

However, assuming that a smooth Lagrangian surface $\Sigma^2 \subset \R^4$ is exact (that is, the Liouville form $\lambda$, which satisfies $\dd \lambda = \omega$, restricts to zero on $\Sigma$), it has a lift $\tilde{\Sigma}$ as a surface in the Heisenberg group $\Heis^2 \simeq \R^5$, and the key insight of \cite{schoen-wolfson-jdg} was to show that this lift satisfies a monotonicity formula assuming just that $\Sigma$ is stationary under Hamiltonian deformations. This is due to the larger richness of variations and competitors available in the Legendrian setting. By extensively using this monotonicity formula and techniques exclusive of dimension two (such as a conformal parametrization, energy methods and holomorphic maps), Schoen and Wolfson succeeded in developing an optimal and rather complete regularity theory for minimizers of the mapping problem among Legendrian surfaces, and applied it to the Lagrangian setting. Their main result states that such minimizers are smooth surfaces except at a finite number of points, which are either of branch type (around which the parametrization is smooth) or of cone type (around which the parametrization is Lipschitz). The Hamiltonian-stationary two-dimensional cones are classified, and among them also the stable ones, but to the best of our knowledge it is still an open question to determine which of them are minimizing.

The extension to higher dimensions of the monotonicity formula for Hamiltonian-stationary Legendrian submanifolds was announced in \cite[Theorem 3]{schoen-wolfson-lapietra}. We provide a counterexample in \cref{sec:no-monotonicity}, showing that a different approach is needed and that the regularity of Hamiltonian-stationary surfaces in higher dimensions is more delicate than the regularity of Legendrian minimizing currents carried out here, where the lower bound in the mass density will come from an isoperimetric inequality.

Even in dimension two, the monotonicity formula of \cite{schoen-wolfson-jdg} is not explicit, and its construction goes through the solution of an auxiliary hyperbolic PDE. Note however that a newer simplified almost monotonicity formula has been recently found by Rivi\`ere \cite{riviere-monotonicity}, still in dimension two, and extended to parametrized Hamiltonian-stationary varifolds by Rivi\`ere and Pigati \cite{pigati-riviere-compactness}. This monotonicity formula has been applied by Rivi\`ere to solve min-max problems for Legendrian surfaces in Sasakian $5$-manifolds \cite{riviere-area-variations}, thus generalizing the pure minimization problem.

In addition to that, the proof of the decay of the excess in \cite{schoen-wolfson-jdg} has the added difficulty of controlling the parametrization of the surface at each step, which makes it quite technical. In fact, an attempt to understand and simplify their arguments with the help of newly available tools was the starting point of this work.

In contrast to the slower progress in the GMT problem in arbitrary dimensions, the last few years have seen many new results about the PDE that governs the local picture for graphs. Since Legendrian graphs are locally lifts of Lagrangian graphs, and Lagrangian variations agree locally with Hamiltonian variations, the corresponding Euler--Lagrange equation is precisely the Hamiltonian-stationary equation. This equation has a very special structure when the ambient manifold is $\C^n$ and has been studied by Bhattacharya, Chen and Warren \cite{chen-warren-Cn, bhattacharya-warren}. More recently, these authors have noticed that partial regularity results can be obtained without relying on this special structure, and therefore hold on arbitrary symplectic manifolds \cite{bhattacharya-chen-warren}.

Even more recently, Bhattacharya and Skorobogatova \cite{hessian-partial-regularity} have shown, still in the graphical setting, that smallness of the average oscillation of the Hessian of the potential $u$ in a ball, which corresponds to smallness of the excess of the graph in a cylinder, is enough to deduce regularity. Nevertheless they need to know a priori that $\|\Dd^2 u\|_{L^\infty} \leq 1 - \eta$ to guarantee that the area functional is convex, and their required smallness of the excess depends on this $\eta$. Inspired by their theorem, in our main result, \cref{thm:eps-reg}, we remove the Lipschitz graph assumption and, working with currents, we prove an $\eps$-regularity theorem for Legendrian area-minimizers analogous to the classical ones of De Giorgi \cite{de-giorgi-frontiere-orientate} and Almgren \cite{almgren-regularity-elliptic} for minimal surfaces.

At the same time, as with high-codimension minimal surfaces, the Euler--Lagrange PDE for graphs lacks an existence theory of weak solutions and a satisfactory functional space to search for minimizers. The framework of Legendrian currents developed here provides a solution to these issues and furthermore allows for boundary data which cannot be expressed as a graph.

\subsection{Main results}

We state here our main existence and regularity results. The relevant notions will be introduced in the next two sections. For now, note that $\Heis^n$ denotes the Heisenberg group of topological dimension $2n+1$ with its standard contact structure and metric on the horizontal planes; a horizontal integral current is just an integral current almost all of whose tangent planes are isotropic, and a Legendrian current is a horizontal current of dimension $n$.

\begin{thm}[Legendrian Plateau problem]
  \label{thm:legendrian-plateau-problem}
  Let $S$ be a compactly supported $(n-1)$-dimensional horizontal integral current in $\Heis^n$ with $\del S = 0$. Then there exists a compactly supported Legendrian integral current $T$ with $\del T = S$ such that $\Mass(T) \leq \Mass(T')$ for any other Legendrian integral current $T'$ with $\del T' = S$.

  Moreover there exists an open set $\mathcal{U} \subset \Heis^n$ such that $\spt T \cap \mathcal{U}$ is a real-analytic embedded $n$-dimensional Legendrian submanifold of $\Heis^n$ and $\spt T \cap \mathcal{U}$ is dense in $\spt T \setminus \spt S$.
\end{thm}

  This theorem can also be stated in terms of Ambrosio--Kirchheim metric currents with respect to the intrinsic Carnot--Carath\'eodory distance of $\Heis^n$, but we have preferred to work with Legendrian Federer--Fleming currents to take advantage of the homological theory which is available for them, of the larger literature that uses them, and of their better compatibility with Lagrangian currents in a symplectic manifold. Nonetheless we have to use metric currents crucially in a portion of our arguments, and we have recast our results into that language in \cref{rmk:recast-metric-currents}.

We also remark that the Legendrian Plateau problem had already been studied by Minicozzi in his PhD thesis \cite{minicozzi-phd} and a penalization approach to solve it is explained in \cite{schoen-wolfson-other}.

The assumption that the ambient manifold is the Heisenberg group is not essential here, and the Plateau problem can be solved as well in compact manifolds or in noncompact manifolds with contolled geometry at infinity. For closed manifolds we have the following result\footnote{Here we assume that the ambient manifold is smooth---see \cref{sec:main-regularity} for more precise versions of the regularity statements under weaker assumptions.}:

\begin{thm}[Legendrian homology minimizers]
  \label{thm:legendrian-homology-problem}
  Let $(M^{2n+1}, \Xi, g)$ be a smooth closed contact manifold with contact distribution $\Xi$ and a smooth subriemannian metric $g$ defined on $\Xi$. Let $\mathfrak{a} \in H_n(M, \Z)$ be a homology class. Then there exists a Legendrian integral current $T$ representing $\mathfrak{a}$ such that $\Mass(T) \leq \Mass(T')$ for any other Legendrian integral current $T'$ in $\mathfrak{a}$.

  Moreover there exists an open set $U \subset M$ such that $\spt T \cap U$ is a smooth embedded $n$-dimensional Legendrian submanifold of $M$ and $\spt T \cap U$ is dense in $\spt T$.
\end{thm}

The proof, both of the existence and (mainly) of the partial regularity, makes use of the isoperimetric and coning inequalities for horizontal currents in $\Heis^n$ which have been recently established by Basso, Wenger and Young \cite{basso-wenger-young-fillings} in the context of metric currents. These results give us the basic ingredients to build up a regularity theory in the absence of a monotonicity formula, as pioneered by Almgren for currents which minimize an elliptic functional \cite{almgren-regularity-elliptic}. In fact, our regularity theorem also works in great generality, without assuming any compatibility condition between the metric and the symplectic form (in particular, an almost complex structure may not exist anywhere), and with minor modifications can be extended to functionals more general than an anisotropic area.

Although the ideas of the regularity proof were pioneered by Almgren and De Giorgi, we follow more closely the simplifications introduced by Schoen--Simon \cite{schoen-simon-anisotropic} and Ambrosio--De Lellis--Schmidt \cite{regularity-currents-hilbert} in the proof of the decay of the excess, approximating our current by (roughly) the graph of the gradient of the solution of a fourth order linear PDE and using classical $L^p$ estimates to control the error terms.

We remark that, whereas the existence theory works as well for horizontal currents of any dimension $k < n$ (see Theorems \ref{thm:existence-plateau} and \ref{thm:existence-homology}), our regularity theory relies crucially on the representation of our current as the graph of a gradient over a Legendrian plane, which is only available when $k = n$. Even if the surface is given as a graph over a horizontal $k$-plane, we are not aware of any results in the literature that deal with the corresponding system of PDE, in contrast with the larger literature on the hamiltonian stationary equation. In the parametrized case, however, Qiu \cite{qiu-isotropic} successfully extended the theory of Schoen and Wolfson to isotropic surfaces in any dimension. In our opinion, the extension of these results to higher ambient dimensions is an interesting problem.

We believe that a more precise description of the singular set should be possible, at least in two dimensions and in the Heisenberg group. The key difficulty here would be to develop a blowup analysis around branch points. It is natural to conjecture that one would recover the result obtained by Schoen and Wolfson in the parametrized setting: a locally finite set of singular points consisting of branch points and Schoen--Wolfson conical singularities.

We also hope to be able to apply our results soon to variational problems for Lagrangian surfaces; in particular there are some issues to be addressed regarding the notions of exactness and Legendrian lift for integral currents. But beyond this application, we believe that our result has many other points of interest. On one hand, contact geometry has become a subject on its own and the study of canonical representatives of their homology is a natural and basic question to study. On the other hand, from the perspective of metric geometry, we provide one of the very few available partial regularity results for high codimension minimal surfaces in non-smooth metric spaces (compare with \cite{regularity-currents-hilbert} for Hilbert spaces and with \cite{stadler-minimal-cat0} for parametrized two-dimensional surfaces in $\operatorname{CAT}(0)$ spaces). In addition to that, the techniques developed here may be applied to other geometric variational problems in which Legendrian currents appear as unit normal bundles---see the survey \cite{bernig-survey-curvature} and the references therein for more details.

\subsection{Structure of the paper}

After introducing some preliminary notation and facts in \cref{sec:preliminaries}, the relevant notions of currents will be presented in \cref{sec:existence}. The existence part of the Theorems \ref{thm:legendrian-plateau-problem} and \ref{thm:legendrian-homology-problem} is also stated and proved in \cref{sec:existence}, see \cref{thm:existence-plateau} and \cref{thm:existence-homology} respectively. The regularity part is made precise in \cref{sec:main-regularity}: in \cref{cor:regularity-open-dense} we establish that the regular set (where the support of the current is $C^{1,1/2}$) is dense in the support. Then the higher regularity stated here follows from \cref{thm:regular-is-smooth} for smooth manifolds and \cref{cor:real-analytic} for $\Heis^n$.

The main theorem of the present article, from which the partial regularity is deduced, is the $\eps$-regularity \cref{thm:eps-reg}. Its proof takes most of the work and spans Sections \ref{sec:height-bound-lipschitz} and \ref{sec:biharmonic-decay}. Finally, \cref{sec:metric-geometry} contains all the statements that we use from metric geometry and serves as a dictionary between our currents and Carnot--Carath\'eodory metric currents, and \cref{sec:no-monotonicity} contains a counterexample of the monotonicity of the mass ratios in dimension larger than two.

\subsection{Acknowledgments}
The author thanks his PhD supervisor, Tristan Rivi\`ere, for his support and encouragement to study the work of Schoen and Wolfson and attempt to simplify and generalize parts of it. He also wishes to thank Federico Franceschini, Filippo Gaia, Matilde Gianocca and Tom Ilmanen for inspiring conversations and for their interest in the work.

\section{Preliminaries and notation}
\label{sec:preliminaries}

\subsection{The Heisenberg group}

Most of our analysis will take place in the Heisenberg group $\Heis^n$, because thanks to the Darboux theorem (\cref{prop:contact-darboux}), it is a local model for the geometry of a contact manifold with good homogeneity properties. Here we summarize the main geometric and analytical facts that we will use in the sequel. For the notation, we use mostly the conventions of \cite{introduction-heisenberg}; we also refer to this book for further information and references.

We will often identify $\Heis^n$ with $\R^{2n+1}$ by using the globally defined exponential coordinates $(\Vec{z}, \varphi) = (\Vec{x}, \Vec{y}, \varphi)$. Here $\Vec{z}$ denotes an element of $\C^n$, which we will identify with $\R^{2n}$, and $\Vec{x}, \Vec{y}$ denote $n$-tuples of real numbers corresponding to the real and imaginary part of $\Vec{z}$, respectively.

In these coordinates, the group operations are
\begin{equation}
  \label{eq:heis-product}
  (\Vec{x}, \Vec{y}, \varphi) \cdot (\Vec{x}', \Vec{y}', \varphi')
  = \left(\Vec{x} + \Vec{x}', \Vec{y} + \Vec{y}', \varphi + \varphi' + \frac{1}{2} \left( \Vec{x} \cdot \Vec{y}' - \Vec{y} \cdot \Vec{x}'\right) \right)
\end{equation}
and
\begin{equation}
  \label{eq:heis-inverse}
  (\Vec{x}, \Vec{y}, \varphi)^{-1} = (-\Vec{x}, -\Vec{y}, -\varphi).
\end{equation}

The Heisenberg group is naturally a contact manifold with global contact form
\begin{equation}
  \label{eq:contact-form-heisenberg}
  \theta := \dd \varphi - \frac{1}{2} (\Vec{x} \cdot \dd \Vec{y} - \Vec{y} \cdot \dd \Vec{x}),
\end{equation}
which defines a hyperplane distribution $\Xi = \Xi_{\Heis}$ by $\Xi_\xi := \ker \theta(\xi)$. As all Lie groups, the Heisenberg group is homogeneous and parallelizable. In particular, left translations allow us to identify all tangent spaces. Moreover, it can be checked that $\theta$ is left-invariant and therefore the horizontal bundle $\Xi_\Heis$ can also be trivialized by using left translations. In what follows, we will use this observation extensively and identify the horizontal spaces at all points with $\Xi_0$.

We endow $\Heis^n$ with a subriemannian metric $g_\Heis$, which we will also write as $\langle \cdot, \cdot \rangle$, defined as follows: take $\langle \cdot, \cdot \rangle$ on $\Xi_0$ to be the standard metric of $\C^n = \R^{2n} \simeq \Xi_0 \subset T_0 \Heis^n$, and extend it to the whole $\Xi_\Heis$ in a left-invariant way.

For any horizontal $k$-plane $\pi$ ($1 \leq k \leq 2n$) and any smooth function $f$, we denote by $\nabla^\pi f$ the unique vector in $\pi$ such that $\langle v, \nabla^\pi f \rangle = \dd f (v)$ for any vector $v \in \pi$. In the case that $\pi$ is the whole horizontal $2n$-space $\Xi$, we denote $\nabla^\pi f$ as $\nabla^H f$. Then $\nabla^\pi f$ is the (Riemannian) orthogonal projection of $\nabla^H f$ onto $\pi$ with respect to the Riemannian metric on $\Xi$.

For the coordinate functions, it is easy to check that
\begin{equation}
  \nabla^H x^i = \der{x^i} - \frac{1}{2} y^i \der{\varphi}
  \qquad  \text{and} \qquad
  \nabla^H y^i = \der{y^i} + \frac{1}{2} x^i \der{\varphi},
\end{equation}
since all vector fields involved are left-invariant and agree at the origin, where they form the standard orthonormal basis of $\R^{2n}$. Thus, by definition, $\{\nabla^H x^i, \nabla^H y^i\}_{i=1,\ldots,n}$ form an orthonormal basis of $\Xi$ for the standard subriemannian metric.

As is customary, we will work with the Folland--Kor\'anyi gauge
\begin{equation}
  \label{eq:fk-gauge}
  \tau(\Vec{z}, \varphi) := \FKnorm{(\Vec{z}, \varphi)} := \left(|\Vec{z}|^4 + 16 \varphi^2 \right)^{\frac{1}{4}}.
\end{equation}
Cygan \cite{cygan-subadditivity-koranyi} proved that $\FKnorm{\xi \xi'} \leq \FKnorm{\xi} + \FKnorm{\xi'}$ and therefore it defines a left-invariant distance
\begin{equation}
  \FKdist(\xi, \xi') := \FKnorm{\xi^{-1} \xi'}
\end{equation}
whose balls will be denoted by $\BHeis_r := \{ \xi \in \Heis^n : \FKnorm{\xi} < r \}$. This distance is (bilipschitz) equivalent to the Carnot--Carath\'eodory distance, but for convenience we will use the former most of the time. In particular, it is smooth away from the origin and enjoys the following useful property: the function $\tau : \Heis^n \to \R$, $\tau(\xi) = \FKnorm{\xi},$ satisfies
\begin{equation}
  \label{eq:hor-grad-tau-bound}
  |\nabla^\pi \tau| \leq |\nabla^H \tau| = \frac{|\Vec{z}|}{\tau} \leq 1
\end{equation}
for every horizontal $k$-plane $\pi$. This follows from a computation:
\begin{align*}
  (4 \tau^3 |\nabla^H \tau|)^2
  &= |\nabla^H \tau^4|^2
  = |\nabla^H |\Vec{z}|^4 + 16 \nabla^H \varphi^2|^2 \\
  &= |\nabla^H |\Vec{z}|^4|^2 + 2^8 |\nabla^H \varphi^2|^2 + 32 \langle \nabla^H |\Vec{z}|^4, \nabla^H \varphi^2 \rangle \\
  &= 4 |\Vec{z}|^4 |\nabla^H |\Vec{z}|^2|^2 + 2^{10} \varphi^2 |\nabla^H \varphi|^2 + 2^7 |\Vec{z}|^2 \varphi \langle \nabla^H |\Vec{z}|^2, \nabla^H \varphi \rangle \\
  &= 16 |\Vec{z}|^4 |\Vec{z}|^2 + 2^{8} \varphi^2 |\Vec{z}|^2
  = 16 |\Vec{z}|^2 \tau^4.
\end{align*}

\subsection{Contact manifolds}
For simplicity, we will follow \cite{blair-contact-manifolds} and will consider only contact manifolds admitting a global contact form. All our arguments generalize easily to manifolds that are not co-orientable. Hence a \emph{contact manifold} is a smooth manifold $M^{2n+1}$ equipped with a $1$-form $\theta$ such that $\theta \wedge (\dd \theta)^n \neq 0$ everywhere. In fact, we will only be interested in the $2n$-dimensional plane distribution $\Xi$ given by $\ker \theta$, so all of the geometric notions discussed below will be invariant under scaling $\theta$ by a nonvanishing function.

It is clear that the $1$-form $\theta_{\Heis^n}$ from \eqref{eq:contact-form-heisenberg} defines a global contact form on $\Heis^n$, since
\[
  \theta_{\Heis^n} \wedge \left(\dd \theta_{\Heis^n}\right)^n
  = \theta_{\Heis^n} \wedge \left(\sum_k \dd x^k \wedge \dd y^k \right)^n
  = n! \, \dd \varphi \wedge \dd x^1 \wedge \dd y^1 \wedge \cdots \wedge \dd x^n \wedge \dd y^n \neq 0.
\]

To study the volume of submanifolds, we will need a Riemannian metric $g$ defined on the planes $\Xi$, a so-called subriemannian structure on $M$.
Then we can state a subriemannian version of Darboux's theorem as follows:

\begin{prop}
  \label{prop:contact-darboux}
  For every point $\xi$ of a contact manifold $(M^{2n+1}, \theta)$ there is a neighborhood $U \ni \xi$, an open set $0 \in \mathcal{U} \subset \Heis^n$ and a diffeomorphism $\phi : U \to \mathcal{U}$ such that $\theta = \phi^* \theta_{\Heis^n}$.

  Moreover, if $M$ is compact and $g$ is a Riemannian metric on the contact distribution $\Xi = \ker \theta$, there exist constants $0 < \lambda \leq \Lambda < \infty$ such that
  \begin{equation}
    \label{eq:comparison-metric}
    \lambda \phi^* g_{\Heis^n} \leq g \leq \Lambda \phi^* g_{\Heis^n}.
  \end{equation}
  \begin{proof}
    The first part is just a restatement of the classical Darboux theorem for contact manifolds, see for example \cite{blair-contact-manifolds}. The second part just follows from the compactness of $M$ by covering it with finitely many Darboux charts.
  \end{proof}
\end{prop}

\begin{rmk}
  Notice that if we assume that the metric $g$ and the symplectic form $\omega := \dd \theta |_\Xi$ are compatible, in the sense that there exists an almost complex structure $J_\xi$ on $\Xi_\xi$ smoothly depending on $\xi$ such that $\omega(u, v) = g(J u, v)$ for every $u, v \in \Xi$, then we can make $\lambda$ and $\Lambda$ as close to one as wished by making $U$ smaller. Indeed, $(\phi^{-1})^* \dd \theta = \dd \theta_{\Heis}$ is the standard symplectic form $\omega_{\Heis}$ on the horizontal plane at $0 \in \Heis^n$, and $(\phi_* J)_0$ is a complex structure compatible with it. Therefore extending the scalar product $\omega_{\Heis}(\cdot , (\phi_* J)_0 \cdot)$ to a left-invariant subriemannian metric on $\Heis^n$ makes the resulting space isometric to the standard subriemannian $(\Heis^n, g_{\Heis})$, and $D \phi_{\phi^{-1}(0)}$ is a linear isometry onto $\Xi_0$.

However we will not need this hypothesis for our regularity result. This makes the problem more anisotropic and our setting can be compared with that of Almgren's regularity theorem for minimizers of elliptic functionals \cite{almgren-regularity-elliptic, schoen-simon-anisotropic}.
\end{rmk}

\subsection{Biharmonic functions}
Here we collect some of the classical estimates for solutions to elliptic fourth order equations with constant coefficients:
\begin{prop}
  \label{prop:biharmonic}
  Let $(a_{ik}^{jl})_{i,j,k,l=1}^n$ be real numbers satisfying the boundedness condition $|a_{ik}^{jl}| \leq \Lambda_0$ and the Legendre ellipticity condition, $a_{ik}^{jl} \sigma_{ij} \sigma_{kl} \geq \lambda_0 |\sigma|^2$ for every $n \times n$ symmetric matrix $(\sigma_{ij})$.

  Then for any $f \in W^{2, \infty}(B_R)$ there exists a unique function $u : B_R \to \R$ which is in $W^{2,p}(B_R)$ for every $1 \leq p < \infty$ and solves the Dirichlet problem
  \begin{equation}
    \begin{cases}
      a_{ik}^{jl} \del_{ijkl} u = 0 & \text{in } B_R \\
      u = f & \text{on } \del B_R \\
      \del_\nu u = \del_\nu f & \text{on } \del B_R.
    \end{cases}
  \end{equation}
  Moreover, $u$ is smooth inside $B_R$ and satisfies the interior estimates
  \begin{equation}
    \label{eq:interior-bound-hessian}
    \sup_{B_{r}} |\Dd^{2+k} u|^2 \leq C_k \frac{1}{(R-r)^{n+2k}} \int_{B_R} |\Dd^2 u|^2
    \qquad \text{for any } 0 < r < R \text{ and any integer } k \geq 0,
  \end{equation}
  \begin{equation}
    \label{eq:interior-closeness-hessian}
    \sup_{B_{r}} |\Dd^2 u - \Dd^2 u(0)|^2 \leq C \frac{r^2}{R^n} \int_{B_R} |\Dd^2 u|^2
    \qquad \text{for any } 0 < r < \frac{R}{2},
  \end{equation}
  the global $L^p$ estimate
  \begin{equation}
    \label{eq:lp-estimate-biharmonic}
    \int_{B_R} |\Dd^2 u|^p \leq C_p \int_{B_R} |\Dd^2 f|^p
    \qquad \text{for any } 1 < p < \infty
  \end{equation}
  and Agmon's maximum estimate for the derivative:
  \begin{equation}
    \label{eq:agmon}
    \sup_{B_R} |\Dd u| \leq C \sup_{\del B_R} |\Dd u| \leq C \sup_{B_R} |\Dd f|.
  \end{equation}
  Here the constants depend only on $n, \Lambda_0 / \lambda_0$ and their subindices.
  \begin{proof}
    The existence, uniqueness and interior estimates are classical, see for example \cite{agmon-book}. The estimates up to the boundary follow from \cite{agmon-douglis-nirenberg}, and Agmon's maximum principle originally appeared in \cite{agmon-maximum}. See also the book \cite{polyharmonic-book} and the references therein.
  \end{proof}
\end{prop}

\section{Horizontal currents and existence of minimizers}
\label{sec:existence}

In order to do calculus of variations with Legendrian submanifolds, as in the classical Plateau problem, we need to enlarge our class of objects to a space that includes weaker objects but enjoys better compactness properties. Based on the success of the Federer--Fleming theory of integral currents \cite{federer-fleming}, which is based on the De Rham complex of differential forms, Franchi, Serapioni and Serra Cassano introduced in \cite{regular-submanifolds-heis} a class of objects in the Heisenberg group that, in low dimensions, generalizes horizontal submanifolds. As with Federer--Fleming currents, these objects are defined in duality with a differential complex introduced by Rumin \cite{rumin} in any contact manifold.

Here we work in an arbitrary contact manifold $M^{2n+1}$ and take our class of horizontal currents to be the Rumin currents of dimension at most $n$. In these dimensions, these currents can be characterized as Federer--Fleming currents in which the contact form and its exterior derivative restrict to zero. In the case of the Heisenberg group, integral Rumin currents correspond to integral metric currents, in the sense of Ambrosio--Kirchheim \cite{ambrosio-kirchheim}, with respect to the Carnot--Carath\'eodory metric (see \cref{sec:metric-geometry} for more details).

For Federer--Fleming currents, we follow the conventions of \cite{simon-book}; in particular we allow our integral currents to have globally infinite mass.

\subsection{The Rumin complex and horizontal currents}

We begin by recalling the Rumin complex in low dimensions. Let $U \subset M$ be an open set in which $\Xi$ is determined by the $1$-form $\theta$, and let $0 \leq k \leq n$. We define the differential ideal generated by $\Xi$ as
\begin{equation}
  \mathscr{I}^k(U) := \{ \theta \wedge \alpha + \dd \theta \wedge \beta : \alpha \in \mathcal{D}^{k-1}(U), \beta \in \mathcal{D}^{k-2}(U) \} \subset \mathcal{D}^k(U).
\end{equation}
Notice that $\mathscr{I}^k$ only depends on $\Xi$ and not on $\theta$, and thus makes sense even if $\Xi$ is not the kernel of a single $1$-form in $U$.
The (lower half of the) Rumin complex is the chain complex
\begin{equation}
  \mathcal{D}^0_\mathrm{R}(U) \xrightarrow{\dd} 
  \mathcal{D}^1_\mathrm{R}(U) \xrightarrow{\dd} 
  \cdots
  \xrightarrow{\dd} \mathcal{D}^{n-1}_\mathrm{R}(U)
    \xrightarrow{\dd} \mathcal{D}^{n}_\mathrm{R}(U)
\end{equation}
where
\begin{equation}
  \mathcal{D}^k_\mathrm{R}(U) := \mathcal{D}^k(U) / \mathscr{I}^k(U)
\end{equation}
and the differential $\dd$ clearly descends to the quotients.

Notice that a linear map $T : \mathcal{D}^k_\mathcal{R}(U) \to \R$ can be identified with a map $T : \mathcal{D}^k(U) \to \R$ which vanishes on $\mathscr{I}^k$. Thus it is natural to define

\begin{defn}[Horizontal current]
  \label{defn:horizontal-current}
  A De Rham current $T \in \mathcal{D}_k(U)$ is called \emph{horizontal} if
  \begin{equation}
    \label{eq:def-horizontal}
    T(\theta \wedge \alpha + \dd \theta \wedge \beta) = 0
  \end{equation}
  for any $\alpha \in \mathcal{D}^{k-1}(U)$ and $\beta \in \mathcal{D}^{k-2}(U)$. We denote by $\mathcal{D}^\hor_k(U)$ the set of horizontal currents in $U$, and we call them \emph{Legendrian} when $k = n$.
\end{defn}

Horizontal submanifolds (or currents) are sometimes called \emph{isotropic} in the literature.
It is immediate from the definition that if $T \in \mathcal{D}^\hor_k(U)$, then $\del T \in \mathcal{D}^\hor_{k-1}(U)$. It is also clear that $\mathcal{D}^\hor_k(U)$ is closed in $\mathcal{D}_k(U)$ with respect to weak convergence. Next we need a subriemannian notion of mass. We start with the following definition:

\begin{defn}[Horizontal mass and comass]
  Let $\omega \in \bigwedge^k T^*_\xi M$. We define its \emph{horizontal comass} as
  \begin{equation}
    \| \omega \|^\hor_* := \sup \{ \omega(v_1, \ldots, v_k) : v_1, \ldots, v_k \in \Xi_\xi \subset T_\xi M, |v_1 \wedge \cdots \wedge v_k|_g \leq 1 \}.
  \end{equation}
  For a current $T \in \mathcal{D}^\hor_k(U)$, we let its \emph{horizontal mass} be
  \begin{equation}
    \Mass^\hor(T) := \sup \{ T(\omega) : \omega \in \mathcal{D}^k(U), \| \omega \|^\hor_* \leq 1 \text{ in } U\}.
  \end{equation}
\end{defn}

In order to make use of the Federer--Fleming theory it will be convenient to extend the metric $g$, only defined on the distribution $\Xi$, to a Riemannian metric $g_0$ on the whole $TM$ and then to embed $(M, g_0)$ into some Euclidean space $\R^L$ using Nash's theorem. The following proposition shows that the relevant notions are independent of this extension.

\begin{prop}
  \label{prop:horizontal-mass}
  The horizontal mass enjoys the following properties:
  \begin{enumerate}[(i)]
    \item For any current $T \in \mathcal{D}_k$, $\Mass(T) \leq \Mass^\hor(T)$ holds.
    \item Given $T \in \mathcal{D}_k$ with $\Mass(T) < \infty$, we have $\Mass^\hor(T) < \infty$ if and only if $T(\theta \wedge \alpha) = 0$ for any $\alpha \in \mathcal{D}^{k-1}$. Moreover, in that case $\Mass(T) = \Mass^\hor(T)$.
    \item A current $T \in \mathcal{D}_k$ satisfies $\Mass^\hor(T) + \Mass^\hor(\del T) < \infty$ if and only if $T \in \mathcal{D}^\hor_k$ and $\Mass(T) + \Mass(\del T) < \infty$.
  \end{enumerate}
  In particular, $\Mass(T)$ and $\Mass(\del T)$ do not depend on the extension $g_0$ of $g$ for $T \in \mathcal{D}^\hor_k$.
  \begin{proof}
    Clearly $\| \omega \|^\hor_* \leq \| \omega \|_*$, so $\Mass(T) \leq \Mass^\hor(T)$ for any current $T$, which proves (i).

    We now show (ii): suppose first that $\Mass^\hor(T) < \infty$. Since $\|\theta \wedge \alpha\|^\hor_* = 0 < 1$, we have that $|T(s \theta \wedge \alpha)| \leq \Mass^\hor(T)$ for any $s > 0$, so by letting $s \to \infty$ it follows that $T(\theta \wedge \alpha) = 0$.

    To prove the opposite implication we will show that $\Mass^\hor(T) \leq \Mass(T)$, which together with (i) yields the equality. First let $X$ be the smooth vector field everywhere $g_0$-orthogonal to $\Xi$ such that $\theta(X) \equiv 1$. Given $\omega \in \mathcal{D}^k$ with $\| \omega \|^\hor_* \leq 1$, consider the $k$-form $\omega' := X \slot (\theta \wedge \omega)$.

    We claim that $\| \omega' \|_* \leq 1$ everywhere: indeed, fix $\xi \in M$ and let $v_1 \wedge \cdots \wedge v_k \in \bigwedge^k T_\xi M$ span the $k$-plane $V$ with $|v_1 \wedge \cdots \wedge v_k|_{g_0} \leq 1$. If $V$ is not horizontal, by choosing an orthonormal basis of the $(k-1)$-plane $V \cap \Xi_\xi$ and completing it to an orthonormal basis of $V$, we can suppose that $v_1, \ldots, v_{k-1} \in \Xi_\xi$ and write $v_k = v_k' + s X$ for some $v_k' \in \Xi_\xi$ and $s \in \R$. Then, since $v_1 \wedge \cdots \wedge v_{k-1} \wedge v_k'$ is $g_0$-orthogonal to $v_1 \wedge \cdots \wedge v_{k-1} \wedge sX$, it holds that
  \[
    |v_1 \wedge \cdots \wedge v_{k-1} \wedge v_k'|_g \leq |v_1 \wedge \cdots \wedge v_k|_{g_0} \leq 1,
  \]
  which implies that
  \[
    (X \slot \theta \wedge \omega)(v_1, \ldots, v_k)
    = (\theta \wedge \omega)(X, v_1, \ldots, v_k)
    = (\theta \wedge \omega)(X, v_1, \ldots, v_k')
    = \theta(X) \omega(v_1, \ldots, v_k'),
  \]
  hence
  \[
    \omega'(v_1, \ldots, v_k)
    \leq |\omega(v_1, \ldots, v_k')|
    \leq |v_1 \wedge \cdots \wedge v_k'|_g
    \leq 1
  \]
  and the claim is proven. Now observe that
  \[
    \omega' = X \slot (\theta \wedge \omega) = (X \slot \theta) \wedge \omega - \theta \wedge (X \slot \omega),
  \]
  which implies that for any $T$ satisfying $T(\theta \wedge \cdot) \equiv 0$,
  \[
    T(\omega)
    = \theta(X) T(\omega) - T(\theta \wedge (X \slot \omega))
    = T(\omega') \leq \Mass(T)
  \]
  and $\Mass^\hor(T) \leq \Mass(T)$ follows. Finally, (iii) is immediate from (ii) and the definition of $\mathcal{D}^\hor_k$.
  \end{proof}
\end{prop}

\begin{defn}
  The groups of normal and integral horizontal currents in $U \subset M$ are
  \begin{equation}
    \NC^\hor_k(U) := \NC_k(U) \cap \mathcal{D}^\hor_k(U)
    \qquad
    \text{and}
    \qquad
    \IC^\hor_k(U) := \IC_k(U) \cap \mathcal{D}^\hor_k(U).
  \end{equation}
\end{defn}

The Federer--Fleming compactness theorem \cite[Theorem~27.3]{simon-book} immediately gives:
\begin{prop}
  \label{prop:ff-compactness}
  If $\{T_j\} \subset \IC^\hor_k(U)$ is a sequence with $\|T_j\|(W) + \|\del T_j\|(W) \leq C(W)$ for every $W \Subset U$, then a subsequence converges weakly to a current $T \in \IC^\hor_k(U)$.
\end{prop}

It is well known that for a rectifiable current $T$ and $\|T\|$-almost every point $p$, if $\eta_{p,\rho}$ denotes the map $q \mapsto \frac{q-p}{\rho}$ in any coordinates around $p$, then $(\eta_{p, \rho})_\# T \weakto \Theta(\|T\|, p) T_p$ as $\rho \searrow 0$, where $T_p$ is the current induced by a plane called the tangent plane at $p$. For horizontal currents we have the following characterization:

\begin{prop}
  \label{prop:tangent-planes-isotropic}
  Let $T$ be a rectifiable $k$-current in $U$. Then $T \in \IC^\hor_k(U)$ if and only if for $\|T\|$-almost every $p \in U$, the approximate tangent plane $T_p$ is isotropic, which means that both $\theta$ and $\dd \theta$ vanish on it.
  \begin{proof}
    If all tangent planes are isotropic, then \eqref{eq:def-horizontal} is immediate by the integral representation formula for $T$. For the reciprocal, we just show that $\theta_p$ vanishes on $T_p$ (the computation for $(\dd \theta)_p$ is similar). Since the statement is local and invariant by diffeomorphisms, we may suppose that $T$ is a current in a neighborhood of $0 \in \R^{m}$ and $p = 0$ (in our case, $m = 2n+1$, and clearly the statement is much more general).
    Let $\alpha \in \mathcal{D}^{k-1}(\R^m)$ and write, for short, $\eta_\rho := \eta_{0,\rho}$. Suppose that $\alpha$ is supported in $B_R$. Then
    \[
      T_p(\theta_p \wedge \alpha)
      = \lim_{\rho \searrow 0} (\eta_\rho)_\# T (\theta_p \wedge \alpha)
      = \lim_{\rho \searrow 0} T (\eta_\rho^*\theta_p \wedge \eta_\rho^*\alpha)
      = \lim_{\rho \searrow 0} T ((\eta_\rho^*\theta_p - \theta) \wedge \eta_\rho^*\alpha)
    \]
    and thus
    \begin{align*}
      \left| T_p(\theta_p \wedge \alpha) \right|
      &\leq \limsup_{\rho \searrow 0} \|T\|(B_{\rho R}) \sup_{B_{\rho R}}|\eta_\rho^*\theta_p - \theta| \sup_{\rho R} |\eta_\rho^*\alpha| \\
      &\leq \limsup_{\rho \searrow 0} \|T\|(B_{\rho R}) C \rho^{1-k}
      \leq \limsup_{\rho \searrow 0} C \rho = 0.
      \qedhere
    \end{align*}
  \end{proof}
\end{prop}

As a consequence, in the coarea formula for horizontal currents we can replace the gradient of the slicing function by its horizontal gradient. Not only that, but thanks to the work of Ambrosio and Kirchheim \cite{ambrosio-kirchheim, ambrosio-kirchheim-rectifiable} and the discussion in \cref{sec:metric-geometry}, we can slice by any function which is Lipschitz with respect to the intrinsic Carnot--Carath\'eodory distance. In order to make our presentation more self-contained, however, we will just state the following result, which follows from the Euclidean coarea formula and will suffice for us.

\begin{prop}
  \label{prop:coarea}
  Let $T$ be a horizontal rectifiable current in $\Heis^n$ with finite mass. Then
  \begin{equation}
    \int_{0}^\infty \Mass(\langle T, \tau, t \rangle)\, \dd t \leq \Mass(T)
  \end{equation}
  where $\tau : \Heis \to \R_{\geq 0}$ is the Folland--Kor\'anyi gauge $\tau(\xi) = \FKnorm{\xi}$.
  \begin{proof}
    Since $\tau$ is smooth and Lipschitz outside of $\BHeis_{r}$ for any $r > 0$, it follows from \eqref{eq:hor-grad-tau-bound} and the Euclidean coarea formula that
  \[
    \int_{r}^\infty \Mass(\langle T, \tau, t \rangle)\, \dd t
    = \int_{\Heis^n \setminus \BHeis_r} |\nabla^{\vec{T}} \tau| \, \dd \|T\|
    \leq \int_{\Heis^n} \, \dd \|T\|
    = \Mass(T).
  \]
  The proposition follows by just letting $r \searrow 0$.
  \end{proof}
\end{prop}

We finally notice that horizontal currents indeed generalize isotropic submanifolds:
\begin{prop}
  Let $\Sigma^k \subset M^{2n+1}$ be a $C^1$ embedded isotropic submanifold with (possibly empty) $C^1$ boundary (recall that this means that $T_\xi\Sigma \subset \Xi_\xi$ for every $\xi \in \Sigma$). Then $\DBrack{\Sigma} \in \IC^\hor_k(M)$.
  \begin{proof}
    We have that for any $\alpha \in \mathcal{D}^{k-1}(M)$,
    \[
      \DBrack{\Sigma}(\theta \wedge \alpha) = \int_\Sigma \langle \theta \wedge \alpha, \vec{\tau}_\Sigma \rangle \, \dd \Haus^k = 0
    \]
    because $\vec{\tau}_\Sigma$ is a wedge product of vectors in $\Xi$. Since $T_\xi(\del \Sigma) \subset T_\xi \Sigma$, it follows that $\del \Sigma$ is isotropic too and thus for any $\beta \in \mathcal{D}^{k-2}(M)$,
    \[
      \DBrack{\Sigma}(\dd \theta \wedge \beta)
      = \DBrack{\Sigma}(\dd (\theta \wedge \beta) + \theta \wedge \dd \beta)
      = \del \DBrack{\Sigma}(\theta \wedge \beta) + \DBrack{\Sigma}(\theta \wedge \dd \beta)
      = \DBrack{\del \Sigma}(\theta \wedge \beta) = 0.
      \qedhere
    \]
  \end{proof}
\end{prop}

\begin{rmk}
  In light of this example, a natural question is whether any integer-rectifiable current $T$ which annihilates $\theta$ is automatically horizontal, that is, also annihilates $\dd \theta$. This is in general false, and a counterexample can be found in \cite{fu-erratum}. However, as shown in the same erratum, answer becomes positive if $T$ is an integral current.
\end{rmk}

\begin{thm}
  Let $T \in \IC_k(M)$ and suppose that for any $\alpha \in \mathcal{D}^{k-1}(M)$, $T(\theta \wedge \alpha) = 0$. Then $T \in \IC^\hor_k(M)$.
  \begin{proof}
    This was already proved by Fu in \cite{fu-erratum}, but we prefer to give the short proof here for the convenience of the reader and to fix some typos. Let $T$ be such a current and let $E := \spt (\del T)$. Then $E$ has locally finite $\Haus^{k-1}$ measure, hence $\Haus^{k}(E) = 0$ which implies that $\|T\|(E) = 0$.

    Let $\rho : M \to \R$ be the function $\dist(\cdot, E)$ with respect to any Riemannian metric. Clearly $\rho$ is Lipschitz, so the slices $\langle T, \rho, r \rangle$ are well-defined integral currents for almost every $r > 0$, and they are oriented by $\vec{T} \res \tfrac{\dd \rho}{|\dd \rho|}$. This implies that $\langle T, \rho, r \rangle(\theta \wedge \beta) = 0$ for each $\beta \in \mathcal{D}^{k-2}(M)$. We compute, for any such $\beta$ and almost every $r > 0$,
    \begin{align*}
      (T \res \{ \rho > r \} )(\dd \theta \wedge \beta)
      &= (T \res \{ \rho > r \} )(\dd (\theta \wedge \beta) + \theta \wedge \dd \beta) \\
      &= \del(T \res \{ \rho > r \} )(\theta \wedge \beta) + (T \res \{ \rho > r \})(\theta \wedge \dd \beta) \\
      &= ((\del T) \res \{ \rho > r \} )(\theta \wedge \beta)
      - \langle T, \rho, r \rangle (\theta \wedge \beta) = 0.
    \end{align*}
    Hence, using the fact that $\|T\|(E) = 0$ and $\|T\|$ is a Radon measure,
    \[
      |T(\dd \theta \wedge \beta)|
      \leq \limsup_{r \searrow 0} |(T \res \{ \rho \leq r \})(\dd \theta \wedge \beta)|
      \leq \sup_M (|\dd \theta| |\beta|) \limsup_{r \searrow 0} \|T\|\left(\{ \rho \leq r \} \cap \spt \beta \right)
      = 0.
      \qedhere
    \]
  \end{proof}
\end{thm}

\subsection{Existence for the Plateau problem in the Heisenberg group}

We first need a lemma to control the supports of a minimizing sequence.

\begin{lemma}
  \label{lem:control-support-plateau}
    Let $1 \leq k \leq n$, $r > 0$ and $m > 0$. Then for any current $T \in \IC^\hor_{k}(\Heis^n)$ with $\spt \del T \subset \BHeis_r$ and $\Mass(T) \leq m$ we can find another current $\hat{T} \in \IC^\hor_k(\Heis^n)$ with $\del \hat{T} = \del T$, $\Mass(\hat{T}) \leq \Mass(T)$ and $\spt \hat{T} \subset \BHeis_{R}$, where $R = r + C m^{\frac{1}{k}}$ for a constant $C = C(n, k) > 0$.
  \begin{proof}
    We argue by contradiction. Let $R_1 = r + \lambda m^{1/k}$, for a dimensional constant $\lambda > 0$ to be determined later, and consider the set 
    \[
      G = \left\{ t \in [r, R_1] : \langle T, \tau, t\rangle \text{ exists and } \Mass(\langle T, \tau, t\rangle) \leq \frac{2m}{R_1 - r} \right\},
    \]
    where $\tau$ is the Folland--Kor\'anyi gauge \eqref{eq:fk-gauge}. Then by \eqref{eq:hor-grad-tau-bound} and the coarea formula we have that
    \[
      \Leb^1([r, R_1] \setminus G)
      \leq \frac{R_1 - r}{2m} \int_{[r, R_1]} \Mass(\langle T, \tau, t\rangle) \, \dd t
      \leq \frac{R_1 - r}{2m} \Mass(T)
      \leq \frac{1}{2} (R_1 - r),
    \]
    so that $\Leb^1(G) \geq \tfrac{1}{2} (R_1 - r) = \tfrac{1}{2} \lambda m^{1/k}$. Now for every $t \in G$, consider the current $Q_t$ given by the isoperimetric inequality, which satisfies
    \[
      \del Q_t = \langle T, \tau, t\rangle,
      \qquad
      \Mass(Q_t) \leq C_\idxiso \Mass(\langle T, \tau, t\rangle)^{\frac{k}{k-1}}
      \qquad \text{and} \qquad
      \spt Q_t \subset \BHeis_{R},
    \]
    where
    \[
      R
      := R_1 + C\left(\frac{m}{R_1-r}\right)^{\frac{1}{k-1}}
      = r + \lambda m^{1/k} + C\lambda^{-\frac{1}{k-1}} m^{1/k}.
    \]
    This exists because $\del \langle T, \tau, t \rangle = -\langle \del T, \tau, t\rangle = 0$. Then let $\hat{T}_t := T \res \BHeis_t - Q_t$, also supported in $\BHeis_R$ and with boundary
    \[
      \del \hat{T}_t
      = \del (T \res \BHeis_t) - \del Q_t
      = \langle T, \tau, t \rangle + (\del T) \res \BHeis_t - \langle T, \tau, t \rangle
      = \del T.
    \]
    Now if the lemma were false we must have $\Mass(T) < \Mass(\hat{T}_t)$. Hence
    \[
      \Mass(T \res \BHeis_t) + \Mass(T \res \BHeis_t^c)
      = \Mass(T)
      \leq \Mass(\hat{T}_t)
      \leq \Mass(T \res \BHeis_t) + \Mass(Q_t)
      \leq \Mass(T \res \BHeis_t) + C_\idxiso \Mass(\langle T, \tau, t\rangle)^{\frac{k}{k-1}}.
    \]
    Letting $g(t) := \Mass(T \res \BHeis_t^c)$ and using the coarea formula with \eqref{eq:hor-grad-tau-bound}, for almost every $t \in G$ we obtain that
    \[
      g(t) = \Mass(T \res \BHeis_t^c)
      \leq C_\idxiso \Mass(\langle T, \tau, t\rangle)^{\frac{k}{k-1}}
      \leq C_\idxiso (-g'(t))^{\frac{k}{k-1}}.
    \]
    or
    \[
      \left(-g^{1/k}\right)' \geq C^{-1}.
    \]
    Since $g^{1/k}$ is monotonically decreasing, we can integrate over all of $[r, R_1]$ and find that
    \begin{align*}
      \frac{1}{2} \lambda m^{1/k}
      &\leq \Leb^1(G)
      \leq C\int_{[r, R_1]} \left(-g(t)^{1/k}\right)' \, \dd t
      \leq C(g(r)^{1/k} - g(R)^{1/k}) \\
      &\leq C g(r)^{1/k}
      \leq C \Mass(T)^{1/k}
      \leq Cm^{1/k},
    \end{align*}
    which is a contradiction if $\lambda$ is chosen large enough.
  \end{proof}
\end{lemma}

Using this lemma the existence part of \cref{thm:legendrian-plateau-problem} follows directly.

\begin{thm}
  \label{thm:existence-plateau}
  Let $1 \leq k \leq n$ and consider a current $S \in \IC^\hor_{k-1}(\Heis^n)$ with $\spt S$ compact. Then there exists a current $T \in \IC^\hor_{k}(\Heis^n)$ also with $\spt T$ compact such that $\Mass(T) \leq \Mass(T')$ for any other current $T' \in \IC^\hor_{k}(\Heis^n)$ with $\del T' = S$.
  \begin{proof}
    Let $(T_j)$ be a minimizing sequence, that is, $T_j \in \IC^\hor_k(\Heis^n)$, $\del T_j = S$ and
    \[
      \Mass(T_j) \longrightarrow \inf \{ \Mass(T') : T' \in \IC^\hor_k(\Heis^n), \del T' = S \}.
    \]
    Let $r > 0$ be such that $\spt S \subset \BHeis_r$. By \cref{lem:control-support-plateau}, since $\Mass(T_j)$ are bounded, we can improve our sequence to another minimizing sequence $( \hat{T}_j )$ with $\spt \hat{T}_j \subset \overline{\BHeis_{R}}$ for some $R > 0$. Then by the compactness theorem and the lower semicontinuity of the mass, a subsequence converges to a current $T$ belonging to the same class and supported in $\overline{\BHeis_{R}}$, thus attaining the minimum.
  \end{proof}
\end{thm}

\subsection{Existence of Legendrian minimizers in a homology class}

The existence part of \cref{thm:legendrian-homology-problem} follows easily using the compactness theorem once we can show that horizontal currents exist in any homology class. This is the content of the following proposition:

\begin{prop}
  Let $(M^{2n+1}, \Xi)$ be a closed contact manifold and $1 \leq k \leq n$. Then any homology class $\mathfrak{a} \in H_k(M, \Z)$ contains a horizontal integral cycle.
  \begin{proof}
    Let $g$ be any subriemannian metric on $\Xi$. Since $M$ is compact, we can cover it with finitely many Carnot--Carath\'eodory geodesic balls $B_R(\xi_1), \ldots B_R(\xi_N)$ such that $B_{2R}(\xi_i)$ admits a Darboux chart $\phi_i : B_{2R}(\xi_i) \to \mathcal{U}_i \subset \Heis^n$ for each $i$. Let $0 < \lambda \leq \Lambda < \infty$ be such that
    \[
      \lambda g_{\Heis} \leq (\phi_i^{-1})^* g \leq \Lambda g_{\Heis},
    \]
    so that in particular the maps $\phi_i$ preserve distances up to a constant.

    Given $\mathfrak{a} \in H_k(M, \Z)$, represent it by a Lipschitz polyhedral chain $\bm{\sigma} = \sum_j m_j \sigma_j$, where $m_j$ are positive integers and $\sigma_j$ are Lipschitz simplices each contained in some ball $B_r(\xi)$, where $r \leq R / C_\idxiterateiso$ and $C_\idxiterateiso$ is a constant to be determined below (depending only on $n$ and $\Lambda/\lambda$).

    We will deform $\bm{\sigma}$ into a horizontal integral current $S$ by inductively replacing its skeleton using the isoperimetric inequality (\cref{thm:isoperimetric-heis}). We first construct, for each $1$-simplex $\tau$ in the $1$-skeleton of $\bm{\sigma}$, a horizontal integral $1$-current $S_\tau$ with the same boundary $\del S_\tau = \del \DBrack{\tau}$, with mass at most $C^{(1)} r$, and with support still contained in a ball of radius $C^{(1)} r$ (for example we can take a minimizing Carnot--Carath\'eodory geodesic). Moreover we record the existence of a (not necessarily horizontal) filling, that is, an integral $2$-current $V_\tau$ with $\del V_\tau = S_\tau - \DBrack{\tau}$ and $\diam(\spt V_\tau) \leq C^{(1)} r$, which exists simply because these $1$-currents are supported in a contractible ball.

    Then, for each $2$-simplex $\tau$ in the $2$-skeleton of $\bm{\sigma}$, we construct a horizontal integral $2$-current $S_\tau$ whose boundary is $\del S_{\tau} = S_{\tau_0} - S_{\tau_1} + S_{\tau_2}$, where the $1$-simplices $\tau_0, \tau_1, \tau_2$ are defined by $\del \tau = \tau_0 - \tau_1 + \tau_2$. We construct this horizontal filling by applying \cref{thm:isoperimetric-heis} on $\Heis^n$ and pushing the currents back and forth by means of the Darboux charts, which preserve horizontality. Note that the resulting currents have $\diam(\spt S_\tau) \leq C^{(2)} r^2$ and $\Mass(S_\tau) \leq C^{(2)} r^2$, where $C^{(2)}$ depends only on $n$ and $\Lambda/\lambda$. Moreover, since
    \[
      \del S_{\tau} = S_{\tau_0} - S_{\tau_1} + S_{\tau_2} = \del (V_{\tau_0} - V_{\tau_1} + V_{\tau_2}) + \DBrack{\tau_0} - \DBrack{\tau_1} + \DBrack{\tau_2} = \del (V_{\tau_0} - V_{\tau_1} + V_{\tau_2} + \DBrack{\tau}),
    \]
    we have that $S_{\tau} = V_{\tau_0} - V_{\tau_1} + V_{\tau_2} + \DBrack{\tau} + \del V_\tau$ for some (not necessarily horizontal) integral $3$-current $V_\tau$ also with $\diam(\spt V_\tau) \leq C^{(2)} r$.

    This procedure can be iterated and in the $k$-th step, for each $k$-simplex $\tau$ of $\bm{\sigma}$ with $\del \tau = \sum_{l=0}^k (-1)^l \tau_l$, we obtain a $k$-dimensional horizontal integral current $S_\tau$ and a $(k+1)$-dimensional filling $V_\tau$ with
    \[
      S_\tau = \DBrack{\tau} + \sum_{l=0}^k (-1)^l \, V_{\tau_l} + \del V_{\tau}
    \]
    because additional terms coming from the fillings of each $\tau_l$ cancel each other, since $\del \del \tau = 0$. We also need that balls of diameter at most $C^{(n)} r$ are contained in a Darboux chart; we can guarantee this by choosing above $C_\idxiterateiso = C^{(n)}$, which only depends on $\Lambda / \lambda$ and $n$.

    After $n$ steps we will obtain horiontal integral $n$-currents $S_{\sigma_j}$ for each $j$, together with integral $(n+1)$-currents $V_{\sigma_j}$, satisfying $S_{\sigma_j} = \DBrack{\sigma_j} + \sum_{l=0}^k (-1)^l \, V_{\tau_{jl}} + \del V_{\sigma_j}$, where $\del \sigma_j = \sum_{l=0}^k (-1)^l \tau_{jl}$. Since $\bm{\sigma} = \sum_j \sigma_j$ is a cycle, its boundary terms cancel each other and we have
    \[
      \sum_{j} S_{\sigma_j}
      = \sum_{j} \DBrack{\sigma_j} + \sum_{l=0}^k (-1)^l \, V_{\tau_{jl}} + \del V_{\sigma_j}
      = \DBrack{\bm{\sigma}} + \del \left( \sum_{j} V_{\sigma_j} \right),
    \]
    which means that the horizontal integral $n$-current $S := \sum_j S_{\sigma_j}$ is a cycle representing $\mathfrak{a}$.
  \end{proof}
\end{prop}

\begin{rmk}
  The fact that there is no additional homological condition for a class to be ``horizontal'' should be compared to the fact that the cohomology of the Rumin complex is isomorphic to the usual De Rham cohomology. The situation changes drastically in the symplectic setting.
\end{rmk}

\begin{thm}
  \label{thm:existence-homology}
  Let $(M^{2n+1}, \Xi, g)$ be a closed subriemannian contact manifold, $1 \leq k \leq n$ and $\mathfrak{a} \in H_k(M, \Z)$ a homology class. Then there exists a cycle $T \in \IC^\hor_k(M)$ representing $\mathfrak{a}$ that minimizes the mass among all such cycles.
  \begin{proof}
    Extend $g$ to a Riemannian metric $g_0$ on $M$ and embed the resulting Riemannian manifold isometrically into some Euclidean space $\R^L$. Let $T_j$ be a minimizing sequence in $\mathfrak{a}$; then, by \cref{prop:ff-compactness}, after extracting a subsequence, $T_j \to T$ in the flat distance for an integral current $T \in \IC^\hor_k(M)$. This means that we can write $T_j = T + R_j + \del S_j$ with $\Mass(R_j), \Mass(S_j) \to 0$. As in \cite[Section~5.7]{morgan-gmt}, since $\del T_j = \del T = 0$, we have that $\del R_j = 0$, so by the Euclidean isoperimetric inequality, $R_j = \del Y_j$ for an integral $k$-current $Y_j$ in $\R^L$ with $\spt Y_j$ uniformly close to $M$ for $j$ large. Hence we may retract $Y_j$ onto $\tilde{Y}_j \in \IC_k(M)$ and have $T_j = T + \del (\tilde{Y}_j + S_j)$. In particular, for $j$ large, $T_j$ is homologous to $T$ in $M$ and the theorem follows. 
  \end{proof}
\end{thm}

\subsection{General properties of local minimizers}

In this section we define local minimizers and prove a strong convergence theorem. Since these concepts are local, we work with a subriemannian metric $h$ defined on an open subset $\mathcal{U}$ of the Heisenberg group.

Denote by $\mathscr{L}_0^k$ the set of oriented isotropic $k$-planes in $\Xi_0$, and set $\mathscr{L}_0 := \mathscr{L}_0^n$, the so-called oriented Lagrangian Grassmannian. We can identify $\mathscr{L}_0^k$ with a subset of $\bigwedge^k \Xi_0$; note that any $\vec{\pi} \in \mathscr{L}_0^k$ has $|\vec{\pi}| = 1$. After identifying all horizontal spaces via left translations, any Lipschitz metric $h$ on $\Xi |_{\mathcal{U}}$ induces a norm on $\mathscr{L}_0^k$ for each $\xi \in \mathcal{U}$ in the usual way:
\[
  |\vec{\pi}|_{h_\xi} = \sqrt{h_\xi(\vec{\pi}, \vec{\pi})}.
\]
The mass with respect to the metric $h$ then takes the form
\begin{equation}
  \label{eq:def-funct}
  \Mass^h_\mathcal{U}(T) = \int_{\mathcal{U}} |\vec{T}(\xi)|_{h_\xi} \, \dd \|T\|(\xi).
\end{equation}

\begin{defn}
  \label{defn:minimizing}
  We will say that a current $T \in \IC^\hor_k(\Heis^n)$ is \emph{$\Mass^h$-minimizing in $\mathcal{U}$} if $\Mass^h_{\mathcal{U}}(T) \leq \Mass^h_{\mathcal{U}}(T')$ whenever $T' \in \IC^\hor_k(\Heis^n)$ has $\spt (T' - T) \Subset \mathcal{U}$ and $T - T'$ is the boundary of a integral $(k+1)$-current supported in $\mathcal{U}$.

Note that the current whose boundary must be $T - T'$ is not required to be horizontal; in particular, if $k = n$, it will never be so. We add this condition just to enforce that homology minimizers are mass-minimizing on the whole manifold. Of course, if $\mathcal{U}$ is contractible, the condition reduces to $\del (T - T') = 0$. 
\end{defn}

As in the classical setting, sequences of mass minimizing currents have improved convergence properties. In particular, we will need the following proposition in the regularity theory. The proof relies crucially on a theorem Wenger \cite{wenger-flat-convergence} about flat convergence of integral currents in metric spaces that we adapt in the Appendix.

\begin{prop}
  \label{prop:strong-convergence}
  Let $\mathcal{U} \subset \Heis^n$ be a bounded contractible open set, $1 \leq k \leq n$, and $T_j \in \IC^\hor_k(\Heis^n)$ a sequence of currents with $(\del T_j) \res \mathcal{U} = 0$ and $\sup_j \|T_j\|(\mathcal{U}) < \infty$ that converges weakly to a current $T$. Suppose that $g_j$ is a sequence of subriemannian metrics on $\Xi|_\mathcal{U}$ converging uniformly to a subriemannian metric $g$, and that $T_j$ is $\Mass^{g_j}$-minimizing in $\mathcal{U}$.

  Then $T$ is $\Mass^g$-minimizing in $\mathcal{U}$ and
  $\|T_j\| \res \mathcal{U} \weakto \|T\| \res \mathcal{U}$ as Radon measures in $\mathcal{U}$, i.e.~against $C^0_c(\mathcal{U})$ functions.
  \begin{proof}
    Choose a Lipschitz function $\rho : \mathcal{U} \to [0, 1]$ such that the level sets $\{ \rho > r \}$ for $r > 0$ are compactly contained in $\mathcal{U}$ and exhaust $\mathcal{U}$. It is a classical fact (see for example \cite[Proposition~8.3]{ambrosio-kirchheim}) that for almost all $r > 0$ there is a subsequence $T_{j'}$ such that all the slices $\langle T_{j'}, \rho, r \rangle$ exist, have uniformly bounded masses and
  \[
    \langle T_{j'}, \rho, r \rangle \weakto \langle T, \rho, r \rangle
    \qquad \text{in } \mathcal{D}_{k-1}.
  \]
  Moreover we may assume that $\spt \, \langle T_{j'}, \rho, r \rangle \subset \{ \rho = r \}$ and $\langle T_{j'}, \rho, r \rangle = -\del(T_{j'} \res \{ \rho > r \})$ (hence $\del \langle T_{j'}, \rho, r\rangle = 0$), that the same holds for $T$, and that $\|T\|(\{ \rho = r \}) = 0$. Now \cref{thm:flat-convergence} implies that there exist $S_{j'} \in \IC^\hor_k(\Heis^n)$ with $\del S_{j'} = \langle T_{j'}, \rho, r \rangle - \langle T, \rho, r \rangle$, with $\spt S_{j'} \subset \BHeis_{s_{j'}}(\{ \rho = r \})$ and with $\Mass(S_{j'}) \searrow 0$ and $s_{j'} \searrow 0$. In particular, for $j'$ large, $\spt S_{j'} \subset \BHeis_{s_{j'}}(\{ \rho = r \}) \subset \{ \rho > r / 2 \} \Subset \mathcal{U}$. Now we have that
  \[
    \del S_{j'}
    = \langle T_{j'}, \rho, r \rangle - \langle T, \rho, r \rangle
    = -\del(T_{j'} \res \{ \rho > r \}) + \del(T \res \{ \rho > r \}),
  \]
  so $(T - T_{j'}) \res \{ \rho > r \} - S_{j'}$ bounds in $\mathcal{U}$.

  Let $Z \in \IC^\hor_k(\Heis^n)$ be another horizontal cycle with $\spt Z \Subset \mathcal{U}$ that bounds a $(k+1)$-dimensional integral current supported in $\mathcal{U}$.
  Thus we may use $(T - T_{j'}) \res \{ \rho > r \} - S_{j'} + Z$ to test the minimality of $T_{j'}$:
  \[
    \Mass^{g_{j'}}_{\mathcal{U}}(T_{j'})
    \leq \Mass^{g_{j'}}_{\mathcal{U}}(T_{j'} + (T - T_{j'}) \res \{ \rho > r \} - S_{j'} + Z).
  \]
  It follows that
  \begin{align*}
    \Mass^{g_{j'}}_{\mathcal{U}}(T_{j'})
    &\leq \Mass^{g_{j'}}_{\mathcal{U}}(T_{j'} \res \{ \rho \leq r \} + T \res \{ \rho > r \} + Z - S_{j'}) \\
    &\leq \Mass^{g_{j'}}_{\mathcal{U}}(T_{j'} \res \{ \rho \leq r \})
    + \Mass^{g_{j'}}_{\mathcal{U}}((T + Z) \res \{ \rho > r \})
    + \Mass^{g_{j'}}_{\mathcal{U}}(S_{j'})
  \end{align*}
  and thus
  \begin{equation}
    \label{eq:strong-convergence-aux}
    \Mass^{g_{j'}}_{\mathcal{U}}(T_{j'} \res \{ \rho > r\})
    \leq \Mass^{g_{j'}}_{\mathcal{U}}((T + Z) \res \{ \rho > r \})
    + \Mass^{g_{j'}}_{\mathcal{U}}(S_{j'}).
  \end{equation}
  Letting $j' \to \infty$ and using the uniform convergence $g_{j'} \to g$, the lower semicontinuity of the mass, and the fact that $\Mass(T \res \{ \rho = r \}) = 0$, we obtain
  \[
    \Mass^{g}_{\mathcal{U}}(T \res \{ \rho > r\})
    \leq \Mass^{g}_{\mathcal{U}}((T + Z) \res \{ \rho > r \}),
  \]
  from which the $\Mass^g$-minimality of $T$ is clear. To prove the convergence of the associated measures, we first prove the strict convergence in the sense of Reshetnyak of the $g$-masses of $T_j$ on appropriate open sets. Choosing $Z = 0$ above, for almost every $r > 0$, \eqref{eq:strong-convergence-aux} gives
  \[
    \Mass^{g_{j'}}_{\mathcal{U}}(T_{j'} \res \{ \rho > r\})
    \leq \Mass^{g_{j'}}_{\mathcal{U}}(T \res \{ \rho > r \})
    + \Mass^{g_{j'}}_{\mathcal{U}}(S_{j'}),
  \]
  which together with the uniform convergence $g_{j'} \to g$ implies that for a sequence $\delta_{j'} \to 0$,
  \[
    \Mass^g_{\mathcal{U}}(T_{j'} \res \{ \rho > r\})
    \leq \Mass^{g}_{\mathcal{U}}(T \res \{ \rho > r \}) + \delta_{j'}.
  \]
  Taking the lim sup we obtain that
  \[
    \limsup_{j' \to \infty} \Mass^g_{\mathcal{U}}(T_{j'} \res \{ \rho > r\})
    \leq \Mass^{g}_{\mathcal{U}}(T \res \{ \rho > r \}),
  \]
  whereas the opposite inequality
  \[
    \liminf_{j' \to \infty} \Mass^g_{\mathcal{U}}(T_{j'} \res \{ \rho > r\})
    \geq \Mass^{g}_{\mathcal{U}}(T \res \{ \rho > r \})
  \]
  holds always thanks to the weak convergence $T_{j'} \weakto T$. Hence the vector-valued measures\footnote{Here we are identifying $\Heis^n$ with $\R^{2n+1}$ and treating the $k$-vectors $\vec{T}_j$ as vectors in the abstract vector space $\bigwedge^k \Xi_0$.} $T_{j'} = \vec{T}_{j'} \|T_{j'}\|_g$ converge strictly to $T = \vec{T} \|T\|_g$ in the open set $\{ \rho > r \}$. Hence, Reshetnyak's continuity theorem (see for example \cite[Theorem~10.3]{rindler-book}) shows that for any function $f \in C^0_c(\{\rho > r\})$,
  \[
    \int_{\{ \rho > r \}} f(\xi) \, \dd \|T_{j'}\|
    = \int_{\{ \rho > r \}} \frac{f(\xi)}{|\vec{T}_{j'}(\xi)|_{g_\xi}} \, \dd \|T_{j'}\|_g
    \longrightarrow \int_{\{ \rho > r \}} \frac{f(\xi)}{|\vec{T}(\xi)|_{g_\xi}} \, \dd \|T\|_g
    = \int_{\{ \rho > r \}} f(\xi) \, \dd \|T\|.
  \]
  Finally, given any $f \in C^0_c(\mathcal{U})$ we can always find some $r > 0$ such that $\spt f \Subset \{ \rho > r \}$, so we have shown that $\|T_{j'}\| \weakto \|T\|$ for a subsequence, and now a standard contradiction argument establishes it for the whole sequence.
  \end{proof}
\end{prop}

\section{Main regularity theorems}
\label{sec:main-regularity}

In this section we state the $\eps$-regularity theorem for Legendrian local area minimizers (\cref{thm:eps-reg}) and derive its main consequences for the partial regularity of area-minimizing Legendrian currents. The proof of \cref{thm:eps-reg} will be carried out in the following two sections. Note that from here onwards we consider only Legendrian currents, that is, $n$-dimensional horizontal currents.

\subsection{Geometric considerations}

Thanks to the homogeneity of the Heisenberg group, we can develop our regularity theory around the origin. Fix a Legendrian plane $\pi$ through $0 \in \Heis^n$, namely a horizontal plane $\pi \subset \Xi_0$ such that $\omega|_{\pi} = 0$. It will be important to keep in mind that there is a canonical smooth action of the unitary group $\mathsf{U}(n)$ on $\Heis^n$ via
\begin{equation}
  \label{eq:action-unitary}
  U \cdot (\Vec{z}, \varphi) = (U \Vec{z}, \varphi), \qquad U \in \mathsf{U}(n)
\end{equation}
that preserves all the structures of the Heisenberg group, that is: acts as group automorphisms, pulls back $\theta$ to itself, and is compatible with the metric $g_\Heis$, the almost complex structure $J$ and the symplectic form $\omega$ on $\Xi$. This group action acts transitively on the set of all Legendrian planes in $\Xi_0$. Therefore, up to applying such an automorphism, there is no loss of generality in supposing that our plane $\pi$ is the plane $\pi_0 := \Span \{\del_{x^1}, \ldots, \del_{x^n} \}$.

Let $\vec{\pi} \in \bigwedge^n \Xi_0$ be an $n$-vector orienting $\pi$. By a slight abuse of notation, we will denote also by $\vec{\pi}$ the left-invariant $n$-vector field that extends $\vec{\pi}$ at the origin. We thus have the explicit form
\begin{equation}
  \vec{\pi}_0 = \vec{\nabla}^H x^1 \wedge \cdots \wedge \vec{\nabla}^H x^n
\end{equation}
for the plane $\pi_0$ with its natural orientation.

There is a canonical projection map $\Pi : \Heis^n \to \R^{2n} \cong \Xi_0$ that ``forgets'' the last coordinate and is a group homomorphism: $\Pi(\Vec{z}, \varphi) = \Vec{z}$. Then for any Lagrangian plane $\pi \subset \Xi_0$ with respect to $\omega$ (we will call such planes \emph{Legendrian}) we can consider the projections $\mathbf{p}^\pi : \Heis^n \to \pi \cong \R^n$, defined as the composition of $\Pi$ and then the $g_\Heis$-orthogonal projection onto $\pi$. We will just write $\mathbf{p}$ when the plane $\pi$ is clear. We will also consider the projection $\mathbf{q} = \mathbf{q}^\pi$ onto the orthogonal of $\pi$ (after applying $\Pi$). For the plane $\pi_0$ this has the coordinate expression
\[
  \mathbf{p}(\Vec{x}, \Vec{y}, \varphi) = \Vec{x},
  \qquad
  \mathbf{q}(\Vec{x}, \Vec{y}, \varphi) = \Vec{y}.
\]

It is clear from the expression for the multiplication in $\Heis^n$ that $\Pi : \Heis^n \to \R^{2n}$ is a group homomorphism. Hence the same is true about $\mathbf{p}^\pi$ and $\mathbf{q}^\pi$ for any Legendrian plane $\pi$.
For $r > 0$ and $x_0 \in \R^n$, let $B_r^\pi(x_0)$ denote the Euclidean open ball $B_r^\pi := \{ x \in \pi : |x - x_0| < r \} \subset \pi$ and define the cylinder $\Cyl_r^\pi(x_0) := (\mathbf{p}^\pi)^{-1}(B_r^\pi(x_0))$. We will omit the superscript $\pi$ when the plane is clear (most of the time it will be $\pi_0$), and also set $B_r := B_r(0)$ and $\Cyl_r := \Cyl_r(0)$. We will also denote for short $\Cyl_r^\pi(\xi) := \Cyl_r^\pi(\mathbf{p}^\pi(\xi))$ for $\xi \in \Heis^n$.

As in \cref{sec:existence} it will be convenient to endow $\Heis^n$ with a left-invariant Riemannian metric $g_0$ such that $g_0|_{\Xi} = g_{\Heis^n}$; this is not absolutely necessary but it will allow us to recycle standard arguments from the Riemannian setting. Nevertheless none of our hypotheses or conclusions will involve this auxiliary metric. In particular, since all the currents that we consider are horizontal, their mass is independent of $g_0$ by \cref{prop:horizontal-mass}.

Observe that with respect to such a metric $g_0$, any of the projections $\mathbf{p}$ satisfies
\begin{equation}
  \label{eq:proj-lipschitz}
  |\nabla \mathbf{p}| \leq 1.
\end{equation}
Indeed, it is enough to check this at $0 \in \Heis^n$ thanks to the left-invariance of $g_0$ and to the fact that $\mathbf{p}$ is a group homomorphism. But at $\nabla \mathbf{p}(0)$ is just the projection onto the $x^i$-coordinates, which clearly has norm at most $1$. A similar argument gives the following: for any oriented horizontal $n$-plane $\vec{\pi}'$ at any point $\xi \in \Heis^n$, its push-forward by $\mathbf{p}$ is
\begin{equation}
  \label{eq:proj-jacobian}
  \mathbf{p}_\# (\vec{\pi}') = \langle \vec{\pi}', \vec{\pi}(\xi)\rangle \vec{\pi}.
\end{equation}
As a consequence we get the classical formula for the excess of mass:
\begin{prop}
  \label{prop:excess-formula}
  Let $T \in \IC^\hor_n(\Heis^n)$ and let $K \subset \R^n$ be a measurable set. Then
  \begin{equation}
    \label{eq:excess-formula}
    \Mass(T \res \mathbf{p}^{-1}(K)) - (\mathbf{p}_\# T \res K)(\vec{\pi}) = \frac{1}{2} \int_{\mathbf{p}^{-1}(K)} |\vec{T} - \vec{\pi}|^2 \, \dd \|T\|.
  \end{equation}
  \begin{proof}
    Notice that $|\vec{T} - \vec{\pi}|^2 = 2 - 2 \langle \vec{T}, \vec{\pi} \rangle$. Thus, using \eqref{eq:proj-jacobian},
    \begin{align*}
      \frac{1}{2} \int_{\mathbf{p}^{-1}(K)} |\vec{T} - \vec{\pi}|^2 \, \dd \|T\|
      &= \|T\|(\mathbf{p}^{-1}(K)) - \int_{\mathbf{p}^{-1}(B_r)} \langle \vec{T}, \vec{\pi} \rangle \, \dd \|T\| \\
      &= \|T\|(\mathbf{p}^{-1}(K)) - (\mathbf{p}_\# T \res K)(\vec{\pi}).
      \qedhere
    \end{align*}
  \end{proof}
\end{prop}

This leads us to a natural notion of (cylindrical) excess, the quantity that will lead the regularity theory. Suppose that $T$ is a horizontal integral current in $\Heis^n$ with $(\del T) \res \Cyl_r^\pi(x_0) = 0$. Since $(\del \mathbf{p}^\pi_\# T) \res B_r^\pi(x_0) = \mathbf{p}^\pi_\# ((\del T) \res \Cyl_r^\pi(x_0)) = 0$, by the constancy theorem we have that $(\mathbf{p}^\pi_\# T) \res B_r^\pi(x_0) = Q \DBrack{B_r^\pi(x_0)}$ for some integer $Q$. By orienting $\vec{\pi}$ appropriately we can and will suppose that $Q \geq 0$, thus $\|\mathbf{p}^\pi_\# T\|(B_r^\pi(x_0)) = Q \omega_n r^n$, where $\omega_n = \Leb^n(B_1)$. Thus we can define:
\begin{defn}
  The \emph{excess} of $T$ in the cylinder $\Cyl_r^\pi(x_0)$ is the quantity
  \begin{equation}
    \label{eq:excess-def}
    \Exc(T, \Cyl_r^\pi(x_0))
    := r^{-n} \left( \|T\|(\Cyl_r^\pi(x_0)) - Q \omega_n r^n \right)
    = \frac{1}{2} r^{-n} \int_{\Cyl_r^\pi(x_0)} |\vec{T} - \vec{\pi}|^2 \, \dd \|T\|.
  \end{equation}
\end{defn}

\subsection{Assumptions and the main regularity theorem}

Although our regularity result works with more general functionals (essentially those parametric elliptic functionals defined on $\mathscr{L}_0$ satisfying the ellipticity assumption from \cite{schoen-simon-anisotropic}), for notational convenience we will restrict ourselves to functionals coming from a subriemannian metric $h$, as in \eqref{eq:def-funct}. We record here the estimates of $h$ that we will use---all of them correspond to suitable ellipticity conditions for the associated integrand:
\begin{assumption}
  \label{ass:metric-h}
  For some $0 < \lambda \leq \Lambda < \infty$ and $0 < \mu < \infty$, the metric $h$ satisfies the following in $\mathcal{U}$:
  \begin{equation}
    \label{eq:ass-mass-comparable}
    \lambda \leq |\vec{\pi}|_h \leq \Lambda
    \qquad \text{for any } \vec{\pi} \in \mathscr{L}_0;
  \end{equation}
  \begin{equation}
    \label{eq:ass-ellipticity}
    \lambda^4 |\vec{\pi}|^2 |\vec{\varpi}|^2
    \leq h(\vec{\pi}, \vec{\pi}) h(\vec{\varpi}, \vec{\varpi}) - h(\vec{\pi}, \vec{\varpi})^2
    \qquad \text{for any } \vec{\pi}, \vec{\varpi} \in \bigwedge\nolimits^n(\Xi_0) \text{ such that } g_{\Heis}(\vec{\pi}, \vec{\varpi}) = 0;
  \end{equation}
  \begin{equation}
    \label{eq:ass-coercivity}
    \lambda |\vec{\pi} - \vec{\varpi}|^2 \leq |\vec{\pi}|_{h} - \frac{h(\vec{\pi}, \vec{\varpi})}{|\vec{\varpi}|_h}
    \qquad \text{for any } \vec{\pi}, \vec{\varpi} \in \mathscr{L}_0;
  \end{equation}
  \begin{equation}
    \label{eq:ass-lipschitz}
    |h_\xi(\vec{\pi}) - h_\zeta(\vec{\pi})| \leq \mu \FKdist(\xi, \zeta)
    \qquad \text{for any } \vec{\pi} \in \mathscr{L}_0 \text{ and any } \xi, \zeta \in \mathcal{U}.
  \end{equation}
\end{assumption}

\begin{prop}
  \cref{ass:metric-h} is satisfied for any subriemannian metric $h$ which is Lipschitz in a neighborhood of a relatively compact set $\mathcal{U}$.
  \begin{proof}
    \eqref{eq:ass-mass-comparable} and \eqref{eq:ass-lipschitz} are clear with some $0 < \lambda_1 \leq \Lambda < \infty$ and some $0 < \mu < \infty$. For \eqref{eq:ass-ellipticity}, note that we may assume that $|\vec{\pi}| = |\vec{\varpi}| = 1$. Then, by Cauchy--Schwarz, the right hand side of \eqref{eq:ass-ellipticity} is positive for any such $(\vec{\pi}, \vec{\varpi})$ and the inequality follows by compactness after possibly making $\lambda_1$ smaller.

    To show \eqref{eq:ass-coercivity}, first assume that $|\vec{\pi} - \vec{\varpi}| \leq \delta$ for some $0 < \delta \leq \tfrac{1}{2}$ and write $\vec{\varpi} = (1-t) \vec{\pi} + \vec{\sigma}$ with $\vec{\sigma} \perp_{g_\Heis} \vec{\pi}$. It is clear that $(1-t)^2 + |\vec{\sigma}|^2 = 1$ and $t^2 + |\vec{\sigma}|^2 \leq \tfrac{1}{4}$, hence $0 \leq t \leq |\vec{\sigma}|^2$. Consider the expression for $|\vec{\varpi}|_h$ and expand it to second order in $|\vec{\sigma}|$:
    \begin{align*}
      |\vec{\varpi}|_h
      &= \sqrt{(1-t)^2 |\vec{\pi}|_h^2 + |\vec{\sigma}|_h^2 + 2(1-t)h(\vec{\pi}, \vec{\sigma})} \\
      &= \sqrt{|\vec{\pi}|_h^2 - (2t-t^2) |\vec{\pi}|_h^2 + |\vec{\sigma}|_h^2 + 2(1-t)h(\vec{\pi}, \vec{\sigma})} \\
      &\geq |\vec{\pi}|_h \sqrt{1 + \frac{2h(\vec{\pi}, \vec{\sigma}) - 2t |\vec{\pi}|_h^2 + |\vec{\sigma}|_h^2}{|\vec{\pi}|_h^2} - C |\vec{\sigma}|^3} \\
      &\geq |\vec{\pi}|_h \left(1 - t + \frac{h(\vec{\pi}, \vec{\sigma})}{|\vec{\pi}|_h^2} + \frac{|\vec{\sigma}|_h^2}{2|\vec{\pi}|_h^2} - \frac{1}{2} \frac{h(\vec{\pi}, \vec{\sigma})^2}{|\vec{\pi}|_h^4} - C |\vec{\sigma}|^3 \right).
    \end{align*}
    Now use \eqref{eq:ass-ellipticity} to bound the last two terms from below:
    \begin{align*}
      |\vec{\pi}|_h |\vec{\varpi}|_h
      &\geq (1 - t) |\vec{\pi}|_h^2 + h(\vec{\pi}, \vec{\sigma}) + \frac{|\vec{\pi}|_h^2 |\vec{\sigma}|_h^2 - h(\vec{\pi}, \vec{\sigma})^2}{2|\vec{\pi}|_h^2} - C |\vec{\sigma}|^3 \\
      &\geq (1 - t) |\vec{\pi}|_h^2 + h(\vec{\pi}, \vec{\sigma}) + \frac{\lambda_1^4 |\vec{\sigma}|^2}{2|\vec{\pi}|_h^2} - C |\vec{\sigma}|^3 \\
      &= h(\vec{\pi}, \vec{\varpi}) + \frac{\lambda_1^4 |\vec{\sigma}|^2}{2|\vec{\pi}|_h^2} (1 - C |\vec{\sigma}|)
      \geq h(\vec{\pi}, \vec{\varpi}) + \frac{\lambda_1^4 |\vec{\sigma}|^2}{4|\vec{\pi}|_h^2} 
    \end{align*}
    provided that we choose $\delta$ small enough. In that case,
    \begin{align*}
      |\vec{\pi}|_h - \frac{h(\vec{\pi}, \vec{\varpi})}{|\vec{\varpi}|_h}
      \geq \frac{\lambda_1^4 |\vec{\sigma}|^2}{4|\vec{\pi}|_h^2|\vec{\varpi}|_h^2} 
      \geq \frac{\lambda_1^4}{4\Lambda^3} |\vec{\sigma}|^2
      \geq \frac{\lambda_1^4}{8\Lambda^3} (|\vec{\sigma}|^2 + t^2)
      = \frac{\lambda_1^4}{8\Lambda^3} |\vec{\varpi} - \vec{\pi}|^2.
    \end{align*}
    On the other hand, by compactness, \eqref{eq:ass-coercivity} holds for $|\vec{\pi} - \vec{\varpi}| \geq \delta$ with some $\lambda' > 0$, so we are done by choosing $\lambda = \min\left\{ \lambda_1, \lambda', \tfrac{\lambda_1^4}{8\Lambda^3} \right\}$.
  \end{proof}
\end{prop}

\begin{rmk}
  In the Riemannian setting, given a point $p \in M$, one can choose a good set of coordinates (for example normal coordinates centered at $p$) such that the error term in \eqref{eq:ass-lipschitz} with $\xi = 0$ is quadratic in $\dist(\zeta, 0)$. This improvement gives directly $C^{1,\alpha}$ regularity of minimizing currents for any $\alpha < 1$. In our case however, the coordinates are additionally required to respect the contact structure, so obtaining a quadratic error in \eqref{eq:ass-lipschitz} with such a constraint is impossible without imposing any compatibility condition\footnote{In a symplectic manifold $(W, \omega)$ with a Riemannian metric $g$, a necessary and sufficient condition for such coordinates to exist at $p \in M$ is that $\nabla^g \omega(p) = 0$, where $\nabla^g$ is the Levi-Civita connection. This condition is weaker than K\"ahler.}. This is the reason for the $C^{1,1/2}$ regularity in \cref{thm:eps-reg} below. Of course, this is not a serious issue since we can obtain higher regularity provided that the manifold is regular enough.
\end{rmk}

We will need the following set of hypotheses to state and prove regularity around $0 \in \Heis^n$:
\begin{assumption}
  \label{ass:scale-invt-hypot}
  The current $T \in \IC^\hor_n(\Heis^n)$ satisfies:
  \begin{equation}
    \label{eq:ass-f-minimizing}
    T \text{ is } \Mass^h \text{-minimizing in } \mathcal{U}
    \text{ (as in \cref{defn:minimizing}) and } \FKdist(\spt T \cap \Cyl_{R/2}^\pi, \Heis^n \setminus \mathcal{U}) \geq \frac{R}{2}
  \end{equation}
  \begin{equation}
    \label{eq:ass-no-bdry}
    \partial T \res \Cyl_R^\pi = 0
  \end{equation}
  \begin{equation}
    \label{eq:ass-q-cover}
    \mathbf{p}^\pi_\# (T \res \Cyl_R^\pi) = Q \DBrack{B_R^\pi}
  \end{equation}
  \begin{equation}
    \label{eq:ass-density-q}
    \Theta^n(\|T\|, \xi) \geq Q \text{ for } \|T\|\text{-a.e.~} \xi \in \Cyl_R^\pi
  \end{equation}
  \begin{equation}
    \label{eq:ass-excess}
    \Exc(T, \Cyl_R^\pi) \leq \eps.
  \end{equation}
\end{assumption}

Note that the density $\Theta^n(\|T\|, \xi)$ here can be computed either with respect to a distance coming from a Riemannian metric, or with respect to the Carnot--Carath\'eodory distance; they agree $\|T\|$-almost everywhere thanks to \cref{thm:rumin-metric-currents}.
Before stating the $\eps$-regularity theorem, we need to introduce an appropriate notion of graph over a Legendrian plane $\pi_0$.

If we require that a graph of the form $\{ (x, g(x), f(x) \}$ is Legendrian, the following condition must be satisfied:
\[
  \dd f = \frac{1}{2} (\Vec{x} \cdot \dd \Vec{g} - \Vec{g} \cdot \dd \Vec{x}).
\]
By considering instead $\tfrac{1}{2} \Vec{x} \cdot \Vec{g} - f$, this becomes
\[
  \dd \left(\frac{1}{2}(\Vec{x} \cdot \Vec{g}) - f\right) =  \Vec{g} \cdot \dd \Vec{x},
\]
hence the function $v(x) = \tfrac{1}{2} \Vec{x} \cdot \Vec{g} - f$ recovers both $f$ and $g$ by the expressions
\[
  f(x) = \frac{1}{2}(\Vec{x} \cdot \Dd v(x)) - v(x),
  \qquad
  g(x) = \Dd v(x).
\]
Note that, in the symplectic $\R^{2n}$, this corresponds to the fact that a Lagrangian graph is the graph of the gradient of a certain potential funtion; in the Legendrian case, in addition, we recover the potential from the last coordinate.

Motivated by this observation, for a $C^1$ function $v : \Omega \subset \pi_0 \to \R$, we define the map $\Phi^v : \Omega \to \mathbf{p}^{-1}(\Omega) \subset \Heis^n$ as
\begin{equation}
  \Phi^v(x) := \left(x, \, \Dd v(x), \, \frac{1}{2} \Vec{x} \cdot \Dd v(x) - v(x) \right)
\end{equation}
and note that if $v \in C^{1,1}$ then $\Phi^v_\# \DBrack{\Omega}$ is a Legendrian current, since
\[
  (\Phi^v)^* \theta
  = \dd \left(\frac{1}{2} \Vec{x} \cdot \Dd v - v \right) - \frac{1}{2}\left( \Vec{x} \cdot \dd \Dd v - \Dd v \cdot \dd \Vec{x} \right)
  = \Dd v \cdot \dd \Vec{x} - \dd v = 0.
\]

Graphs of this form constitute intrinsic low-dimensional graphs in the sense of \cite{regular-submanifolds-heis}. The mass of such a graph can be computed as follows: it is immediate to see that
\begin{equation}
  \label{eq:pushfwd-phi}
  \Dd \Phi^v(x) e_i = \nabla^H x^i + \sum_{j=1}^n \del_{ij} v(x) \nabla^H y^i,
\end{equation}
where $e_1, \ldots, e_n$ is the standard basis of $\R^n$ and $\{ \nabla^H x^i, \nabla^H y^i \}$ is a $g_\Heis$-orthonormal basis of $\Xi$ adapted to $\pi_0$. Thus the $h$-area of $\Phi^v$ over any measurable set $K \subset \R^n$ is
\begin{equation}
  \label{eq:area-legendrian-graph-metric}
  \Mass^h(\Phi^v_\# \DBrack{K})
  = \int_{K} \sqrt{\det\nolimits_{ij} (g^v)_{ij}} \, \dd x.
\end{equation}
where $((g^v)_{ij})$ are the coefficients of the pullback metric $g^v = (\Phi^v)^* h$. In particular, for the metric $g_\Heis$ we have
\begin{equation}
  \label{eq:area-legendrian-graph}
  \Mass(\Phi^v_\# \DBrack{K})
  = \int_{K} \sqrt{\det(\id + \Dd^2 v \cdot \Dd^2 v)} \, \dd x.
\end{equation}
We can now state our main regularity theorem.

\begin{thm}[$\eps$-regularity]
  \label{thm:eps-reg}
  Let $\mathcal{U} \subset \Heis^n$ be an open set, $h$ a subriemannian metric on $\mathcal{U}$ satisfying \cref{ass:metric-h}, and $T \in \IC^\hor_n(\Heis^n)$ with $0 \in \spt T$. Then there exist constants $\epsreg, \bar{C} \in (0, 1)$, depending only on $n, \Lambda / \lambda$ and $Q$ such that, if $T$ satisfies \cref{ass:scale-invt-hypot} with $\eps \leq \epsreg$ and $\pi = \pi_0$, and moreover $\mu R \leq \epsreg$, then
  \begin{equation}
    \label{eq:eps-reg-graph}
    T \res \Cyl_{R/72} = Q \, \Phi^f_\# \DBrack{B_{R/72}}
  \end{equation}
  for a function $f \in C^{2,1/2}(B_{R/72}, \R)$ with
  \begin{equation}
    \label{eq:estimates-eps-reg}
    R^{-2} \sup_{B_{R/72}} |f| + R^{-1} \sup_{B_{R/72}} |\Dd f| + \sup_{B_{R/72}} |\Dd^2 f| + R^{1/2} \left[ \Dd^2 f \right]_{C^{1/2}(B_{R/72})}
    \leq \bar{C} \left( \Exc(T, \Cyl_R^{\pi_0}) + \mu R \right)^{1/2}.
  \end{equation}
  In particular, $T \res \Cyl_{R/72}$ is a $C^{1,1/2}$ Legendrian graph over $\pi_0$.
\end{thm}

The \hyperlink{proof-eps-reg}{Proof of \cref{thm:eps-reg}} will occupy Sections \ref{sec:height-bound-lipschitz} and \ref{sec:biharmonic-decay}. In the rest of this section we present the main consequences of \cref{thm:eps-reg}.

\subsection{Partial regularity} \cref{thm:eps-reg} motivates the following definition: we say that a point $\xi \in \spt T \setminus \spt \del T$ is regular if $\spt T$ is a $C^{1,1/2}$ submanifold in a neighborhood of $\xi$; otherwise we call it singular. This determines a partition $\spt T \setminus \spt \del T = \operatorname{Reg} T \cup \operatorname{Sing} T$.

For a Radon measure $\nu$ on $\Heis^n$, we define its $n$-dimensional density with respect to the Carnot--Carath\'eodory distance as
\begin{equation}
  \label{eq:cc-density}
  \Theta^n_{d_\mathrm{CC}}(\nu, \xi) := \lim_{r \searrow 0} \frac{\nu(\mathscr{B}_r(\xi))}{\omega_n r^n}
\end{equation}
whenever the limit exists. Here $\mathscr{B}_r(\xi) := \{ \zeta \in \Heis^n : d_\mathrm{CC}(\zeta, \xi) < r \}$ denotes a ball with respect to the Carnot--Carath\'eodory distance. \cref{thm:rumin-metric-currents}, in particular Equation \eqref{eq:density-cc-riem}, implies that for a horizontal integral current $T \in \IC^\hor_n(\Heis^n)$ and for $\|T\|$-almost every $\xi$, $\Theta^n_{d_\mathrm{CC}}(\|T\|, \xi)$ exists and agrees with the usual density $\Theta^n_{d_0}(\|T\|, \xi)$ computed with respect to a Riemannian distance $d_0$. As above, when its precise value does not matter on sets of $\|T\|$-measure zero, we will not specify with respect to which distance it is computed.

\begin{thm}
  \label{thm:main-reg}
  Let $\mathcal{U} \subset \Heis^n$ be an open set, $h$ a subriemannian metric on $\mathcal{U}$ satisfying \cref{ass:metric-h}, and $T \in \IC^\hor_n(\Heis^n)$ a $\Mass^h$-minimizing current in $\mathcal{U}$. Suppose that a point $\xi \in \mathcal{U}$ has a neighborhood $\mathcal{V} \subset \mathcal{U}$ such that $\spt \del T \cap \mathcal{V} = \varnothing$, $\Theta^n(\|T\|, \zeta) \geq Q$ for $\|T\|$-a.e.~$\zeta \in \mathcal{V}$ and
  \[
    (\delta_{1/\rho_j} \circ \ell_{\xi^{-1}})_\# T \weakto Q \DBrack{\vec{\pi}_0} \qquad \text{in } \mathcal{D}_n(\Heis^n) \text{ for a sequence } \rho_j \searrow 0,
  \]
  where $Q = \Theta^k_{d_\mathrm{CC}}(\|T\|, \xi)$ is a positive integer and $\vec{\pi}_0$ is a Legendrian plane. Then $\xi \in \operatorname{Reg}(T)$.
  \begin{proof}
    Let $\eta_{\xi, \rho_{j}}(\zeta) = \delta_{1/\rho_j}(\ell_{\xi^{-1}}(\zeta))$ and consider the currents $T_j := (\eta_{\xi, \rho_{j}})_\# T \weakto Q \DBrack{\vec{\pi}_0}$. Elementary scaling considerations show that these are local $\Mass^{h_j}$-minimizers on $\mathcal{V}_j := \eta_{\xi, \rho_j}(\mathcal{V})$ for the metrics $h_j := h \circ (\eta_{\xi, \rho_j}^{-1})$ (here we are viewing the metrics as $\Sym^2(\Xi_0^*)$-valued functions). Note that being an $\Mass^h$-minimizer is invariant under multiplication of $h$ by a constant.

    It is clear that the sets $\mathcal{V}_j$ exhaust $\Heis^n$, that the metrics $h_j$ converge to $h_\xi$ locally uniformly (since $h$ is in fact Lipschitz) and that, in their domains, the metrics $h_j$ satisfy \cref{ass:metric-h} with the same constants $\Lambda, \lambda$ as $h$ and with $\mu_j = \mu \rho_j \to 0$.

    We claim that for well chosen $1 < R < 2$, $6 < r < 7$ and $j$ large enough, the currents $\tilde{T}_j := T_j \res \BHeis_r$ satisfy the assumptions of \cref{thm:eps-reg} with $\mathcal{U} := \BHeis_5$. Since $T_j$ agrees with $\tilde{T}_j$ in a neighborhood of $0$, this will establish that $0$ is a regular point of $T_j$ and hence $\xi$ is a regular point of $T$.

    First of all, the condition $\mu_j R \leq \epsreg = \epsreg(n, \Lambda/\lambda, Q)$ is clear for $j$ large. Since $T_j$ is a sequence of minimizers with uniformly bounded masses in $\BHeis_{10}$ (thanks to the existence and finiteness of $\Theta^k_{d_\mathrm{CC}}(\|T\|, \xi)$), by \cref{prop:strong-convergence} we have that $\| T_j \| \res \BHeis_{10} \weakto Q \Haus^n \res \pi_0 \res \BHeis_{10}$ in the sense of measures. In particular, testing with the compact set $\{ |\Vec{y}| + |\varphi| \geq 1 \} \cap \overline{\BHeis_9}$ we have that
    \[
      \|T_j\|\left(\{ |\Vec{y}| + |\varphi| \geq 1 \} \cap \overline{\BHeis_9} \right) \xrightarrow{j \to \infty} 0.
    \]
    Now applying \cref{lem:lower-density-bound} below we deduce that, for $j$ large,
    \begin{equation}
      \label{eq:convergence-hausdorff}
      \spt T_j \cap \{ |\Vec{y}| + |\varphi| \geq 2 \} \cap \overline{\BHeis_8} = \varnothing,
    \end{equation}
    since any point $\zeta$ in the latter set must satisfy $c \leq \|T_j\|(\BHeis_{1/100}(\zeta)) \leq \|T_j\|(\{ |\Vec{y}| + |\varphi| \geq 1 \} \cap \overline{\BHeis_9})$ for a constant $c > 0$ independent of $j$.

    We are ready to check the Assumptions \ref{ass:scale-invt-hypot} for $\tilde{T}_j$: to show \eqref{eq:ass-f-minimizing}, it is clear that $\tilde{T}_j$ is a minimizer in $\BHeis_5$ since $\BHeis_5 \subset \BHeis_r$. Moreover, given a point $\zeta \in \spt \tilde{T}_j \cap \Cyl_{R/2}$, since $\zeta \in \spt T_j \cap \overline{\BHeis_8}$, by \eqref{eq:convergence-hausdorff} we must have that $|\Vec{y}| + |\varphi| < 2$, in particular $|\Vec{x}|, |\Vec{y}|, |\varphi| \leq 2$ and then a computation shows that $\FKnorm{\zeta} \leq 4 < 5 - \tfrac{R}{2}$, so the additional condition about the support in \eqref{eq:ass-f-minimizing} is satisfied too.

    For \eqref{eq:ass-no-bdry}, notice that $\spt (\del \tilde{T}_j) \subset \spt T_j \cap \del \BHeis_r \subset \spt T_j \cap \overline{\BHeis_8}$. Hence any point $\zeta \in \spt (\del \tilde{T}_j) \cap \Cyl_R$ must have, on one hand $|\Vec{x}| < R < 2$, and on the other hand, again by \eqref{eq:convergence-hausdorff}, $|\Vec{y}|, |\varphi| < 2$, which again implies that $\FKnorm{\zeta} \leq 4 < r$. Thus such a point cannot exist and this proves \eqref{eq:ass-no-bdry}.

    So far we have that $(\del \tilde{T}_j) \res \Cyl_R = 0$; this implies that $\del (\mathbf{p}_\# \tilde{T}_j) \res B_R = \mathbf{p}_\# (\del \tilde{T}_j \res \Cyl_R) = 0$ and therefore, by the constancy theorem, we can write $\mathbf{p}_\# (\tilde{T}_j \res \Cyl_R) = (\mathbf{p}_\# \tilde{T}_j) \res B_R = Q_j \DBrack{B_R}$ for a sequence of integers $Q_j$. On the other hand, if we choose an appropriate $r$, we have that
    \begin{equation}
      \label{eq:weak-convergence-cut-current}
    \tilde{T}_j = T_j \res \BHeis_r \weakto Q \DBrack{\vec{\pi}_0} \res \BHeis_r = Q \DBrack{B_r}
    \end{equation}
    and, if in addition we choose a suitable $R$, it also holds that $\tilde{T}_j \res \Cyl_R \weakto Q \DBrack{B_r} \res \Cyl_R = Q \DBrack{B_R}$. Then $Q_j \DBrack{B_R} = \mathbf{p}_\# (\tilde{T}_j \res \Cyl_R) \weakto Q \DBrack{B_R}$, which implies that eventually $Q_j = Q$ and \eqref{eq:ass-q-cover} holds.

    Trivially \eqref{eq:ass-density-q} is inherited from the assumption for $T$ on $\mathcal{V}$. Finally by \eqref{eq:weak-convergence-cut-current} and \cref{prop:strong-convergence} we have that $\Mass(\tilde{T}_j \res \Cyl_R) \to \Mass(Q \DBrack{B_r} \res \Cyl_R) = Q \omega_n R^n$ and therefore
    \[
      \Exc(\tilde{T}_j, \Cyl_R) = R^{-n} \left( \Mass(\tilde{T}_j, \Cyl_R) - Q \omega_n R^n \right) \leq \epsreg
    \]
    for $j$ large, which is \eqref{eq:ass-excess}.
  \end{proof}
\end{thm}

\begin{cor}
  \label{cor:regularity-open-dense}
  Let $\mathcal{U}, h$ and $T$ be as in \cref{thm:main-reg}. Then $\operatorname{Reg}(T)$ is dense in $\spt T \setminus \spt \del T \cap \mathcal{U}$.
  \begin{proof}
    Let $\mathcal{V} \subset \mathcal{U} \setminus \spt \del T$ be any open set with $\spt T \cap \mathcal{V} \neq \varnothing$, and let $Q \geq 1$ be the essential minimum of $\Theta^n(\|T\|, \cdot)$ on $\spt T \cap \mathcal{V}$, that is, the greatest positive integer $Q$ such that $\Theta^n(\|T\|, \xi) \geq Q$ for $\|T\|$-almost every $\xi$ in $\mathcal{V}$. It is clear that the set of $\xi \in \mathcal{V} \cap \spt T$ such that $\Theta^n(\|T\|, \xi) = Q$ has positive $\|T\|$-measure, and hence by \eqref{eq:density-cc-riem} and \cref{prop:blowup-heisenberg} in the Appendix we may find some such $\xi$ with
  \[
    (\delta_{1/\rho} \circ \ell_{\xi^{-1}})_\# T \weakto Q \DBrack{\vec{T}(\xi)} \qquad \text{in } \mathcal{D}_n(\Heis^n) \text{ as } \rho \searrow 0
  \]
  and $Q = \Theta^n_{d_\mathrm{CC}}(\|T\|, \xi)$. Thus the hypothesis of \cref{thm:main-reg} are satisfied and therefore $\xi \in \operatorname{Reg}(T)$. This shows that any nonempty relatively open subset of $\spt T \cap \mathcal{U} \setminus \spt \del T$ intersects $\operatorname{Reg} T$, i.e.~$\operatorname{Reg} T$ is dense in $\spt T \cap \mathcal{U} \setminus \spt \del T$.
  \end{proof}
\end{cor}

In particular, we get almost everywhere regularity under the additional (very strong) assumption of multiplicity one. This includes nevertheless all Legendrian currents which are graphs of gradients of $C^{1,1}$ functions, without any a priori bound on their $W^{2,\infty}$ norm, and even of $W^{2,n}$ functions. This answers partially a question from \cite{hessian-partial-regularity}. Our argument however assumes that the current is minimizing with respect to competitors that are not necessarily graphical, so we do not see how to use our results to show almost everywhere regularity of $W^{2,n}$ minimizing solutions of the Hamiltonian stationary equation.

\begin{cor}
  \label{cor:regularity-density-one}
  Let $\mathcal{U}, h$ and $T$ be as in \cref{thm:main-reg}. If in addition, $\Theta^n(\|T\|, \cdot) = 1$ $\|T\|$-almost everywhere, then $\Haus^n(\operatorname{Sing} T) = 0$.
  \begin{proof}
    This follows, as in \cref{cor:regularity-open-dense}, by applying \cref{thm:main-reg} in conjunction with \cref{prop:blowup-heisenberg} with $Q = 1$, but now we can do this at $\|T\|$-almost every point.
  \end{proof}
\end{cor}

\subsection{Higher regularity}

Higher regularity around points where $\spt T$ is $C^{1,1/2}$ follows from the work of \cite{bhattacharya-chen-warren}. Here we assume that the ambient metric is smooth, as they do; intermediate results for $C^{k,\alpha}$ metrics should follow by analyzing their proofs.

\begin{thm}
  \label{thm:regular-is-smooth}
  Suppose that $\mathcal{U}$, $T$ and $h$ are as in \cref{thm:eps-reg}, and moreover $h$ is smooth. Then for any regular point $\xi \in \operatorname{Reg}(T)$, $\spt T$ is in fact a smooth Legendrian submanifold in a neighborhood of $\xi$.
  \begin{proof}
    Let $\xi \in \operatorname{Reg}(T)$, so that $T = Q \Phi^f_\# \DBrack{B_r^{\pi_0}}$ for some $f : C^{2,1/2}(B_r, \R) $ on a Legendrian plane $\pi_0$. It is clear that for any test function $\psi \in C^\infty_c(B_r)$ and $t \in \R$ close enough to zero, the currents $Q \Phi^{f+t\psi}_\# \DBrack{B_r}$ are admissible competitors for the $\Mass^h$-minimality of $T$. Therefore by \eqref{eq:area-legendrian-graph-metric}
    \[
      \ddt \Mass^h\left(\Phi^{f+t\psi}_\# \DBrack{B_r}\right) = \ddt \int_{B_r} F(x, f, \Dd f, \Dd^2 f) \, \dd x = 0
    \]
    where we have set $F(x, f, \Dd f, \Dd^2 f) := \sqrt{\det\nolimits_{ij} \left(g^{f+t\psi})_{ij}\right)}\, \dd x = 0$ with $g^v := (\Phi^v)^* h$. This kind of functionals satisfy the requirements\footnote{
      The authors show this precisely in \cite[Section 3]{bhattacharya-chen-warren} assuming that $h$ does not depend on $f$ and that $h$ is compatible with the symplectic structure. However, their proof of higher regularity, based on difference quotients, is very robust and works as well in this setting. We remark that, in our main case of interest, when the contact manifold is a fibration over a sympletic manifold with a compatible metric, their computations apply directly.
    } of the higher regularity theory developed in \cite{bhattacharya-chen-warren} and in particular their Theorem 2.3 implies that $f$ is smooth.
  \end{proof}
\end{thm}

\begin{cor}
  \label{cor:real-analytic}
  If $T \in \IC^\hor_n(\Heis^n)$ is a Legendrian local mass minimizer for the standard metric, then $\operatorname{Reg}(T)$ is a real-analytic Legendrian submanifold.
  \begin{proof}
    This follows from the work of Morrey \cite{morrey-analytic} as in \cite{chen-warren-Cn}.
  \end{proof}
\end{cor}

\section{Height bound and Lipschitz approximation}
\label{sec:height-bound-lipschitz}

In this section we make the first steps towards the proof of \cref{thm:eps-reg} following \cite{schoen-simon-anisotropic}. We will work with a fixed Legendrian plane $\pi = \pi_0 = \Span\{ \del_{x^1}, \ldots, \del_{x^n} \}$. Here we are only assuming that $h$ satisfies \eqref{eq:ass-mass-comparable} from the assumptions.

\subsection{Area control and oscillation bounds}

The following lower bound on the mass of a minimizing current is well known in spaces which enjoy isoperimetric inequalities. However here we will need to adapt the proof to be able to control the support of the fillings.
\begin{lemma}[density lower bound]
  \label{lem:lower-density-bound}
  Let $T \in \IC^\hor_n(\Heis^n)$ be $\Mass^h$-minimizing in $\mathcal{U}$, and let $\xi \in \mathcal{U}$ and $r > 0$ be such that $\BHeis_{r}(\xi) \subset \mathcal{U}$ and $(\del T) \res \BHeis_r(\xi) = 0$. Then
  \begin{equation}
    \label{eq:lower-density-bound}
    \|T\|(\BHeis_r(\xi)) \geq c r^n
  \end{equation}
  for a constant $c = c(n, \Lambda/\lambda) > 0$.
  \begin{proof}
    We may assume that $\xi = 0$ and let $\tau(\zeta) = \FKnorm{\zeta}$. As in the proof of \cref{lem:control-support-plateau}, let
    \[
      G := \left\{ 0 < s < \frac{r}{2} : \langle T, \tau, s \rangle \text{ exists and }
      \Mass(\langle T, \tau, s \rangle) \leq \frac{4m}{r} \right\}.
    \]
    where $m := \|T\|(\BHeis_{r/2})$. We may suppose that
    \begin{equation}
      \label{eq:aux-bound-m}
      2 C_\idxiso^{n-1} m < \left(\frac{r}{2}\right)^{n},
    \end{equation}
    where $C_\idxiso$ denotes the constant from the isoperimetric inequality (\cref{thm:isoperimetric-heis}), otherwise we are done. By the coarea formula (\cref{prop:coarea}), 
    \[
      \Leb^2([0, \tfrac{r}{2}] \setminus G)
  \leq \frac{r}{4m} \int_{0}^{\frac{r}{2}} \Mass(\langle T, \tau, s \rangle) \, \dd s
      \leq \frac{r}{4m} \|T\|(\BHeis_{r/2})
      = \frac{1}{4} r,
    \]
    hence $\Leb^1(G) \geq \tfrac{1}{4} r$. For each $s \in G$ we may apply \cref{thm:isoperimetric-heis} and obtain a current $Q_s$ with
    \[
      \del Q_s = \langle T, \tau, s \rangle, \qquad
      \Mass(Q_s) \leq C_1 \Mass(\langle T, \tau, s \rangle)^{\frac{n}{n-1}}
    \]
    and additionally $\spt Q_s \subset \BHeis_{r} \subset \mathcal{U}$ thanks to \eqref{eq:aux-bound-m}. Now for $T_s = Q_s + T \res \BHeis_s^c$ we have that 
    \[
      \del T_s
      = \langle T, \tau, s \rangle + \del(T \res \BHeis_s^c)
      = \del (T \res \BHeis_s) - (\del T) \res \BHeis_s + \del(T \res \BHeis_s^c),
      = \del T
    \]
    so $T_s$ is a valid competitor for $\Mass^h$ and we get
    \[
      \|T\|(\BHeis_s)
      \leq \lambda^{-1} \Mass^h(T \res \BHeis_s)
      \leq \lambda^{-1} \Mass^h(Q_s)
      \leq \frac{\Lambda}{\lambda} \Mass(Q_s)
      \leq \frac{\Lambda}{\lambda} C_1 \Mass(\langle T, \tau, s \rangle)^{\frac{n}{n-1}}.
    \]
    Letting $g(s) := \|T\|(\BHeis_s)$ and using again the coarea formula, we see that
    \[
      g(s)^{\frac{n-1}{n}}
      \leq C g'(s).
    \]
    Thus, since $g(s) > 0$ for every $s > 0$, we have that $g^{\frac{1}{n}}(s)' = \frac{1}{n} g'(s) g(s)^{\frac{1-n}{n}} \geq c$ for every $s \in G$. Thanks to the monotonicity of $g(s)^{1/n}$ we may integrate and get
    \[
      g^{\frac{1}{n}}(r)
      \geq g^{\frac{1}{n}}\left(\tfrac{r}{2}\right)
      \geq \int_{G} \left(g^{\frac{1}{n}}(s)\right)' \, \dd s
      \geq c \Leb^1(G)
      \geq c r.
      \qedhere
    \]
  \end{proof}
\end{lemma}

As a consequence we obtain a bound on the vertical oscillation of the current, as in \cite[Lemma 5.3.4]{federer-gmt} or in the Appendix of \cite{schoen-simon-anisotropic}. The bound for the $\mathbf{q}$ components is standard, but since the proof of the bound for $\varphi$ uses the same strategy, we have decided to present them together for completeness.
\begin{lemma}[height bounds]
  \label{lem:height-bound}
  There are constants $\epshb, C > 0$ depending only on $n, \Lambda/\lambda$ and $Q$ such that, if $T$ satisfies \cref{ass:scale-invt-hypot} at scale $R = 1$ with $\eps \leq \epshb$, then 
  \begin{equation}
    \label{eq:height-bound-y}
    \sup_{\xi_1, \xi_2 \in \spt T \cap \Cyl_{1/2}} |\mathbf{q}(\xi_1) - \mathbf{q}(\xi_2)| \leq C \Exc(T, \Cyl_1)^{\frac{1}{2n}}.
  \end{equation}
  Moreover, if $\mathbf{q}(\xi_0) = 0$ for some $\xi_0 \in \spt T \cap \Cyl_{1/2}$, then also
  \begin{equation}
    \label{eq:height-bound-phi}
    \sup_{\xi_1, \xi_2 \in \spt T \cap \Cyl_{1/4}} |\varphi(\xi_1) - \varphi(\xi_2)| \leq C \Exc(T, \Cyl_1)^{\frac{1}{n(n+1)}}.
  \end{equation}
  \begin{proof}
    Denote $E := \Exc(T, \Cyl_1)$ and choose $\epshb \leq \tfrac{1}{2} Q \omega_n$ so that
    \begin{equation}
      \label{eq:total-mass-bound}
      \|T\|(\Cyl_1)
      = \Mass(\mathbf{p}_\# (T \res \Cyl_1)) + E
      \leq \frac{3}{2} Q \omega_n.
    \end{equation}
    Let $f : \Heis^n \to \R$ be a smooth function and for $-\infty \leq t_1 \leq t_2 \leq +\infty$ define
    \[
      W(t_1, t_2) := \{ \xi \in \Cyl_1 : t_1 < f(\xi) < t_2 \}
    \]
    and $T_{t} := T \res W(t, \infty)$ for $t \in \R$. Let $t_0$ be a median of $f$, in the sense that
    \[
      \max \{\|T\|(W(-\infty, t_0)), \|T\|(W(t_0, \infty)) \} \leq \frac{1}{2} \|T\|(\Cyl_1).
    \]
    We proceed in three steps.
    \begin{steps}
      \item We show that whenever $t_1 > t_0$ is such that $\Mass(T_{t_1}) \geq 2E$,
        \begin{equation}
          \label{eq:partial-height-bound}
          (t_1 - t_0) \Mass(T_{t_1})^{\frac{n-1}{n}}
          \leq C \int_{W(t_0, t_1)} \left| \nabla^{\vec{T}} f \right| \, \dd \|T\|
        \end{equation}
        for a constant $C = C(n, \Lambda/\lambda, Q) > 0$. First observe that $\Mass(T_{t_0}) \leq \tfrac{3}{4} Q \omega_n$. By the standard coarea formula and \eqref{eq:proj-lipschitz} we have
        \begin{align*}
          \Leb^n \{x \in B_1 : \Theta^n(\|\mathbf{p}_\# T_t\|, x) > 0 \}
          &\leq \Haus^n \{\xi \in \Cyl_1 : \Theta^n(\|T_t\|, \xi) > 0 \} \\
          &= \Haus^n \{\xi \in \Cyl_1 : \Theta^n(\|T_t\|, \xi) \geq Q \}
          \leq \frac{1}{Q} \Mass(T_t)
          \leq \frac{1}{Q} \Mass(T_{t_0})
          \leq \frac{3}{4} \omega_n.
        \end{align*}
        Notice that here is the only place where we use the assumption \eqref{eq:ass-density-q}. Thus, thanks to the fact that $\Leb^n \{x \in B_1 : \Theta^n(\|\mathbf{p}_\# T_t\|, x) = 0 \} \geq \tfrac{1}{4} \omega_n$, we may apply the Poincar\'e--Sobolev inequality (see \cite[5.3.2]{federer-gmt}) to the integer-valued BV functions corresponding to (the signed multiplicity of) $\mathbf{p}_\# T_t$:
        \[
          \|\mathbf{p}_\# T_t\|(B_1)^{\frac{n-1}{n}}
          \leq C \|\del \mathbf{p}_\# T_t\|(B_1)
          = C \|\mathbf{p}_\# \del T_t\|(B_1)
          = C \|\mathbf{p}_\# \del T_t\|(B_1).
        \]
        Now for $\Leb^1$-almost every $t > t_0$ we may apply \eqref{eq:proj-lipschitz} and use slicing to deduce
        \[
          \|\mathbf{p}_\# T_t\|(B_1)^{\frac{n-1}{n}}
          \leq C \|\mathbf{p}_\# \del T_t\|(B_1)
          \leq C \|\del T_t\|(\Cyl_1)
          = C \|\langle T, f, t \rangle\|(\Cyl_1).
        \]
        On the other hand, by \cref{prop:excess-formula}, 
        \[
          \|T_t\|(\Cyl_1)
          \leq \|\mathbf{p}_\# T_t\|(B_1) + \frac{1}{2} \int_{\Cyl_1} |\vec{T} - \vec{\pi}_0|^2 \, \dd \|T_t\|
          \leq \|\mathbf{p}_\# T_t\|(B_1) + E.
        \]
        Assuming also that $t < t_1$, we have that $\|T_t\|(\Cyl_1) \geq \|T_{t_1}\|(\Cyl_1) \geq 2E$, thus $E \leq \|\mathbf{p}_\# T_t\|(B_1)$ and we get
        \[
          \|T_t\|(\Cyl_1)^{\frac{n-1}{n}}
          \leq \left( 2 \|\mathbf{p}_\# T_t\|(B_1) \right)^{\frac{n-1}{n}}
          \leq C \|\langle T, f, t \rangle\|(\Cyl_1).
        \]
        Integrating gives finally the desired estimate:
        \[
          (t_1 - t_0) \|T_{t_1}\|(\Cyl_1)^{\frac{n-1}{n}}
          \leq C \int_{t_0}^{t_1} \|\langle T, f, t \rangle\|(\Cyl_1) \, \dd t
          = C \int_{W(t_0, t_1)} \left|\nabla^{\vec{T}} f(\xi)\right| \, \dd \|T\|(\xi).
        \]
      \item We now show the classical height bound \eqref{eq:height-bound-y}. Choose $f(\Vec{x}, \Vec{y}, \varphi) = y^k$ for a fixed $k \in \{1, \ldots, n\}$. We claim that there exists $\overline{t} \geq t_0$ such that $\|T\|(W(\overline{t}, \infty)) \leq \sqrt{E}$ and $\overline{t} - t_0 \leq CE^{\frac{1}{2n}}$. If $\|T\|(W(t_0, \infty)) \leq \sqrt{E}$ we may just choose $\overline{t} = t_0$, otherwise define $\overline{t}$ to be the supremum of all $t_1 > t_0$ such that $\|T\|(W(t_1, \infty)) \geq \sqrt{E}$. Then it is easy to see that $\|T\|(W(\overline{t}, \infty)) \leq \sqrt{E}$. To get the bound on $\overline{t} - t_0$, we apply \eqref{eq:partial-height-bound} for all such $t_1$ and pass to the limit $t_1 \nearrow \overline{t}$:
        \begin{equation}
          \label{eq:partial-height-bound-y}
          (\overline{t} - t_0) E^{\frac{n-1}{2n}}
          \leq C \int_{W(t_0, \overline{t})} \left| \nabla^{\vec{T}} y^k \right| \, \dd \|T\|
        \end{equation}
        (note that this is justified since $\epshb \leq \tfrac{1}{4}$ guarantees that $\sqrt{E} \geq 2E$). We now claim that
        \begin{equation}
          \label{eq:projected-gradient-ineq}
          \left| \nabla^{\vec{T}} y^k \right|^2 \leq 1 - |\langle \vec{T}, \vec{\pi}_0\rangle|^2.
        \end{equation}
        Recall that $\vec{\pi}_0 = \nabla^{H} x^1 \wedge \cdots \wedge \nabla^H x^n$ and $\{ \nabla^H x^i, \nabla^H y^i\}$ are a global orthonormal frame of $\Xi$. Thus
        \[
          1 =|\vec{T}|^2 \geq |\langle \vec{T}, \nabla^{H} x^1 \wedge \cdots \wedge \nabla^H x^n\rangle|^2 + \sum_{\mathbf{e}} \left| \left\langle \vec{T}, \nabla^H y^k \wedge \mathbf{e} \right\rangle \right|^2
        \]
        where the sum is over the orthonormal basis of $\bigwedge^{n-1}(\Xi)$ induced by our orthonormal frame. On the other hand,
        \[
          |\nabla^{\vec{T}} y^k|^2
          = |\vec{T} \res \dd y^k|^2
          = \sum_{\mathbf{e}} |\langle \vec{T} \res \dd y^k, \mathbf{e} \rangle|^2
          = \sum_{\mathbf{e}} |\langle \vec{T}, \nabla^H y^k \wedge \mathbf{e} \rangle|^2
        \]
        and \eqref{eq:projected-gradient-ineq} follows. Putting this together with the standard inequality $1 - \langle \vec{T}, \vec{\pi}_0 \rangle^2 \leq |\vec{T} - \vec{\pi}_0|^2$ we have
        \[
          \left| \nabla^{\vec{T}} y^k \right|^2 \leq |\vec{T} - \vec{\pi}_0|^2.
        \]
        Using this after applying Cauchy--Schwarz and the bound \eqref{eq:total-mass-bound} on \eqref{eq:partial-height-bound-y}, we get
        \[
          (\overline{t} - t_0) E^{\frac{1}{2} - \frac{1}{2n}}
          \leq C \left(\int_{W(t_0, \overline{t})} \left| \nabla^{\vec{T}} y^k \right|^2 \, \dd \|T\|\right)^{1/2}
          \leq C E^{\frac{1}{2}},
        \]
        from which the estimate $\overline{t} - t_0 \leq CE^{\frac{1}{2n}}$ follows. Finally let $\xi \in \spt T \cap \Cyl_{1/2}$, suppose that $y^k(\xi) > \overline{t}$ and let $r = \min \left\{\tfrac{1}{2}, y^k(\xi) - \overline{t} \right\}$. Since $\BHeis_r(\xi) \subset W(\overline{t}, \infty)$, the lower bound \eqref{eq:lower-density-bound} shows that
        \[
          c r^n \leq \|T\|(W(\overline{t}, \infty)) \leq \sqrt{E} \leq \sqrt{\epshb}.
        \]
        By choosing $\epshb$ small enough we can guarantee that $r < \tfrac{1}{2}$; in particular, $r = y^k(\xi) - \overline{t}$ and it follows that $y^k(\xi) \leq \overline{t} + CE^{\frac{1}{2n}} \leq t_0 + C E^{\frac{1}{2n}}$. Repeating the same argument for the level sets of $f$ below $t_0$ the oscillation bound for $y^k$ follows.
        
      \item We finally show the bound for $\varphi$. In order to use the bound from Step 2 we will have to apply Step 1 in the smaller cylinder $\Cyl_{1/2}$; we can do this thanks to the scale invariance of the hypotheses (see \cref{cor:rescaled-height-bound} below). In this case we choose $f(\Vec{x}, \Vec{y}, \varphi) := \varphi$ and, arguing as before, we will find $\overline{t} \geq t_0$ such that
        \[
          \overline{t} - t_0 \leq C E^{\frac{1}{n(n+1)}} \qquad \text{and} \qquad
          \|T\|(\Cyl_{1/2} \cap \{ \varphi > \overline{t} \}) \leq E^{\frac{1}{2(n+1)}}.
        \]
        To see this, we argue as above and define $\overline{t}$ as the supremum of those $t_1$ for which the second inequality does not hold. Then using \eqref{eq:partial-height-bound} (after making $\epshb$ small enough) and passing to the limit $t_1 \nearrow \overline{t}$ we get
        \[
          (\overline{t} - t_0) E^{\frac{n-1}{2n(n+1)}}
          \leq C \int_{\Cyl_{1/2} \cap \{ t_0 < \varphi < \overline{t} \}} \left| \nabla^{\vec{T}} \varphi \right| \, \dd \|T\|.
        \]
        Since $T$ is horizontal, $\|T\|$-almost everywhere it holds that
        \[
          \nabla^{\vec{T}} \varphi = \frac{1}{2} \left( \Vec{x} \cdot \nabla^{\vec{T}} \Vec{y} - \Vec{y} \cdot \nabla^{\vec{T}} \Vec{x} \right),
        \]
        hence $|\nabla^{\vec{T}} \varphi| \leq \frac{1}{2} \sum_k |x^k| |\nabla^{\vec{T}} y^k| + |y^k| |\nabla^{\vec{T}} x^k|$. Now \eqref{eq:height-bound-y} together with the additional hypothesis on $\spt T$ implies that $|y^k| \leq C E^{\frac{1}{2n}}$ for every point in $\spt T \cap \Cyl_{1/2}$, and moreover we have the trivial bounds $|x^k| \leq 1$ and $|\nabla^{\vec{T}} x^k| \leq 1$. Therefore
        \begin{align*}
          (\overline{t} - t_0) E^{\frac{n-1}{2n(n+1)}}
          &\leq C \int_{\Cyl_{1/2} \cap \{ t_0 < \varphi < \overline{t} \}} \left(E^{\frac{1}{2n}} + \sum_k |\nabla^{\vec{T}} y^k|\right) \, \dd \|T\| \\
          &\leq C E^{\frac{1}{2n}} + C \sum_k \left(\int_{\Cyl_{1/2}}  |\nabla^{\vec{T}} y^k|^2 \, \dd \|T\|\right)^{1/2} \\
          &\leq C E^{\frac{1}{2n}} + C E^{\frac{1}{2}} \\
          &\leq C E^{\frac{1}{2n}}
        \end{align*}
        and the bound $\overline{t} - t_0 \leq C E^{\frac{1}{n(n+1)}}$ follows. Finally consider a point $\xi_0 \in \spt T \cap \Cyl_{1/4} \cap \{ \varphi > \overline{t} \}$ with coordinates $(\Vec{x}_0, \Vec{y}_0, \varphi_0)$ and let $0 < r \leq \tfrac{1}{4}$. Given a point $\xi = (\Vec{x}, \Vec{y}, \varphi) \in \spt T \cap \BHeis_r(\xi_0)$, we have that
        \[
          \xi_0^{-1} \xi = (\Vec{x}', \Vec{y}', \varphi') = \left(\Vec{x} - \Vec{x}_0, \Vec{y} - \Vec{y}_0, \varphi - \varphi_0 + \frac{1}{2}(-\Vec{x}_0 \cdot \Vec{y} + \Vec{y}_0 \cdot \Vec{x}) \right) \in \BHeis_r(0),
        \]
        so it is clear that $|\Vec{x}| \leq \tfrac{1}{2}$ and, again thanks to \eqref{eq:height-bound-y},
        \[
          \varphi
          > \varphi_0 - \frac{1}{2}(-\Vec{x}_0 \cdot \Vec{y} + \Vec{y}_0 \cdot \Vec{x}) -\frac{r^2}{4}
          \geq \varphi_0 - \frac{1}{2}(|\Vec{x}_0| |\Vec{y}| + |\Vec{y}_0| |\Vec{x}|) -\frac{r^2}{4}
          \geq \varphi_0 - \overline{C} E^{\frac{1}{2n}} -\frac{r^2}{4}
          \geq \overline{t}
        \]
        whenever $0 < r \leq r_0$, where $r_0$ is the solution of $\tfrac{1}{4} r_0^2 + \overline{C} E^{\frac{1}{2n}} = \varphi_0 - \overline{t}$. Arguing as in the end of Step 2, let $r = \min\left\{\tfrac{1}{4}, r_0 \right\}$ and apply \eqref{eq:lower-density-bound} to $\BHeis_r(\xi_0)$, noting that $\BHeis_r(\xi_0) \cap \spt T \subset \Cyl_{1/2} \cap \{ \varphi > \overline{t} \}$:
        \[
          c r^n
          \leq \|T\|(\BHeis_r(\xi_0))
          = \|T\|(\BHeis_r(\xi_0) \cap \spt T)
          \leq \|T\|(\Cyl_{1/2} \cap \{ \varphi > \overline{t} \})
          \leq C E^{\frac{1}{2(n+1)}}.
        \]
        After making $\epshb$ smaller if necessary, it follows that $r < \tfrac{1}{4}$, so $r_0 = r \leq C E^{\frac{1}{2n(n+1)}}$ and finally
        \[
          \varphi_0 - t_0
          \leq (\varphi_0 - \overline{t}) + (\overline{t} - t_0)
          \leq \frac{1}{4} r_0^2 + \overline{C} E^{\frac{1}{2n}} + C E^{\frac{1}{n(n+1)}}
          \leq C E^{\frac{1}{n(n+1)}}
        \]
        The oscillation bound now follows by arguing analogously below $t_0$. \qedhere
    \end{steps}
  \end{proof}
\end{lemma}

By using the scaling automorphisms of $\Heis^n$ we get with little effort the following scaled version of the height estimates:

\begin{cor}
  \label{cor:rescaled-height-bound}
  For the same constants $\epshb, C > 0$ as in \cref{lem:height-bound} the following holds: suppose that $T$ satisfies \cref{ass:scale-invt-hypot} for some $R > 0$ and some $\eps \leq \epshb$, and that $\{\mathbf{q} = 0\} \cap \spt T \cap \Cyl_{R/2} \neq \varnothing$. Then
  \begin{equation}
    \label{eq:rescaled-height-bound-y}
    \sup_{\xi_1, \xi_2 \in \spt T \cap \Cyl_{R/2}} |\mathbf{q}(\xi_1) - \mathbf{q}(\xi_2)| \leq C \Exc(T, \Cyl_R)^{\frac{1}{2n}} R
  \end{equation}
  and
  \begin{equation}
    \label{eq:rescaled-height-bound-phi}
    \sup_{\xi_1, \xi_2 \in \spt T \cap \Cyl_{R/4}} |\varphi(\xi_1) - \varphi(\xi_2)| \leq C \Exc(T, \Cyl_R)^{\frac{1}{n(n+1)}} R^2.
  \end{equation}
  \begin{proof}
    Let $\delta_R : \Heis^n \to \Heis^n$ be the automorphism $\delta_R(\Vec{z}, \varphi) = (R \Vec{z}, R^2 \varphi)$. If $T$ satisfies \cref{ass:scale-invt-hypot} for some $R > 0$, then $\tilde{T} := (\delta_{R^{-1}})_\# T$ is also in $\IC^\hor_n(\Heis^n)$ and satisfies the same assumptions with $R = 1$, except possibly for \eqref{eq:ass-f-minimizing}. This is clear for \eqref{eq:ass-no-bdry}, \eqref{eq:ass-q-cover} and for \eqref{eq:ass-density-q} (since $\delta_R$ is a diffeomorphism). To get \eqref{eq:ass-excess} we use the definition
    \[
      \Exc((\delta_{R^{-1}})_\# T, \Cyl_1)
      = \|(\delta_{R^{-1}})_\# T\|(\Cyl_1) - Q \omega_n
      = R^{-n} \|T\|(\Cyl_R) - Q \omega_n
      = \Exc(T, \Cyl_R);
    \]
    here we have made use of the formula for the push-forward of currents,
    \[
      \|(\delta_{R^{-1}})_\# T\|(\Cyl_1)
      = \int_{\Cyl_R} \left|\left.{\wedge}^n (\delta_{R^{-1}})_\# \right|_{\xi} (\vec{T})\right| \, \dd \|T\|(\xi)
      = R^{-n} \int_{\Cyl_R} \, \dd \|T\|
      = R^{-n} \|T\|(\Cyl_R),
    \]
    since the tangent map $(\delta_{R^{-1}})_\#$ restricts to $\Xi_{\xi} \longrightarrow \Xi_{R^{-1} \xi}$ as multiplication by $R^{-1}$ (under the identification of all tangent planes to $\Heis^n$ via left translations).

    The scaled back version of \eqref{eq:ass-f-minimizing} asserts that $\tilde{T}$ is $\Mass^{\tilde{h}}$-minimizing in $\tilde{\mathcal{U}} := \delta_{R^{-1}}(\mathcal{U})$ for the metric $\tilde{h} := \delta_R^* h$. This follows from a similar computation:
    Since the quotient $\Lambda / \lambda$ is unaltered for $\tilde{h}$, we may apply the estimates \eqref{eq:height-bound-y} and \eqref{eq:height-bound-phi} to $\tilde{T}$. Then \eqref{eq:rescaled-height-bound-y} and \eqref{eq:rescaled-height-bound-phi} follow by simply scaling back.
  \end{proof}
\end{cor}

\subsection{Approximation by the graph of a Lipschitz gradient}

In order to embed most of the support of our minimizing current into such a graph, we will need an appropriate version of the Whitney--Glaeser $C^{1,1}$ extension theorem:
\begin{thm}
  \label{thm:whitney}
  Let $G \subset \overline{B_{R}} \subset \R^n$ be any set, $L_0, L_1, L_2 > 0$ given constants and $f : G \to \R$, $\Vec{g} : G \to \R^n$ given functions with $|f| \leq L_0$ and $|\Vec{g}| \leq L_1$ on $G$. Suppose that the following Whitney conditions are satisfied: for any $x, y \in G$,
  \begin{equation}
    |\Vec{g}(x) - \Vec{g}(y)| \leq L_2 |x - y|
    \quad \text{and} \quad
    |f(x) - f(y) - \Vec{g}(y) \cdot (\Vec{x} - \Vec{y})| \leq L_2 |x - y|^2.
  \end{equation}
  Then there exists a function $\tilde{f} \in C^{1,1}(\overline{B_R})$ such that $\tilde{f} = f$ and $\Dd \tilde{f} = \Vec{g}$ on $G$, and which satisfies the following bounds:
  \begin{equation}
    \label{eq:bounds-c11-ext}
    \| \tilde{f} \|_{L^\infty} \leq C L_0, \qquad
    \| \Dd \tilde{f} \|_{L^\infty} \leq C L_0 + C R L_2, \qquad
    \| \Dd^2 \tilde{f} \|_{L^\infty} \leq C L_2
  \end{equation}
  for a constant $C = C(n) > 0$.
  \begin{proof}
    This follows from the proof of Theorem VI.2.4 from \cite{stein-singular-integrals}. Notice that, in the definition of the extension through a sum on Whitney cubes, we do not need to restrict ourselves to cubes near $G$ since we are working on a bounded domain. Moreover, we may assume that $G$ is closed if we extend $f$ and $\Vec{g}$ by continuity. The more precise, scale-invariant estimates given here follow by closely examining the proof.
  \end{proof}
\end{thm}

Now the approximation of $T$ by a Lipschitz graph follows the strategy of \cite{schoen-simon-anisotropic}.

\begin{lemma}[Lipschitz approximation]
  \label{lem:lip-approx}
  Suppose that a current $T \in \IC^\hor_n(\Heis^n)$ satisfies \cref{ass:scale-invt-hypot} with $\eps \leq \epshb$ and has $0 \in \spt T$. Then for every $0 < \gamma \leq 1$ there exists a $C^{1,1}$ function $v : B_{R/9} \to \R$ with the following properties:
  \begin{equation}
    \label{eq:c11-bound-lip-approx}
    \sup_{B_{R/9}} | \Dd^2 v | \leq \gamma,
  \end{equation}
  \begin{equation}
    \label{eq:bounds-lip-approx}
    \sup_{B_{R/9}} |\Dd v| \leq C \left( \Exc(T, \Cyl_R)^{\frac{1}{2n}} + \gamma \right) R, \qquad
    \sup_{B_{R/9}} |v| \leq C \Exc(T, \Cyl_R)^{\frac{1}{n(n+1)}} R^2,
  \end{equation}
  \begin{equation}
    \label{eq:lip-bound-bad-set}
    T \res \mathbf{p}^{-1}(K) = T^v \res \mathbf{p}^{-1}(K)
    \text{ for a set } K \subset B_{R/9} \text{ with }
    \Leb^n (B_{R/9} \setminus K) \leq C \gamma^{-n(n+1)} R^n \Exc(T, \Cyl_R)
  \end{equation}
  and
  \begin{equation}
    \label{eq:lip-bound-total-mass-diff}
    \| T - T^v \|(\Cyl_{R/9}) \leq C \gamma^{-n(n+1)} R^n \Exc(T, \Cyl_R),
  \end{equation}
  where $T^v$ is the current $T^v := Q (\Phi^v)_\# \DBrack{B_{R/9}}$ and the constant $C$ depends only on $n$, $\Lambda/\lambda$ and $Q$.

  Moreover we have the following estimate for the Dirichlet energy of $v$:
  \begin{equation}
    \label{eq:energy-lip}
    \int_{B_{R/9}} |\Dd^2 v|^2 \, \dd x \leq C \gamma^{-n(n+1)} R^n \Exc(T, \Cyl_R).
  \end{equation}
  \begin{proof}
    Fix $0 < \eta < \epshb$ to be determined later and consider
    \begin{equation}
      \label{eq:good-set-lip}
      G_\gamma := \left\{ x \in B_{R/9} : \Exc(T, \Cyl_r(x)) \leq \eta \; \forall r \in (0, 8R/9) \right\}.
    \end{equation}
    For each $x \in B_{R/9} \setminus G_\gamma$ we can find $0 < r < 8R/9$ such that $\Exc(T, \Cyl_r(x)) > \eta$. Vitali's covering Lemma gives us a countable disjoint subcollection of such balls $B_{r_1}(x_1), B_{r_2}(x_2), \ldots \subset B_R$ such that $\{B_{5r_i}(x_i)\}$ cover $B_{R/9} \setminus G_\gamma$. Then
    \begin{equation}
      \label{eq:lip-covering-estimate}
      \Leb^n(B_{R/9} \setminus G_\gamma)
      \leq \sum_i \omega_n (5r_i)^n
      \leq \frac{\omega_n 5^n}{\eta} \sum_i \frac{1}{2} \int_{\Cyl_{r_i}(x_i)} |\vec{T} - \vec{\pi}_0|^2 \, \dd \|T\|
      \leq \frac{\omega_n 5^n}{\eta} R^n \Exc(T, \Cyl_{R}).
    \end{equation}
    On the other hand, let $\xi_1, \xi_2 \in \spt T \cap \mathbf{p}^{-1}(G_\gamma)$ and denote $\xi_i = (\Vec{x}_i, \Vec{y}_i, \varphi_i)$. Consider the current $\tilde{T} := (\ell_{\xi_1^{-1}})_\# T$, so that for any $r \in (4 |\Vec{x}_1 - \Vec{x}_2|, 8R/9)$ we have $0, \xi_1^{-1} \xi_2 \in \spt \tilde{T} \cap \Cyl_{r/4}$ and $\Exc(\tilde{T}, \Cyl_r) = \Exc(T, \Cyl_r(x_1)) \leq \eta$. The coordinates of the translated point are
    \[
      \xi_1^{-1} \xi_2 = \left(\Vec{x}_2 - \Vec{x}_1, \Vec{y}_2 - \Vec{y}_1, \varphi_2 - \varphi_1 + \frac{1}{2} (-\Vec{x}_1 \cdot \Vec{y}_2 + \Vec{y}_1 \cdot \Vec{x}_2) \right)
    \]
    and thus the height bound from \cref{cor:rescaled-height-bound} gives
    \[
      |\Vec{y}_2 - \Vec{y}_1| \leq C \eta^{\frac{1}{2n}} |\Vec{x}_2 - \Vec{x}_1|
      \qquad \text{and} \qquad
      \left|\varphi_2 - \varphi_1 + \frac{1}{2} (-\Vec{x}_1 \cdot \Vec{y}_2 + \Vec{y}_1 \cdot \Vec{x}_2)\right| \leq C \eta^{\frac{1}{n(n+1)}} |\Vec{x}_2 - \Vec{x}_1|^2
    \]
    after letting $r \searrow 4 |\Vec{x}_2 - \Vec{x}_1|$. In particular, if $\Vec{x}_1 = \Vec{x}_2$ it follows that $\Vec{y}_1 = \Vec{y}_2$ and $\varphi_1 = \varphi_2$, so $\spt T \cap \mathbf{p}^{-1}(x)$ consists of exactly one point\footnote{There is at least one point since otherwise $\FKdist(\spt T, \mathbf{p}^{-1}(x)) > 0$ and thus $\Mass(\mathbf{p}_\# T \res B_s(x)) = 0$ for some small $s > 0$, contradicting the assumption that $Q \geq 1$.} for every $x \in G_\gamma$. Hence we may define functions $f : G_\gamma \to \R, \Vec{g} : G_\gamma \to \R^n$ as follows: if $(\Vec{x}, \Vec{y}, \varphi)$ are the coordinates of the unique point in $\spt T \cap \mathbf{p}^{-1}(x)$, we set $f(x) := \tfrac{1}{2} (\Vec{x} \cdot \Vec{y}) - \varphi$ and $\Vec{g}(x) := \Vec{y}$.

    Let us show that we can apply \cref{thm:whitney} with $L_0 = C E^{\frac{1}{n(n+1)}} R^2$, $L_1 = C E^{\frac{1}{2n}} R$ and $L_2 = C\eta^{\frac{1}{n(n+1)}}$, where we are using $E := \Exc(T, \Cyl_R)$ for short. Since $0 \in \spt T$ and $E \leq \epshb$, the bounds $|f| \leq L_0$ and $|\Vec{g}| \leq L_1$ are clear. For the conditions involving $L_2$, just notice that
    \begin{align*}
      |f(x_1) - f(x_2) - \Vec{g}(x_2) \cdot (\Vec{x}_1 - \Vec{x}_2)|
      &= \left|\frac{1}{2}(\Vec{x}_1 \cdot \Vec{y}_1) - \varphi_1 - \frac{1}{2}(\Vec{x}_2 \cdot \Vec{y}_2) + \varphi_2 - \Vec{y}_2 \cdot (\Vec{x}_1 - \Vec{x}_2)\right| \\
       &= \left|\varphi_2 - \varphi_1 + \frac{1}{2} (-\Vec{x}_1 \cdot \Vec{y}_2 + \Vec{y}_1 \cdot \Vec{x}_2) \right| + \frac{1}{2} |(\Vec{y}_1 - \Vec{y}_2) \cdot (\Vec{x}_1 - \Vec{x}_2)| \\
       &\leq C \eta^{\frac{1}{n(n+1)}} |\Vec{x}_1 - \Vec{x}_2|^2 + C \eta^{\frac{1}{2n}} |\Vec{x}_1 - \Vec{x}_2|^2 \\
       &\leq C \eta^{\frac{1}{n(n+1)}} |\Vec{x}_1 - \Vec{x}_2|^2.
    \end{align*}
    whereas the condition for $\Vec{g}$ is clear. Thus we obtain a function $v \in C^{1,1}(B_R, \R)$ satisfying the bounds
  \[
    \| v \|_{L^\infty} \leq C E^{\frac{1}{n(n+1)}} R^2, \qquad
    \| \Dd v \|_{L^\infty} \leq C E^{\frac{1}{2n}} R + C \eta^{\frac{1}{n(n+1)}} R, \qquad
    \| \Dd^2 v \|_{L^\infty} \leq C \eta^{\frac{1}{n(n+1)}}
  \]
  and such that $\spt T \cap \mathbf{p}^{-1}(x) = \{\Phi^v(x)\}$ for every $x \in G_\gamma$. Choosing $\eta = c \gamma^{n(n+1)}$ for $c$ small enough proves \eqref{eq:bounds-c11-ext} and \eqref{eq:bounds-lip-approx}.

  Now note that for $\Leb^n$-a.e.~$x \in G_\gamma$, the slice $\langle T, \mathbf{p}, x \rangle$ is a $0$-dimensional integral current supported in $\mathbf{p}^{-1}(x) \cap \spt T = \{\Phi^v(x)\}$ and satisfies $\mathbf{p}_\# \langle T, \mathbf{p}, x \rangle = Q \DBrack{x}$, so we must have $\langle T, \mathbf{p}, x \rangle = Q \DBrack{\Phi^v(x)}$. Therefore, for a set $K \subset G_\gamma$ with $\Leb^n(G_\gamma \setminus K)$ we have that $T \res \mathbf{p}^{-1}(K) = T^v \res \mathbf{p}^{-1}(K)$, which together with \eqref{eq:lip-covering-estimate} and our choice of $\eta$ gives \eqref{eq:lip-bound-bad-set}. From this we get \eqref{eq:lip-bound-total-mass-diff} as a consequence of the formula for the area \eqref{eq:area-legendrian-graph} and \eqref{eq:excess-formula}:
  \begin{align*}
    \| T - T^v \|(B_{R/9})
    &\leq \| T - T^v \|(\mathbf{p}^{-1}(K)) + \|T\|(\mathbf{p}^{-1}(B_{R/9} \setminus K)) + \|T^v\|(\mathbf{p}^{-1}(B_{R/9} \setminus K)) \\
    &= \|T\|(\mathbf{p}^{-1}(B_{R/9} \setminus K)) + Q \int_{B_{R/9} \setminus K} \sqrt{\det(\id + (\Dd^2 v)^2)} \, \dd \Leb^n \\
    &\leq C \Leb^n(B_{R/9} \setminus K) + \frac{1}{2} \int_{\mathbf{p}^{-1}(B_{R/9} \setminus K)} |\vec{T} - \vec{\pi}_0|^2 \, \dd \|T\| \\
    &\leq C \Leb^n(B_{R/9} \setminus G_\gamma) + R^n \Exc(T, \Cyl_{R}) \\
    &\leq C \gamma^{-n(n+1)} R^n \Exc(T, \Cyl_R).
  \end{align*}

  Finally \eqref{eq:energy-lip} is a consequence of the following two estimates: on one hand,
  \[
    \int_{B_{R/9} \setminus K} |\Dd^2 v|^2
    \leq C \gamma^2 \cdot \gamma^{-n(n+1)} R^n \Exc(T, \Cyl_R)
    \leq C \gamma^{-n(n+1)} R^n \Exc(T, \Cyl_R),
  \]
  and on the other hand, again using \eqref{eq:area-legendrian-graph}, \eqref{eq:excess-formula} and the Taylor expansion of the area integrand
  \[
    \sqrt{\det(\id + A^2)} \geq 1 + c |A|^2 \qquad \text{for } |A| \leq 1,
  \]
  we find
  \begin{align*}
    \int_{K} |\Dd^2 v|^2
    &\leq C \int_{K} \left(\sqrt{\det(\id + (\Dd^2 v)^2)} - 1 \right)
    = \frac{C}{Q} \left( \| T^v \|(\mathbf{p}^{-1}(K)) - Q \Leb^n(K) \right) \\
    &= \frac{C}{Q} \left( \|T\|(\mathbf{p}^{-1}(K)) - Q \Leb^n(K) \right)
    \leq C R^n \Exc(T, \Cyl_R)
    \leq C \gamma^{-n(n+1)} R^n \Exc(T, \Cyl_R).
    \qedhere
  \end{align*}
\end{proof}
\end{lemma}

\begin{rmk}
  By using interpolation, we can obtain a better bound on the derivative of the $C^{1,1}$ function $v$ even when $\gamma$ is relatively large. Indeed, the estimate
  \begin{equation}
    \label{eq:interpolation-bound}
    \sup_{B_{R/9}} |\Dd v| \leq C \Exc(T, \Cyl_{R})^{\frac{1}{2n(n+1)}} R
  \end{equation}
  follows from interpolating between the bounds for $\|v\|_{L^\infty}$ and $\|\Dd^2 v\|_{L^\infty}$ from \cref{lem:lip-approx} (if necessary by considering a $C^{1,1}$ extension on a larger ball) and using that $\gamma \leq 1$.
\end{rmk}

\section{Biharmonic approximation and excess decay}
\label{sec:biharmonic-decay}

The main goal of this section is to prove the following lemma, which will be the main ingredient in the proof of \cref{thm:eps-reg}:

\begin{lemma}[excess improvement]
  \label{lem:excess-improvement}
  There exists a constant $C_\idximp = C_\idximp(n, \Lambda / \lambda, Q)$ with the following property: for any $0 < \rho \leq \tfrac{1}{72}$ there is $\eps_\rho > 0$ such that if $T \in \IC^\hor_n(\Heis^n)$ satisfies \cref{ass:scale-invt-hypot} with $\pi_0$ in place of $\pi$ and $\eps \leq \eps_\rho$, and in addition $0 \in \spt T$ and $\mu R \leq \eps_\rho$, then there exists an oriented Legendrian plane $\vec{\pi}_1 \in \mathscr{L}_0$ such that
  \begin{equation}
    \label{eq:rotation-excess-improvement}
    |\vec{\pi}_1 - \vec{\pi}_0|^2 \leq C_\idximp \left(\Exc(T, \Cyl_R^{\pi_0}) + \mu R \right)
  \end{equation}
  and
  \begin{equation}
    \label{eq:excess-improvement}
    \Exc(T \res \Cyl_{R/2}^{\pi_0}, \Cyl_{\rho R}^{\pi_1}) \leq C_\idximp \left( \rho^2 \Exc(T, \Cyl_R^{\pi_0}) + \rho^{-n} R \mu \right).
  \end{equation}
\end{lemma}

The strategy to prove \cref{lem:excess-improvement} is inspired by \cite{schoen-simon-anisotropic} but actually uses some technical simplifications of \cite{regularity-currents-hilbert} that allow us to avoid regularization and instead use classical $L^p$ estimates for solutions to elliptic fourth order constant coefficient equations in a ball.
In particular, our choice of a good radius comes from \cite{regularity-currents-hilbert}, but since our setting is fundamentally anisotropic, their strategy does not work for us and instead we use two distinct Lipschitz approximations as in \cite{schoen-simon-anisotropic}.

The idea (when $h$ is the standard Heisenberg metric) is to approximate the current $T$ by the graph of $\Phi^u$, where $u$ is the biharmonic function whose boundary data matches that of a Lipschitz approximation $\Phi^v$ on a suitably chosen ball. As already mentioned, our method actually works for more general functionals than those induced by metrics, but for simplicity we will restrict to these.

Before starting will need some computations from multilinear algebra.
For any $C^{1,1}$ function $w$, the $n$-vector orienting $T^w$ is
\begin{equation}
  \label{eq:tangent-graph}
  \vec{T}^w
  = \frac{(\Phi^w)_\# \vec{\pi}_0}{\left|(\Phi^w)_\# \vec{\pi}_0\right|}
  = \frac{(\Phi^w)_\# \vec{\pi}_0}{\Jac(\Phi^w)},
\end{equation}
where we have introduced the Jacobian
\begin{equation}
  \Jac(\Phi^w) = \left|(\Phi^w)_\# \vec{\pi}_0\right|
\end{equation}
and
\begin{equation}
  \label{eq:expansion-pushfwd}
  \begin{split}
    (\Phi^w)_\# \vec{\pi}_0
  &= (\Phi^w)_\# (e_1 \wedge \cdots \wedge e_n)
  = ((\Phi^w)_\# e_1) \wedge \cdots \wedge ((\Phi^w)_\# e_n) \\
  &= \left(\nabla^H x^1 + \sum_{j=1}^n \del_{1j}w \nabla^H y^j \right) \wedge \cdots \wedge \left(\nabla^H x^n + \sum_{j=1}^n \del_{nj}w \nabla^H y^j \right) \\
  &= \vec{\pi}_0 + \sum_{i,j=1}^n \del_{ij} w \, \vec{\pi}_i^j + \Ord^\perp(|\Dd^2 w|^2 + |\Dd^2 w|^n).
  \end{split}
\end{equation}
Here we have introduced the $n$-vectors
\[
  \vec{\pi}_i^j = \nabla^H x^1 \wedge \cdots \wedge \nabla^H x^{i-1} \wedge \nabla^H y^j \wedge \nabla^H x^{i+1} \wedge \cdots \wedge \nabla^H x^n
\]
and the notation $\Ord^\perp(f)$ indicates a term of size $\Ord(f)$ that belongs to the $g_\Heis$-orthogonal of $\Span \{ \vec{\pi}_0, \vec{\pi}_i^j \} \subset \bigwedge^n (\Xi)$.
It follows that
\[
  1
  \leq \Jac(\Phi^w)^2
  \leq 1 + C |\Dd^2 w|^2 + C |\Dd^2 w|^{2n}
\]
and hence
\begin{equation}
  \label{eq:estimate-jacobian}
  1
  \leq \Jac(\Phi^w)
  \leq 1 + C (|\Dd^2 w|^2 + |\Dd^2 w|^n).
\end{equation}
Therefore, if $|\Dd^2 w| \leq 1$, by expanding in \eqref{eq:tangent-graph} we find that
\begin{equation}
  \label{eq:tilt-gradient-comparable}
  | \vec{T}^w - \vec{\pi}_0 |^2 = |\Dd^2 w|^2 \left(1 + \Ord(|\Dd^2 w|^2)\right).
\end{equation}
In general we have the estimate
\begin{equation}
  \label{eq:Taylor-inverse-volume-element}
  \frac{1}{|\Phi^w_\# \vec{\pi}_0|_{h_0}} = \frac{1}{|\vec{\pi}_0|_{h_0}} \left( 1 - \sum_{i,j} \frac{\del_{ij}(w)\, h_0(\vec{\pi}_i^j, \vec{\pi}_0)}{|\vec{\pi}_0|_{h_0}^2} + \Ord(|\Dd^2 w|^2) \right).
\end{equation}
Indeed, for $|\Dd^2 w| \leq 1$ this comes from the Taylor expansion $|\Phi^w_\# \vec{\pi}_0|_{h_0}^2 = |\vec{\pi}_0|_{h_0}^2 + 2 \sum_{ij} \del_{ij}w \, h_0(\vec{\pi}_0, \vec{\pi}_i^j) + \Ord(|\Dd^2 w|^2)$, whereas for $|\Dd^2 w| \geq 1$ this is trivial given the lower bound
\begin{equation}
  \label{eq:lower-bound-volume-element}
  |\Phi^w_\# \vec{\pi}_0|_{h_0} \geq \lambda |\Phi^w_\# \vec{\pi}_0| \geq \lambda.
\end{equation}
We also have the expression
\begin{equation}
  \label{eq:jacobian-as-det}
  \Jac(\Phi^w) = \sqrt{\det(\id + (\Dd^2 w)^2)},
\end{equation}
which follows from \eqref{eq:expansion-pushfwd} using the usual formula for the Jacobian of a linear map.

\begin{proof}[Proof of \cref{lem:excess-improvement}]

  It is clear that we may suppose $R = 1$ and $E := \Exc(T, \Cyl_1) \leq \epshb$. We fix two parameters $0 < \delta \ll \eta \ll 1$ to be determined later but depending explicitly on $n$ only, and let $v_0$ and $v_\delta$ be the functions $v_0, v_\delta \in C^{1,1}(B_{1/9}, \R)$ provided by \cref{lem:lip-approx} with the choices $\gamma = 1$ and $\gamma = E^{\delta}$, respectively.  We record the estimates for the energy of $v_0$ and $v_\delta$ coming from \eqref{eq:energy-lip}:
\begin{equation}
  \label{eq:energy-v0}
  \int_{B_{1/9}} |\Dd^2 v_0|^2 \, \dd x \leq C E,
\end{equation}
\begin{equation}
  \label{eq:energy-vdelta}
  \int_{B_{1/9}} |\Dd^2 v_\delta|^2 \, \dd x \leq C E^{1-n(n+1)\delta}.
\end{equation}
It will be convenient to denote the $1$-Lipschitz approximation $L := T^{v_0} = Q (\Phi^{v_0})_\# \DBrack{B_{1/9}}$ and also, for any $r > 0$, $T_r := T \res \Cyl_r$ and $L_r := L \res \Cyl_r$.

\begin{steps}
\item \emph{The fourth order linear PDE.}
We define the coefficients
\[
  \tilde{a}_{ik}^{jl} := h_0(\vec{\pi}_i^j, \vec{\pi}_k^l),
  \qquad
  b_i^j := h_0(\vec{\pi}_0, \vec{\pi}_i^j)
  \qquad \text{and} \qquad
  a_{ik}^{jl} := \tilde{a}_{ik}^{jl} - \frac{b_i^j b_k^l}{|\vec{\pi}_0|_{h_0}^2},
\]
where the indices $i, j, k, l$ range from $1$ to $n$. It follows from \eqref{eq:ass-ellipticity} that $(a_{ik}^{jl})$ satisfies the Legendre ellipticity condition: given a symmetric matrix $(\sigma_{ij})$, we have that
\begin{align*}
  a_{ik}^{jl} \sigma_{ij} \sigma_{kl}
  &= \frac{1}{|\vec{\pi}_0|_{h_0}^2}(h_0(\vec{\pi}_i^j, \vec{\pi}_k^l) h_0(\vec{\pi}_0, \vec{\pi}_0) - h_0(\vec{\pi}_0, \vec{\pi}_i^j) h_0(\vec{\pi}_0, \vec{\pi}_k^l)) \sigma_{ij} \sigma_{kl} \\
  &= \frac{1}{|\vec{\pi}_0|_{h_0}^2}(h_0(\sigma_{ij} \vec{\pi}_i^j, \sigma_{kl} \vec{\pi}_k^l) h_0(\vec{\pi}_0, \vec{\pi}_0) - h_0(\vec{\pi}_0, \sigma_{ij} \vec{\pi}_i^j) h_0(\vec{\pi}_0, \sigma_{kl} \vec{\pi}_k^l)) \\
  &= \frac{1}{|\vec{\pi}_0|_{h_0}^2}(h_0(\vec{\sigma}, \vec{\sigma}) h_0(\vec{\pi}_0, \vec{\pi}_0) - h_0(\vec{\pi}_0, \vec{\sigma})^2)
  \geq \frac{\lambda^4}{\Lambda^2} |\vec{\sigma}|^2,
\end{align*}
where we have used $\vec{\sigma} = \sum_{ij} \sigma_{ij} \vec{\pi}_i^j \perp_{g_\Heis} \vec{\pi}_0$. We remark that these coefficients arise in the second order Taylor expansion of the function $(\sigma_{ij}) \mapsto \left|\vec{\pi}_0 + \sigma_{ij} \vec{\pi}_i^j \right|_{h_0}$.

Now for a radius $\tfrac{1}{18} < \sigma < \tfrac{1}{9}$ to be fixed soon, we consider the solution $u : B_\sigma \to \R$ of
\begin{equation}
  \label{eq:definition-u}
  \begin{cases}
    a_{ik}^{jl} \del_{ijkl} u = 0 & \text{in } B_\sigma \\
    u = v_\delta & \text{on } \del B_\sigma \\
    \del_\nu u = \del_\nu v_\delta & \text{on } \del B_\sigma
  \end{cases}
\end{equation}
provided by \cref{prop:biharmonic}, and recall that $u \in W^{2,p}(B_\sigma)$ for every $1 \leq p < \infty$ and $u \in C^\infty_\mathrm{loc}(B_\sigma)$.

  \item \emph{Definition of the comparison current.}
Let $S$ be the current $S := Q (\Phi^u)_\# \DBrack{B_{\sigma}}$. Some care is needed in defining $S$ up to the boundary and checking that $S \in \IC^\hor_n(\Heis^n)$, as $\Phi^u$ is not necessarily globally Lipschitz. This is a local issue and independent of the Heisenberg setting, so we may work in coordinates. Since $\Phi^u \in W^{1,n}(B_\sigma)$ by \eqref{eq:lp-estimate-biharmonic}, the current $(\Phi^u)_\# \DBrack{B_\sigma}$ has finite mass and is rectifiable (recall that $u$ is smooth inside $B_\sigma$). By the boundary rectificability theorem and the fact that $\Phi^{v_\delta}$ is Lipschitz, it is enough to show that $\del (\Phi^u)_\# \DBrack{B_\sigma} = (\Phi^{v_\delta})_\# \DBrack{\del B_\sigma}$. This identity follows from compensation properties of determinants analogous to those from \cite{brezis-nguyen-jacobian}. More precisely, given a smooth $(n-1)$-form $\alpha$,
\begin{align*}
  \del (\Phi^u)_\# \DBrack{B_\sigma}(\alpha)
  &= \int_{B_\sigma} (\Phi^u)^* (\dd \alpha)
  = \int_{B_\sigma} \dd ((\Phi^u)^* \alpha)
  = \sum_{|I| = |J| = n-1} \int_{B_\sigma} \dd \left( \det [\Dd (\Phi^u)]^I_J \, \alpha_I \circ \Phi^u \, \dd x^J \right) \\
  &= \sum_{|I| = n - 1} \sum_{k=1}^n (-1)^{k-1} \int_{B_\sigma} \det [\Dd (\Phi^u)]^I_{1 \ldots \hat{k} \ldots n} \, \der{x^k} \left( \alpha_I \circ \Phi^u \right) \, \dd x
\end{align*}
since the derivatives that act on the determinant cancel each other. Now observe that $(\Phi^u)^I$ is a $W^{1,n}$ extension of $(\Phi^{v_\delta})^I|_{\del B_\sigma}$, and $\alpha_I \circ \Phi^u$ is also a $W^{1,n}$ extension of $\alpha_I \circ \Phi^{v_\delta} |_{\del B_\sigma}$. The same argument as in \cite[Lemma 4 and Proposition 3]{brezis-nguyen-jacobian} but using H\"older's inequality with $L^{n}$ in all factors, instead of $L^{n-1}$ and $L^\infty$ (see for example \cite{lenzmann-schikorra} for this estimate in another context), shows that 
\[
  \sum_{k=1}^n \int_{B_\sigma} (-1)^{k-1} \det [\Dd (\Phi^u)]^I_{1 \ldots \hat{k} \ldots n} \, \der{x^k} \left( \alpha_I \circ \Phi^u \right)
  = \int_{\del B_\sigma} \det [\Dd (\Phi^{v_\delta})]^I \, (\alpha_I \circ \Phi^{v_\delta})
\]
for each multi-index $I$ of degree $n-1$, hence
\[
  \del (\Phi^u)_\# \DBrack{B_\sigma}(\alpha)
  = \sum_{|I| = n - 1} \int_{\del B_\sigma} \det [\Dd (\Phi^{v_\delta})]^I \, (\alpha_I \circ \Phi^{v_\delta})
  = \int_{\del B_\sigma} (\Phi^{v_\delta})^* \alpha
  = (\Phi^{v_\delta})_\# \DBrack{\del B_\sigma} (\alpha).
\]
See also \cite{regularity-currents-hilbert} for an alternative argument by approximation.

In the following, it will be convenient to extend $\vec{S}$ vertically to a $\bigwedge^n \Xi_0$-valued function on $\Cyl_{\sigma}$, that is, $\vec{S}(\xi) := \vec{S}(\Phi^u(\mathbf{p}(\xi)))$.
We also define the horizontal $n$-forms $\vec{S}^{h_0}$ and $\vec{\pi}_0^{h_0}$ by
\begin{equation}
  \vec{S}^{h_0} := \frac{h_0(\vec{S}, \cdot)}{|\vec{S}|_{h_0}}
  \qquad \text{and} \qquad
  \vec{\pi}_0^{h_0} := \frac{h_0(\vec{\pi}_0, \cdot)}{|\vec{\pi}_0|_{h_0}}.
\end{equation}

\item \emph{Choice of a good radius.}
  We claim that, as long as $\delta < \eta < \tfrac{1}{n(n+1)+n+2}$, we can choose $\sigma \in \left( \tfrac{1}{18}, \tfrac{1}{9} \right)$ that guarantees the following bounds for the current $S = Q \Phi^u_\# \DBrack{B_\sigma}$, if $u$ is the solution of \eqref{eq:definition-u} on $B_\sigma$:
  \begin{equation}
    \label{eq:boundary-S-estimate}
    \Mass(\del (T_\sigma - S)) \leq C E^{1-n(n+1)\delta},
  \end{equation}
  \begin{equation}
    \label{eq:good-radius-estimate}
    \left| (L_\sigma - T_\sigma)(\vec{S}^{h_0} - \vec{\pi}_0^{h_0}) \right| \leq C E^{1 + \eta},
  \end{equation}
  \begin{equation}
    \label{eq:energy-lip-near-boundary}
    \int_{B_\sigma \setminus B_{\sigma - E^\eta}} |\Dd^2 v_\delta|^2 + |\Dd^2 v_0|^2 \leq C E^{1 - n(n+1)\delta + \eta}.
  \end{equation}
  Here the constant $C$ should depend only on $n, \Lambda/\lambda$ and $Q$. To prove this, choose an integer $N$ between $\tfrac{E^{-\eta}}{72}$ and $\tfrac{E^{-\eta}}{36}$ and consider the sequence of radii $r_i = \tfrac{1}{18}\left(1 + \tfrac{i}{N} \right)$ for $i = 0, 1, \ldots, N$, so that $\tfrac{1}{18} = r_0 < r_1 < \cdots < r_N = \tfrac{1}{9}$. Notice that
    the measure
    \[
      \nu := E^{-1} \mathbf{p}_\# \| L - T \| + E^{-(1-n(n+1)\delta)} \mathbf{p}_\# \| T^{v_\delta} - T \| + E^{-(1-n(n+1)\delta)} \left( |\Dd^2 v_0|^2 + |\Dd^2 v_\delta|^2 \right) \, \dd x
    \]
    satisfies $\nu(B_{1/9}) \leq C$, hence for some $0 \leq i < N$ we have
    \[
      \| L - T \|(\Cyl_{r_{i+1}} \setminus \Cyl_{r_i})
      \leq \frac{C E}{N}
      \leq C E^{1 + \eta},
      \qquad
      \| T^{v_\delta} - T \|(\Cyl_{r_{i+1}} \setminus \Cyl_{r_i})
      \leq C E^{1 - n(n+1)\delta + \eta},
    \]
    and also 
    \[
      \int_{B_{r_{i+1}} \setminus B_{r_i}} |\Dd^2 v_\delta|^2 + |\Dd^2 v_0|^2 \, \dd x \leq C E^{1 - n(n+1)\delta + \eta}.
    \]
    Fix that $i$ and notice that, by slicing,
    we can find some $\tfrac{r_i + r_{i+1}}{2} \leq \sigma \leq r_{i+1}$ such that
    \[
      \Mass(\del(T^{v_\delta}_\sigma - T_\sigma))
      \leq \frac{2}{r_{i+1} - r_i} \| T^{v_\delta} - T \|(\Cyl_{r_{i+1}} \setminus \Cyl_{r_i})
      \leq C N \cdot \frac{1}{N} C E^{1 - n(n+1)\delta}
      \leq C E^{1 - n(n+1)\delta}.
    \]
    Now \eqref{eq:boundary-S-estimate} is clear since $\del T^{v_\delta}_\sigma = \del S$. Also \eqref{eq:energy-lip-near-boundary} follows since $\sigma - E^\eta \geq r_i$. On the other hand, we have the estimate
    \begin{align*}
      \sup_{B_{r_i}} |\Dd^2 u|^2
      &\leq \frac{C}{(\sigma - r_i)^n} \int_{B_\sigma} |\Dd^2 u|^2
      \leq C E^{-n\eta} \int_{B_\sigma} |\Dd^2 u|^2 \\
      &\leq C E^{-n\eta} \int_{B_\sigma} |\Dd^2 v_\delta|^2
      \leq C E^{-n\eta} E^{1-n(n+1)\delta}
      = C E^{1-n(n+1)\delta - n\eta}
    \end{align*}
    by \eqref{eq:lp-estimate-biharmonic} with $p = 2$ and \eqref{eq:energy-vdelta}. Now using the bound \eqref{eq:tilt-gradient-comparable} and the fact that $|\Dd^2 v_0| \leq 1$ almost everywhere, we get
    \begin{align*}
      \left|(L_\sigma - T_\sigma)(\vec{S}^{h_0} - \vec{\pi}_0^{h_0}) \right|
      &\leq C \| L - T \|(\Cyl_{\sigma} \setminus \Cyl_{r_i}) + C \| L - T \|(\Cyl_{r_i}) \cdot \sup_{B_{r_i}} |\vec{S} - \vec{\pi}_0| \\
      &\leq 2 \| L - T \|(\Cyl_{r_{i+1}} \setminus \Cyl_{r_i}) + C \| L - T \|(\Cyl_{1/9}) \cdot \sup_{B_{r_i}} |\Dd^2 u| \\
      &\leq C E^{1+\eta} + C E \cdot E^{\frac{1}{2}(1-n(n+1)\delta-n\eta)} \\
      &\leq C E^{1+\eta} + C E^{1 + \frac{1}{2}(1-n(n+1)\eta-n\eta)} \\
      &\leq C E^{1+\eta}
    \end{align*}
    thanks to the upper bound on $\eta$.

  \item \emph{Using the minimality of $T$.}
Once we have chosen a good radius, our goal is to show that the tangent planes to $S$ approximate those of $T$ very closely. More precisely, in this step we exploit the $\Mass^h$-minimality of $T$ to prove the following preliminary estimate:
  \begin{equation}
    \label{eq:T-S-filling-estimate}
    \frac{1}{2} \int_{\Cyl_\sigma} |\vec{T} - \vec{S}|^2 \, \dd \|T\|
    \leq C E^{1+\eta} + C \mu + C \left| (S - L_\sigma)(\vec{S}^{h_0} - \vec{\pi}_0^{h_0}) \right|
  \end{equation}
  Here the constant $C$ depends only on $n, \Lambda, \lambda$ and $Q$, and $\eta$ should be smaller than a dimensional constant.
Notice that the remaining term involves only the graphical currents corresponding to the functions $u$ and $v_0$. This will allow us to exploit the PDE satisfied by $u$ in the next step.

    To prove this, we begin with using the isoperimetric inequality (\cref{thm:isoperimetric-heis}) and \eqref{eq:boundary-S-estimate} to produce a current $V \in \IC^\hor_n(\Heis^n)$ with $\del V = \del (T_\sigma - S)$ and
    \[
      \Mass(V)
      \leq C \Mass(\del (T_\sigma - S))^{\frac{n}{n-1}}
      \leq C E^{\frac{n}{n-1}(1 - n(n+1) \delta)}
      \leq C E^{1+\eta}
    \]
    (the last inequality uses that $\eta$ is small enough). Moreover, by making $\eps$ smaller if necessary we may assume that $\spt V \subset \BHeis_{1/2} \subset \mathcal{U}$. Furthermore $\spt T_\sigma \subset \BHeis_{1/2}$, and thanks to Agmon's estimate \eqref{eq:agmon}, also $\spt S \subset \BHeis_{1/2}$. We may now fill $S + V - T_\sigma$ and use the minimality of $T$ to get
    \[
      \Mass^h(T_\sigma)
      \leq \Mass^h(S + V)
      \leq \Mass^h(S) + \Mass^h(V),
    \]
    which implies (since the total masses of $T$ and $S$ are bounded by a constant)
    \[
      \Mass^0(T_\sigma)
      \leq \Mass^h(T_\sigma) + C \mu
      \leq \Mass^h(S) + \Mass^h(V) + C \mu
      \leq \Mass^0(S) + C \Mass(V) + C \mu.
    \]
    Here we are denoting $\Mass^0 := \Mass^{h_0}$, with $h_0$ extended to $\Xi$ by left translations. Assumption \eqref{eq:ass-coercivity} now gives
    \begin{align*}
      \lambda \int_{\Cyl_\sigma} |\vec{T} - \vec{S}|^2 \, \dd \|T\|
      &\leq \int_{\Cyl_\sigma} \left(|\vec{T}|_{h_0} - \frac{h_0(\vec{T}, \vec{S})}{|\vec{S}|_{h_0}}\right) \, \dd \|T\|
      = \Mass^0(T_{\sigma}) - T_\sigma(\vec{S}^{h_0}) \\
      &= \Mass^0(T_{\sigma}) - \Mass^0(S) + S(\vec{S}^{h_0}) - T_\sigma(\vec{S}^{h_0}) \\
      &\leq C \Mass(V) + C \mu + (S - T_\sigma)(\vec{S}^{h_0}).
    \end{align*}
    Notice that $(S + V - T_\sigma)(\vec{\pi}_0^{h_0}) = 0$ since $\vec{\pi}_0^{h_0}$ is constant. Thus
    \begin{align*}
      \lambda \int_{\Cyl_\sigma} |\vec{T} - \vec{S}|^2 \, \dd \|T\|
      &\leq C \Mass(V) + C \mu + (S - T_\sigma)(\vec{S}^{h_0} - \vec{\pi}_0^{h_0}) - V(\vec{\pi}_0^{h_0}) \\
      &\leq C \Mass(V) + C \mu + (S - T_\sigma)(\vec{S}^{h_0} - \vec{\pi}_0^{h_0}).
    \end{align*}
    We finally estimate the last term by \eqref{eq:good-radius-estimate} and use the bound on $\Mass(V)$:
    \begin{align*}
      \int_{\Cyl_\sigma} |\vec{T} - \vec{S}|^2 \, \dd \|T\|
      &\leq C E^{1+\eta} + C \mu + C |(S - T_\sigma)(\vec{S}^{h_0} - \vec{\pi}_0^{h_0})| \\
      &\leq C E^{1+\eta} + C \mu + C |(S - L_\sigma)(\vec{S}^{h_0} - \vec{\pi}_0^{h_0})|.
    \end{align*}

  \item \emph{Linearization and $L^p$ estimates.}
    We proceed to estimate the last term in \eqref{eq:T-S-filling-estimate}. In the proof we will need to estimate simultaneously the energy of $u$, so we record it here as well. More precisely, we claim that if the $\eps$ in \eqref{eq:ass-excess} is small enough, then
  \begin{equation}
    \label{eq:estimate-deviation-T-S}
    \frac{1}{2} \int_{\Cyl_\sigma} |\vec{T} - \vec{S}|^2 \, \dd \|T\| \leq C E^{1+\delta/2} + C \mu
  \end{equation}
  and
  \begin{equation}
    \label{eq:energy-u}
    \int_{B_\sigma} |\Dd^2 u|^2 \, \dd x \leq C E + C \mu
  \end{equation}
  for a constant $C = C(n, \Lambda, \lambda, Q)$. To prove this, let us examine the last term in \eqref{eq:T-S-filling-estimate}:
    \begin{align*}
      (S - L_\sigma)(\vec{S}^{h_0} - \vec{\pi}_0^{h_0})
      &= Q \int_{B_\sigma} (\vec{S}^{h_0} - \vec{\pi}_0^{h_0})(\Phi^u_\# \vec{\pi}_0 - \Phi^{v_0}_\# \vec{\pi}_0) \\
      &= Q \int_{B_\sigma} h_0 \left(\frac{\vec{S}}{|\vec{S}|_{h_0}} - \frac{\vec{\pi}_0}{|\vec{\pi}_0|_{h_0}}, \Phi^u_\# \vec{\pi}_0 - \Phi^{v_0}_\# \vec{\pi}_0 \right) \\
      &= Q \int_{B_\sigma} h_0 \left(\frac{\Phi^u_\# \vec{\pi}_0}{|\Phi^u_\# \vec{\pi}_0|_{h_0}} - \frac{\vec{\pi}_0}{|\vec{\pi}_0|_{h_0}}, \del_{kl}(u - v_0) \vec{\pi}_k^l + \Ord(|\Dd^2 u|^2 + |\Dd^2 u|^n + |\Dd^2 v_0|^2) \right).
    \end{align*}
    We expand the first factor using \eqref{eq:Taylor-inverse-volume-element}:
    \begin{align*}
      \frac{\Phi^u_\# \vec{\pi}_0}{|\Phi^u_\# \vec{\pi}_0|_{h_0}} - \frac{\vec{\pi}_0}{|\vec{\pi}_0|_{h_0}}
      &= \frac{\vec{\pi}_0 + u_{ij} \vec{\pi}_i^j + \Ord(|\Dd^2 u|^2 + |\Dd^2 u|^n)}{|\vec{\pi}_0|_{h_0}}\left( 1 - \frac{u_{ij} b_i^j}{|\vec{\pi}_0|_{h_0}^2} + \Ord(|\Dd^2 u|^2) \right) - \frac{\vec{\pi}_0}{|\vec{\pi}_0|_{h_0}} \\
      &= \frac{1}{|\vec{\pi}_0|_{h_0}} \left( u_{ij} \vec{\pi}_i^j - u_{ij} b_i^j \frac{\vec{\pi}_0}{|\vec{\pi}_0|_{h_0}^2}\right) + \Ord(|\Dd^2 u|^2 + |\Dd^2 u|^{n+2}).
    \end{align*}
    Putting all together and analyizing the error terms (in particular we use that $2|\Dd^2 u| |\Dd^2 v| \leq |\Dd^2 u|^2 + |\Dd^2 v|^2$) we get:
    \begin{align*}
      |(S - L_\sigma)(\vec{S}^{h_0} - \vec{\pi}_0^{h_0})|
       & \leq \frac{Q}{|\vec{\pi}_0|_{h_0}} \left| \int_{B_\sigma} h_0 \left( u_{ij} \vec{\pi}_i^j - u_{ij} b_i^j \frac{\vec{\pi}_0}{|\vec{\pi}_0|_{h_0}^2}, \del_{kl}(u - v_0) \vec{\pi}_k^l \right) \, \dd x\right| & ( =: I)\\
       &\quad + C \int_{B_\sigma} |\Dd^2 u|^3 + |\Dd^2 u|^{2n+2} + |\Dd^2 u||\Dd^2 v_0|^2 \, \dd x. & (=: II)
    \end{align*}

    We first examine the term $(I)$. To begin with, we identify the coefficients $a_{ij}^{kl}$ in this expression:
    \begin{align*}
      \int_{B_\sigma} h_0 \left( u_{ij} \vec{\pi}_i^j - u_{ij} b_i^j \frac{\vec{\pi}_0}{|\vec{\pi}_0|_{h_0}^2}, \del_{kl}(u - v_0) \vec{\pi}_k^l \right) \, \dd x
      &= \int_{B_\sigma} \left( \tilde{a}_{ik}^{jl} - \frac{1}{|\vec{\pi}_0|_{h_0}^2} b_i^j b_k^l \right) u_{ij} \del_{kl}(u - v_0) \, \dd x \\
      &= \int_{B_\sigma} a_{ik}^{jl} \del_{ij} (u) \del_{kl}(u - v_0) \, \dd x.
    \end{align*}
    Integrating by parts twice on \eqref{eq:definition-u} against $u - v_\delta$ and using the fact that $u = v_\delta$ and $\Dd u = \Dd v_\delta$ on $\del B_\sigma$, we obtain that
    \[
      \int_{B_\sigma} a_{ik}^{jl} \del_{ij} (u) \del_{kl}(u - v_\delta) \, \dd x = 0.
    \]
    Therefore our task is to estimate
    \[
      (I) \leq C \left| \int_{B_\sigma} a_{ik}^{jl} \del_{ij} (u) \del_{kl}(v_\delta - v_0) \, \dd x \right|.
    \]
    Let $\chi \in C^\infty_c(B_\sigma)$ be a smooth radial cutoff function such that $\chi \equiv 1$ in $B_{\sigma - E^\eta}$, $\chi \equiv 0$ outside $B_{\sigma - E^\eta/2}$ and $|\Dd \chi| \leq C E^{-\eta}$. We have
    \begin{equation}
      \label{eq:aux-estimate-near-boundary}
      \begin{split}
        \left| \int_{B_\sigma} (1 - \chi) a_{ik}^{jl} \del_{ij} (u) \del_{kl}(v_\delta - v_0) \, \dd x \right|
    & \leq C \int_{B_\sigma \setminus B_{\sigma - E^\eta}} |\Dd^2 u| (|\Dd^2 v_\delta| + |\Dd^2 v_0|) \, \dd x \\
    & \leq C \int_{B_\sigma \setminus B_{\sigma - E^\eta}} E^{\eta/2} |\Dd^2 u|^2 + E^{-\eta/2} (|\Dd^2 v_\delta|^2 + |\Dd^2 v_0|^2) \, \dd x \\
    & \leq C E^{\eta/2} \int_{B_\sigma} |\Dd^2 u|^2 \, \dd x + C E^{1 - n(n+1)\delta + \eta / 2}
      \end{split}
    \end{equation}
    using Young's inequality and \eqref{eq:energy-lip-near-boundary}. For the remainder, we integrate by parts and use the interior estimates of \eqref{eq:interior-bound-hessian} and Cauchy--Schwarz:
    \begin{align*}
      \left| \int_{B_\sigma} \chi a_{ik}^{jl} \del_{ij} (u) \del_{kl}(v_\delta - v_0) \, \dd x \right|
      &= \left| -\int_{B_\sigma} \del_k \left(\chi \, a_{ik}^{jl} \del_{ij} (u)\right) \del_{l}(v_\delta - v_0) \, \dd x \right| \\
      &= \left| -\int_{B_\sigma} \left(\del_k \chi \, a_{ik}^{jl} \del_{ij} (u) + \chi \, a_{ik}^{jl} \del_{ijk} (u)\right) \del_{l}(v_\delta - v_0) \, \dd x \right| \\
      &\leq C \int_{B_{\sigma - E^\eta / 2}} \left(E^{-\eta} |\Dd^2 u| + |\Dd^3 u| \right) |\Dd v_\delta - \Dd v_0| \, \dd x \\
      &\leq C \sup_{B_{\sigma - E^\eta / 2}} \left(E^{-\eta} |\Dd^2 u| + |\Dd^3 u| \right) \left( \int_{B_\sigma} |\Dd v_\delta - \Dd v_0|^2 \, \dd x  \right)^{1/2} \\
      &\leq C E^{-\eta (n+2) / 2} \left(\int_{B_\sigma} |\Dd^2 u|^2 \, \dd x\right)^{1/2} \left(\int_{B_\sigma} |\Dd v_\delta - \Dd v_0|^2 \, \dd x\right)^{1/2}.
    \end{align*}
    Now recall that in a set $K = K_0 \cap K_\delta \cap B_\sigma$ such that $\Leb^n(B_\sigma \setminus K) \leq C E^{1-n(n+1)\delta}$ we have that $T^{v_0} \res \mathbf{p}^{-1}(K) = T^{v_\delta} \res \mathbf{p}^{-1}(K)$, so in particular $\Dd v_0 = \Dd v_\delta$ there, and everywhere we have the height bound $|\Dd v_\delta - \Dd v_0| \leq C E^{\frac{1}{2n(n+1)}}$ from \eqref{eq:interpolation-bound}. As a result,
    \[
      \int_{B_\sigma} |\Dd v_\delta - \Dd v_0|^2 \, \dd x
      \leq C E^{1 - n(n+1) \delta + \frac{1}{n(n+1)}}
    \]
    and
    \begin{equation}
      \label{eq:aux-estimate-chi}
      \left| \int_{B_\sigma} \chi a_{ik}^{jl} \del_{ij} (u) \del_{kl}(v_\delta - v_0) \, \dd x \right|
      \leq C \left(\int_{B_\sigma} |\Dd^2 u|^2 \, \dd x\right)^{1/2} \left(E^{1 - n(n+1) \delta + \frac{1}{n(n+1)} - (n+2) \eta} \right)^{1/2}.
    \end{equation}
    Now fix $\delta = \tfrac{\eta}{4n(n+1)}$ so that the error term in \eqref{eq:aux-estimate-near-boundary} is at most $C E^{1+\eta/4}$, and then choose $\eta$ small enough so that the exponent in the second factor of \eqref{eq:aux-estimate-chi} is at least $1 + \eta$. Note that we can do this thanks to the summand $\tfrac{1}{n(n+1)}$ which is independent of $\eta$. This will be our choice of parameters for the rest of the proof. Adding up \eqref{eq:aux-estimate-near-boundary} and \eqref{eq:aux-estimate-chi} and applying Young's inequality we get
    \begin{align*}
      (I)
      &\leq C \left| \int_{B_\sigma} \chi a_{ik}^{jl} \del_{ij} (u) \del_{kl}(v_\delta - v_0) \, \dd x \right| + C \left| \int_{B_\sigma} (1 - \chi) a_{ik}^{jl} \del_{ij} (u) \del_{kl}(v_\delta - v_0) \, \dd x \right| \\
      &\leq C E^{\eta/2} \left(\int_{B_\sigma} |\Dd^2 u|^2 \, \dd x\right)^{1/2} E^{1/2}
      + C E^{\eta/2} \int_{B_\sigma} |\Dd^2 u|^2 \, \dd x + C E^{1 + \eta / 4} \\
      &\leq C E^{\eta/2} \int_{B_\sigma} |\Dd^2 u|^2 \, \dd x + C E^{1 + \eta / 4}.
    \end{align*}
    Now we have to estimate the integrals that appear in $(I) + (II)$. Following \cite{regularity-currents-hilbert}, we split the integrals according to whether $|\Dd^2 u|$ is small or not. On one hand, choose $p = 1 + \tfrac{2}{\delta} = p(n) > 2n+2$ and compute
    \begin{equation}
      \label{eq:estimate-big-hessian}
      \begin{split}
      &\int_{ \{ |\Dd^2 u| > E^{\delta / 2} \} } |\Dd^2 u|^3 + |\Dd^2 u|^{2n+2} + |\Dd^2 u||\Dd^2 v_0|^2 + CE^{\eta/2} |\Dd^2 u|^2 \, \dd x \\
      &\qquad \leq C \int_{ \{ |\Dd^2 u| > E^{\delta / 2} \} } |\Dd^2 u| + |\Dd^2 u|^{2n+2} \, \dd x \\
      &\qquad \leq C \int_{B_\sigma} (E^{-(p-1)\delta/2} + E^{-(p-2n-2)\delta/2}) |\Dd^2 u|^p \, \dd x \\
      &\qquad \leq C E^{-(p-1)\delta/2} \int_{B_\sigma} |\Dd^2 u|^p \, \dd x
      \leq C_p E^{-(p-1)\delta/2} \int_{B_\sigma} |\Dd^2 v_\delta|^p \, \dd x \\
      &\qquad \leq C_p E^{p\delta-(p-1)\delta/2}
      \leq C E^{1+\delta}.
      \end{split}
    \end{equation}
    On the other hand, by \eqref{eq:energy-v0},
    \begin{equation}
      \label{eq:estimate-small-hessian}
      \begin{split}
      &\int_{ \{ |\Dd^2 u| \leq E^{\delta / 2} \} } |\Dd^2 u|^3 + |\Dd^2 u|^{2n+2} + |\Dd^2 u||\Dd^2 v_0|^2 + C E^{\eta/2} |\Dd^2 u|^2 \, \dd x \\
      & \qquad \leq C \int_{ \{ |\Dd^2 u| \leq E^{\delta / 2} \} } |\Dd^2 u|^3 + |\Dd^2 u||\Dd^2 v_0|^2 + E^{\eta / 2} |\Dd^2 u|^2 \, \dd x \\
      & \qquad \leq C E^{\delta / 2} \int_{\{ |\Dd^2 u| \leq E^{\delta / 2} \}} |\Dd^2 u|^2 + |\Dd^2 v_0|^2 \, \dd x 
      \leq C E^{\delta / 2} \int_{\{ |\Dd^2 u| \leq E^{\delta / 2} \}} |\Dd^2 u|^2 \, \dd x + C E^{1 + \delta / 2}.
      \end{split}
    \end{equation}
    Thanks to \eqref{eq:tilt-gradient-comparable}, by assuming that $\eps$ is small enough, we can estimate the first summand by
    \[
      \int_{\{ |\Dd^2 u| \leq E^{\delta/2} \}} |\Dd^2 u(x)|^2 \, \dd x
      \leq \int_{\{ |\Dd^2 u| \leq E^{\delta/2} \}} 2 |\vec{T}^u(x) - \vec{\pi}_0|^2 \, \dd x
      \leq 2 \int_{B_\sigma} |\vec{S}(x) - \vec{\pi}_0|^2 \, \dd x.
    \]
    Now recall that $\mathbf{p}_\# (T \res \Cyl_\sigma) = Q \DBrack{B_\sigma}$. Therefore the coarea formula and \eqref{eq:proj-lipschitz} give
    \[
      \int_{B_\sigma} |\vec{S}(x) - \vec{\pi}_0|^2 \, \dd x
      = \frac{1}{Q} \int_{B_\sigma} |\vec{S}(x) - \vec{\pi}_0|^2 \, \dd \|\mathbf{p}_\# T\|(x)
      \leq \frac{1}{Q} \int_{\Cyl_\sigma} |\vec{S}(\xi) - \vec{\pi}_0|^2 \, \dd \|T\|(\xi),
    \]
    since we defined $\vec{S}(\xi) = \vec{S}(\mathbf{p}(\xi))$. Furthermore,
    \[
      \int_{\Cyl_\sigma} |\vec{S} - \vec{\pi}_0|^2 \, \dd \|T\|
      \leq 2 \int_{\Cyl_\sigma} |\vec{S} - \vec{T}|^2 \, \dd \|T\| + 2 \int_{\Cyl_\sigma} |\vec{T} - \vec{\pi}_0|^2 \, \dd \|T\|
      \leq 2 \int_{\Cyl_\sigma} |\vec{S} - \vec{T}|^2 \, \dd \|T\| + 2 \Exc(T, \Cyl_1)
    \]
    and as a result,
    \begin{equation}
      \label{eq:energy-estimate-small-hessian}
      \int_{\{ |\Dd^2 u| \leq E^{\delta/2} \}} |\Dd^2 u|^2
      \leq C \int_{\Cyl_\sigma} |\vec{S} - \vec{T}|^2 \, \dd \|T\| + C E.
    \end{equation}
    Putting this together with \eqref{eq:estimate-small-hessian}, \eqref{eq:estimate-big-hessian} and \eqref{eq:T-S-filling-estimate} we get
    \[
      \frac{1}{2} \int_{\Cyl_\sigma} |\vec{T} - \vec{S}|^2 \, \dd \|T\|
      \leq
      C E^{\delta/2} \int_{\Cyl_\sigma} |\vec{T} - \vec{S}|^2 \, \dd \|T\|
      + C E^{1+\delta/2}
      + C \mu
    \]
    and, by choosing $\eps$ small enough, we can absorb the first term and \eqref{eq:estimate-deviation-T-S} is proven. Finally \eqref{eq:energy-u} follows easily from \eqref{eq:energy-estimate-small-hessian}, \eqref{eq:estimate-big-hessian} and \eqref{eq:estimate-deviation-T-S}:
    \begin{align*}
      \int_{B_\sigma} |\Dd^2 u|^2
      &\leq \int_{\{ |\Dd^2 u| \leq E^{\delta/2} \}} |\Dd^2 u|^2 + \int_{\{ |\Dd^2 u| > E^{\delta/2} \}} |\Dd^2 u|^2 \\
      &\leq C \int_{\Cyl_\sigma} |\vec{S} - \vec{T}|^2 \, \dd \|T\| + C E + C \int_{\{ |\Dd^2 u| > E^{\delta/2} \}} |\Dd^2 u| + |\Dd^2 u|^{2n+2} \\
      &\leq C E^{1+\delta/2} + C \mu + C E + C E^{1+\delta} \leq C E + C \mu.
    \end{align*}

  \item \emph{Excess improvement by tilting.}
    We finally show the bounds \eqref{eq:rotation-excess-improvement} and \eqref{eq:excess-improvement}. We start by showing that
    \begin{equation}
    \label{eq:sup-tilt-S}
    \sup_{B_{2\rho}} |\vec{S} - \vec{S}(0)|^2 \leq C \rho^2 (E + \mu)
  \end{equation}
  holds for any $0 < \rho < \tfrac{\sigma}{4}$, provided that $E, \mu \leq \eps_\rho$, with $C$ independent of $\rho$.
  Recall the estimates \eqref{eq:interior-bound-hessian} and \eqref{eq:interior-closeness-hessian} in the ball of radius $2\rho < \tfrac{\sigma}{2}$ in light of the bound \eqref{eq:energy-u}:
  \[
    \sup_{B_{2\rho}} |\Dd^2 u|^2
    \leq C \int_{B_\sigma} |\Dd^2 u|^2
    \leq C (E + \mu)
  \]
  \[
    \sup_{B_{2\rho}} |\Dd^2 u - \Dd^2 u(0)|^2
    \leq C \rho^2 \int_{B_\sigma} |\Dd^2 u|^2
    \leq C \rho^2 (E + \mu).
  \]
  After making $\eps_\rho$ smaller if necessary, we may assume that $|\Dd^2 u|^2 \leq \tfrac{1}{2}$ in $B_{2\rho}$. Hence \eqref{eq:estimate-jacobian} implies that
  \[
    \Jac(\Phi^u(x)) = 1 + \Ord(|\Dd^2 u(x)|^2) = 1 + \Ord(E + \mu)
  \]
  in $B_{2\rho}$. Since $\vec{S}(x) = \vec{T}^u(x)$ where $\vec{T}^u$ is as in \eqref{eq:tangent-graph}, we have that
  \[
    \vec{S}(x)
    = \frac{(\Phi^u)_\# \vec{\pi}_0(x)}{\Jac(\Phi^u)(x)}
    = \frac{(\Phi^u)_\# \vec{\pi}_0(x)}{1 + \Ord(E + \mu)}
    = \left( (\Phi^u)_\# \vec{\pi}_0(x) \right) (1 + \Ord(E + \mu)),
  \]
  so by \eqref{eq:expansion-pushfwd} this is
  \[
    \vec{S}(x)
    = \left( \vec{\pi}_0 + \sum_{i,j=1}^n \del_{ij}u(x) \vec{\pi}_i^j + \Ord(E + \mu) \right) (1 + \Ord(E + \mu)).
  \]
  Hence
  \[
    \vec{S}(x) - \vec{S}(0)
    = (\vec{\pi}_0 - \vec{\pi}_0) + \left( \sum_{i,j=1}^n (\del_{ij}u(x) - \del_{ij}u(0)) \vec{\pi}_i^j \right) + \Ord(E + \mu),
  \]
  which, after choosing $\eps_\rho$ small enough, gives
  \[
    |\vec{S}(x) - \vec{S}(0)|
    \leq C \rho \sqrt{E + \mu} + C (E + \mu)
    \leq C \rho \sqrt{E + \mu}
  \]
  and \eqref{eq:sup-tilt-S} immediately follows. Next we estimate, using \eqref{eq:estimate-deviation-T-S},
  \begin{align*}
    \frac{1}{2} \int_{\Cyl_{2\rho}} |\vec{T} - \vec{S}(0)|^2 \, \dd \|T\|
      &\leq \int_{\Cyl_{2\rho}} |\vec{T} - \vec{S}|^2 \, \dd \|T\| + \int_{\Cyl_{2\rho}} |\vec{S} - \vec{S}(0)|^2 \, \dd \|T\| \\
      &\leq \int_{\Cyl_{\sigma}} |\vec{T} - \vec{S}|^2 \, \dd \|T\| + \sup_{B_{2\rho}} |\vec{S} - \vec{S}(0)|^2 \|T\|(\Cyl_{2\rho}) \\
      &\leq C E^{1+\delta/2} + C \mu + C \rho^2 (E + \mu) (\rho^n + E) \\
      &\leq C E^{1+\delta/2} + C \mu + C \rho^{2+n} E.
  \end{align*}
  Thus
  \begin{equation}
    \label{eq:excess-wrong-cylinder}
    \rho^{-n} \frac{1}{2} \int_{\Cyl_{2\rho}} |\vec{T} - \vec{S}(0)|^2 \, \dd \|T\|
    \leq C (\rho^2 + \rho^{-n} E^{\delta/2}) E + C \rho^{-n} \mu.
  \end{equation}
  On the other hand, thanks to \eqref{eq:tilt-gradient-comparable}, if we define $\vec{\pi}_1 := \vec{S}(0)$, we have
  \[
    |\vec{\pi}_1 - \vec{\pi}_0|^2
    = |\vec{S}(0) - \vec{\pi}_0|^2
    \leq C |\Dd^2 u(0)|^2 
    \leq C (E + \mu),
  \]
  which is \eqref{eq:rotation-excess-improvement}. Next we need to show that $\spt T_{1/2} \cap \Cyl_{\rho}^{\pi_1} \subset \Cyl_{2\rho}^{\pi_0}$. If $\xi \in \spt T_{1/2}$, by the height bound \eqref{eq:height-bound-y} it holds that $|\mathbf{q}^{\pi_0}(\xi)| \leq \rho$ as long as we take $\eps_\rho$ small enough. Then \cref{lem:modulus-from-projections} below shows that $|\Pi(\xi)| \leq 2\rho$, so in particular $\xi \in \Cyl_{2\rho}^{\pi_0}$. Combining this with \eqref{eq:excess-wrong-cylinder} we get
  \[
    \rho^{-n} \frac{1}{2} \int_{\Cyl_{\rho}^{\pi_1}} |\vec{T} - \vec{\pi}_1|^2 \, \dd \|T_{1/2}\|
    \leq C (\rho^2 + \rho^{-n} E^{\delta/2}) E + C \rho^{-n} \mu.
  \]
  From this, \eqref{eq:excess-improvement} finally follows by choosing $\eps_\rho \leq \rho^{2(n+2)/\delta}$. \qedhere
\end{steps}
\end{proof}

We still need to prove:
\begin{lemma}
  \label{lem:modulus-from-projections}
  Let $\vec{\pi}_0, \vec{\pi}_1$ be two oriented Legendrian $n$-planes such that $|\mathbf{p}^{\pi_1} - \mathbf{p}^{\pi_0}|^2 \leq \tfrac{1}{8}$ and let $r > 0$. Suppose that $\xi \in \Heis^n$ satisfies $|\mathbf{p}^{\pi_1}(\xi)| \leq r$ and $|\mathbf{q}^{\pi_0}(\xi)| \leq r$. Then $|\Pi(\xi)| \leq 2r$.
  \begin{proof}
    Write $z = \Pi(\xi) \in \R^{2n}$ and compute
    \begin{align*}
      |z|^2
      = |\mathbf{p}^{\pi_0}(z)|^2 + |\mathbf{q}^{\pi_0}(z)|^2
      &\leq \left( |\mathbf{p}^{\pi_0} - \mathbf{p}^{\pi_1}| |z| + |\mathbf{p}^{\pi_1} (z)|\right)^2 + |\mathbf{q}^{\pi_0}(z)|^2 \\
      &\leq 2 |\mathbf{p}^{\pi_0} - \mathbf{p}^{\pi_1}|^2 |z|^2 + 2 |\mathbf{p}^{\pi_1} (\xi)|^2 + |\mathbf{q}^{\pi_0}(\xi)|^2 \\
      &\leq \frac{1}{4} |z|^2 + 3 r^2.
    \end{align*}
    This immediately implies the inequality.
  \end{proof}
\end{lemma}

We get the power decay of the excess by combining \cref{lem:excess-improvement} with a delicate control of the tilting of the cylinders, for which we follow essentially \cite{regularity-currents-hilbert} with some simplifications.

\begin{lemma}[excess decay]
  \label{lem:excess-decay}
  There exist constants $0 < \epsed < 1$ and $C_\idxed > 0$ depending only on $n, \Lambda / \lambda$ and $Q$ such that if $T \in \IC^\hor_n(\Heis^n)$ satisfies \cref{ass:scale-invt-hypot} with $\pi_0$ in place of $\pi$ and $\eps \leq \epsed$, and in addition $0 \in \spt T$ and $\mu R \leq \epsed$, then for any $\xi \in \spt T \cap \Cyl_{R/8}^{\pi_0}$ and any $0 < r \leq R/8$ it holds that
  \begin{equation}
    \label{eq:excess-pi0-controlled}
    \Exc(T, \Cyl_{r}^{\pi_0}(\xi)) \leq C_\idxed \left( \Exc(T, \Cyl_R^{\pi_0}(0)) + \mu R \right)
  \end{equation}
  and
  \begin{equation}
    \label{eq:excess-decay}
    \Exc(T \res \Cyl_{R/2}^{\pi_0}, \Cyl_{r}^{\pi(\xi, r)}(\xi)) \leq C_\idxed \frac{r}{R} \left(\Exc(T, \Cyl_R^{\pi_0}(0)) + \mu R \right)
  \end{equation}
  for a plane $\vec{\pi}(\xi, r) \in \mathscr{L}_n$ satisfying
  \begin{equation}
    \label{eq:nearly-horizontal-plane}
    |\vec{\pi}(\xi, r) - \vec{\pi}_0|^2 \leq C_\idxed (\Exc(T, \Cyl_R^{\pi_0}) + \mu R).
  \end{equation}

  \begin{proof}
    Let $\rho := \min \left\{ \tfrac{1}{2 C_\idximp}, \tfrac{1}{72} \right\}$, where $C_\idximp$ is the constant from \cref{lem:excess-improvement}, and let $\eps_\rho$ be the corresponding parameter. Let also $T_0 := T \res \Cyl_{R/2}^{\pi_0}(0)$ and observe that by \eqref{eq:rescaled-height-bound-y} we can make sure that $|\mathbf{q}^{\pi_0}(\zeta)| \leq \tfrac{R}{16}$ for every $\zeta \in \spt T_0$. Fix any $\xi \in \spt T \cap \Cyl_{R/8}^{\pi_0}$ and let $T_\xi := (\ell_{\xi^{-1}})_\# T_0$. Since $\mathbf{q}^{\pi_0}$ is a homomorphism, the previous bound and the triangle inequality imply that
    \begin{equation}
      \label{eq:height-bound-aux}
      |\mathbf{q}^{\pi_0}(\zeta)| \leq \frac{R}{8} \qquad \forall \zeta \in \spt T_\xi.
    \end{equation}
    We define the quantity
    \begin{equation}
      \mathcal{E}(\xi, r, \vec{\pi})
      := \max \left\{ \Exc(T_0, \Cyl_r^{\pi}(\xi)), 2 C_\idximp \rho^{-n-1} r \mu \right\}
      = \max \left\{ \Exc(T_\xi, \Cyl_r^{\pi}), 2 C_\idximp \rho^{-n-1} r \mu \right\}
    \end{equation}
    for any $0 < r \leq R/8$ and any oriented Legendrian $n$-plane $\vec{\pi}$. It is clear that if we take $\epsed$ small enough, then $\Exc(T, \Cyl_R^{\pi_0}(0)) \leq \epsed \leq 8^{-n} \eps_\rho$ and $2 C_\idximp \rho^{-n-1} \tfrac{R}{8} \mu \leq \eps_\rho$, hence by using \eqref{eq:excess-def} we have
    \begin{equation}
      \label{eq:starting-excess}
      \mathcal{E}(\xi, R/8, \vec{\pi}_0) \leq 8^n \epsed \leq \eps_\rho.
    \end{equation}
    Let $0 < \theta < \tfrac{\pi}{4}$ be such that $\tan \theta = \tfrac{\rho}{8}$, and define the set
    \[
      \mathbf{K}_r := \left\{ \zeta \in \Heis^n : |\mathbf{p}^{\pi_0}(\zeta)| \leq r + \tan \theta \, |\mathbf{q}^{\pi_0}(\zeta)| \right\}.
    \]
    Note that $\mathbf{K}_r$ is asymptotic to a cone with opening angle $\theta$. For each integer $k \geq 0$ let $r_k := \rho^k R / 8$. We will prove by induction that for every $k \geq 1$ there exists a plane $\vec{\pi}_k$ such that
    \begin{equation}
      \label{eq:gen-excess-decay}
      \mathcal{E}(\xi, r_k, \vec{\pi}_k)
      \leq 
      \rho^k \mathcal{E}(\xi, R/8, \vec{\pi}_0),
    \end{equation}
    \begin{equation}
      \label{eq:induction-inclusion}
      \spt T_\xi \cap \mathbf{K}_{r_{k}/8} \subset \Cyl_{r_k / 2}^{\pi_k}(0)
    \end{equation}
    and
    \begin{equation}
      \label{eq:induction-tilt}
    |\vec{\pi}_{k} - \vec{\pi}_{k-1}| \leq C_\idxtilt \rho^{k/2} \mathcal{E}(\xi, R/8, \vec{\pi}_0)^{1/2}.
    \end{equation}
    Equation \eqref{eq:gen-excess-decay} is clear for $k = 0$, and \eqref{eq:induction-inclusion} follows from the following computation: for any $\zeta \in \spt T_\xi \cap \mathbf{K}_{r_0/8}$,
    \[
      |\mathbf{p}^{\pi_0}(\zeta)|
      \leq \frac{r_0}{8} + \tan \theta \,|\mathbf{q}^{\pi_0}(\zeta)|
      \leq \frac{r_0}{8} + \frac{\rho}{8} \frac{R}{8}
      \leq \frac{r_0}{8} + \frac{1}{8} r_0
      < \frac{r_0}{2},
    \]
    thanks to \eqref{eq:height-bound-aux}.
    To prove the claims \eqref{eq:gen-excess-decay}-\eqref{eq:induction-tilt} for $k$, we need to check, using that they hold up to $k - 1$, that the assumptions of \cref{lem:excess-improvement} (which contain those of \cref{cor:rescaled-height-bound}) are satisfied for $T_\xi$ in the cylinder $\Cyl_{r_{k-1}}^{\pi_{k-1}}$. First, \eqref{eq:ass-f-minimizing} and \eqref{eq:ass-density-q} are clear with $\mathcal{U} = \BHeis_{R/2}$. Next observe that, summing \eqref{eq:induction-tilt} for $j = 1, \ldots, k-1$ in place of $k$ gives
    \begin{equation}
      \label{eq:plane-tilt}
      |\vec{\pi}_{k-1} - \vec{\pi}_{0}|
      \leq C \mathcal{E}(\xi, R/8, \vec{\pi}_0)^{1/2}
      \leq C \sqrt{\epsed}
    \end{equation}
    for a constant $C$ independent of $k$. Moreover, \eqref{eq:gen-excess-decay} and \eqref{eq:starting-excess} trivially give that
    \[
      \Exc \left(T_\xi, \Cyl_{r_{k-1}}^{\pi_{k-1}}(0) \right)
      \leq \mathcal{E}(\xi, r_{k-r}, \vec{\pi}_{k-1})
      \leq 8^n \epsed
      \leq \eps_\rho,
    \]
    which is \eqref{eq:ass-excess}. Next consider a continuous path $\vec{\pi} : [0, 1] \to \mathscr{L}_0$ with $\vec{\pi}(0) = \vec{\pi}_0, \vec{\pi}(1) = \vec{\pi}_{k-1}$ and $|\mathbf{p}^{\pi(t)} - \mathbf{p}^{\pi_0}|^2 \leq \tfrac{1}{8}$ for all $t$, which exists by \eqref{eq:plane-tilt} if $\epsed$ is small enough. The height bound \eqref{eq:height-bound-aux} and \cref{lem:modulus-from-projections} then imply that
    \begin{equation}
      \spt T_\xi \cap \Cyl_{r_{k-1}}^{\pi(t)}
      \subset \spt T_\xi \cap \Cyl_{R/8}^{\pi(t)}
      \subset \left\{ \zeta \in \Heis^n : |\Pi(\zeta)| < \frac{R}{4} \right\}
      \subset \Cyl_{R/4}^{\pi_{0}}.
    \end{equation}
    On the other hand we have that $\spt \del T_0 \cap \Cyl_{R/2}^{\pi_0} = \varnothing$ and hence $\spt \del T_\xi \cap \Cyl_{R/4}^{\pi_0} = \varnothing$. These two conditions imply that $\spt \del T_\xi \cap \Cyl_{r_{k-1}}^{\pi(t)} = \varnothing$, and setting $t = 1$ gives \eqref{eq:ass-no-bdry}. Finally, by the constancy theorem, for each $t \in [0, 1]$ we have that
    \[
      \mathbf{p}^{\pi(t)}_\# \left(T_\xi \res \Cyl_{r_{k-1}}^{\pi(t)}\right) = Q(t) \DBrack{B_{r_{k-1}}^{\pi(t)}}
    \]
    for some integer $Q(t)$. It is easy to see that $Q(t)$ is continuous and thus constant, proving \eqref{eq:ass-q-cover}. We have thus shown that \cref{ass:scale-invt-hypot} is satisfied, and it is clear that \cref{ass:metric-h} holds as well with the same constants $\lambda, \Lambda, \mu$.

    Now we may apply \cref{lem:excess-improvement} to obtain a plane $\vec{\pi}_k$ satisfying
    \begin{equation}
      \label{eq:excess-estimate-aux}
      \Exc \left(T_\xi \res \Cyl_{r_{k-1}/2}^{\pi_{k-1}}, \Cyl_{r_{k}}^{\pi_{k}} \right)
      \leq C_\idximp \rho^2 \Exc(T_\xi, \Cyl_{r_{k-1}}^{\pi_{k-1}}) + C_\idximp \rho^{-n} r_{k-1} \mu
    \end{equation}
    and $|\vec{\pi}_k - \vec{\pi}_{k-1}| \leq C \mathcal{E}(\xi, r_{k-1}, \vec{\pi}_{k-1})^{1/2}$, which implies \eqref{eq:induction-tilt} by induction hypothesis.

    Next we show \eqref{eq:induction-inclusion} by making use of \cref{cor:rescaled-height-bound}: if $\zeta \in \spt T_\xi \cap \mathbf{K}_{r_k/8}$, then obviously $\zeta \in \spt T_\xi \cap \mathbf{K}_{r_{k-1}/8} \subset \Cyl_{r_{k-1}/2}^{\pi_{k-1}}$ by induction hypothesis. Now using \eqref{eq:rescaled-height-bound-y}, as long as $\epsed$ is small enough, we have that
  \begin{equation}
    \label{eq:height-previous-pi}
    |\mathbf{q}^{\pi_{k-1}}(\zeta)| \leq \frac{\rho}{4} r_{k-1} = \frac{r_k}{4}.
  \end{equation}
  In particular,
  \[
    |\Pi(\zeta)|
    \leq |\mathbf{p}^{\pi_{k-1}}(\zeta)| + |\mathbf{q}^{\pi_{k-1}}(\zeta)|
    \leq \frac{r_{k-1}}{2} + \frac{\rho}{4} r_{k-1} \leq r_{k-1},
  \]
  so that $|\mathbf{q}^{\pi_0}(\zeta)| \leq r_{k-1}$. Now since $\zeta \in \mathbf{K}_{r_k/8}$,
  \begin{equation}
    \label{eq:width-pi0}
    |\mathbf{p}^{\pi_0}(\zeta)|
    \leq \frac{r_k}{8} + \tan \theta \, |\mathbf{q}^{\pi_0}(\zeta)|
    \leq \frac{r_k}{8} + \frac{\rho}{8} r_{k-1}
    = \frac{r_k}{4}.
  \end{equation}
  We may finally apply \cref{lem:modulus-from-projections}, together with \eqref{eq:height-previous-pi}, \eqref{eq:width-pi0} and \eqref{eq:plane-tilt} to deduce that $|\Pi(\zeta)| \leq r_k / 2$ and hence $\zeta \in \Cyl_{r_k / 2}^{\pi_{k}}$. This finishes the proof of \eqref{eq:induction-inclusion}.

  We finally show \eqref{eq:gen-excess-decay}: starting from \eqref{eq:excess-estimate-aux}, we have
    \begin{align*}
      \Exc \left(T_\xi \res \Cyl_{r_{k-1}/2}^{\pi_{k-1}}, \Cyl_{r_{k}}^{\pi_{k}} \right)
      &\leq C_\idximp \rho^2 \Exc(T_\xi, \Cyl_{r_{k-1}}^{\pi_{k-1}}) + C_\idximp \rho^{-n} r_{k-1} \mu \\
      &\leq \frac{1}{2} \rho \Exc(T_\xi, \Cyl_{r_{k-1}}^{\pi_{k-1}}) + \frac{1}{2} \rho \cdot 2 C_\idximp \rho^{-n-1} r_{k-1} \mu \\
      &\leq \rho \max \left\{ \Exc(T_\xi, \Cyl_{r_{k-1}}^{\pi_{k-1}}), 2 C_\idximp \rho^{-n-1} r_{k-1} \mu \right\} \\
      &= \rho \mathcal{E}(\xi, r_{k-1}, \vec{\pi}_{k-1}),
    \end{align*}
    therefore in order to prove that $\mathcal{E}(\xi, r_k, \vec{\pi}_{k}) \leq \rho \mathcal{E}(\xi, r_{k-1}, \vec{\pi}_{k-1})$ it is enough to show that $\spt T_\xi \cap \Cyl_{r_{k}}^{\pi_{k}} \subset \Cyl_{r_{k-1}/2}^{\pi_{k-1}}$, since the linear decay of $r_k$ is clear. Observe that by \eqref{eq:induction-inclusion} for $k-1$, this will follow from the inclusion
    \begin{equation}
      \label{eq:cylinder-in-K}
      \Cyl_{r_k}^{\pi_k} \subset \mathbf{K}_{r_{k-1} / 8}.
    \end{equation}
    To show \eqref{eq:cylinder-in-K}, let $\zeta \in \Cyl_{r_k}^{\pi_k}$ and compute
    \begin{align*}
      |\mathbf{p}^{\pi_0}(\zeta)|
      &\leq |\mathbf{p}^{\pi_0} - \mathbf{p}^{\pi_k}| |\zeta| + |\mathbf{p}^{\pi_k}(\zeta)| \\
      &\leq |\mathbf{p}^{\pi_0} - \mathbf{p}^{\pi_k}| (|\mathbf{p}^{\pi_0}(\zeta)| + |\mathbf{q}^{\pi_0}(\zeta)|) + r_k,
    \end{align*}
    hence, using \eqref{eq:plane-tilt} and \eqref{eq:induction-tilt},
    \begin{align*}
      |\mathbf{p}^{\pi_0}(\zeta)|
      &\leq \frac{|\mathbf{p}^{\pi_0} - \mathbf{p}^{\pi_k}|}{1 - |\mathbf{p}^{\pi_0} - \mathbf{p}^{\pi_k}|} |\mathbf{q}^{\pi_0}(\zeta)| + \frac{1}{1 - |\mathbf{p}^{\pi_0} - \mathbf{p}^{\pi_k}|} r_k \\
      &\leq \frac{\rho}{8} \, |\mathbf{q}^{\pi_0}(\zeta)| + \frac{1}{8\rho} r_k
      = \tan \theta \, |\mathbf{q}^{\pi_0}(\zeta)| + \frac{r_{k-1}}{8}
    \end{align*}
    provided that $\epsed$ is small enough. This shows \eqref{eq:cylinder-in-K} and closes the induction step.

    Finally we show the bounds of the statement: given any $0 < r < R / 8$, if $k \geq 1$ is the integer such that $r_k < r \leq r_{k-1}$, we have that
    \[
      \Exc(T_\xi, \Cyl_{r}^{\pi_{k-1}})
      \leq \rho^{-n} \Exc(T_\xi, \Cyl_{r_{k-1}}^{\pi_{k-1}})
      \leq \rho^{-n} \rho^{k-1} \mathcal{E}(\xi, R/8, \vec{\pi}_0)
      \leq 8^{n+1} \rho^{-n-1} \frac{r}{R} (\Exc(T, \Cyl_{R}^{\pi_0}(\xi)) + C \mu R),
    \]
    which is \eqref{eq:excess-decay} with $\vec{\pi}(\xi, r) = \vec{\pi}_{k-1}$. The bound on \eqref{eq:nearly-horizontal-plane} is immediate from \eqref{eq:plane-tilt}. To show \eqref{eq:excess-pi0-controlled}, observe that by \eqref{eq:induction-inclusion},
    \[
      \spt T_\xi \cap \Cyl_{r_{k}/8}^{\pi_0}
      \subset \spt T_\xi \cap \mathbf{K}_{r_{k}/8}
      \subset \spt T_\xi \cap \Cyl_{r_k / 2}^{\pi_k}
      \subset \spt T_\xi \cap \Cyl_{r_k}^{\pi_k}
    \]
    and $\Cyl_r^{\pi_0}(\xi) \subset \Cyl_{R/2}^{\pi_0}(0)$ for any $r \leq R / 8$.
    Therefore we can compute, for any $k \geq 0$,
    \begin{align*}
      \Exc(T, \Cyl_{r_k/8}^{\pi_0}(\xi))
      &= \Exc(T_0, \Cyl_{r_k/8}^{\pi_0}(\xi))
      = \Exc(T_\xi, \Cyl_{r_k/8}^{\pi_0}(0)) \\
      &= \frac{8^n}{2r_k^n} \int_{\Cyl_{r_k/8}^{\pi_0}} |\vec{T}_\xi - \vec{\pi}_0|^2 \, \dd \|T_\xi \| \\
      &\leq \frac{8^n}{2r_k^n} \int_{\Cyl_{r_k}^{\pi_k}} |\vec{T}_\xi - \vec{\pi}_0|^2 \, \dd \|T_\xi \| \\
      &\leq \frac{8^n}{r_k^n} \int_{\Cyl_{r_k}^{\pi_k}} |\vec{T}_\xi - \vec{\pi}_k|^2 \, \dd \|T_\xi \| + \frac{8^n}{r_k^n} \int_{\Cyl_{r_k}^{\pi_k}} |\vec{\pi}_k - \vec{\pi}_0|^2 \, \dd \|T_\xi \| \\
      &= 2 \cdot 8^n \Exc(T_\xi, \Cyl_{r_k}^{\pi_k}) + 8^n |\vec{\pi}_k - \vec{\pi}_0|^2 \frac{\|T_\xi\|(\Cyl_{r_k}^{\pi_k})}{r_k^n} \\
      &\leq 2 \cdot 8^n \Exc(T_\xi, \Cyl_{r_k}^{\pi_k}) + 8^n C (\Exc(T, \Cyl_R^{\pi_0}) + \mu R) (Q \omega_n + \Exc(T_\xi, \Cyl_{r_k}^{\pi_k})) \\
      &\leq C (\Exc(T, \Cyl_{R}^{\pi_0}) + \mu R),
    \end{align*}
    where we have used \eqref{eq:plane-tilt}, \eqref{eq:excess-def} and \eqref{eq:excess-decay}. Then \eqref{eq:excess-pi0-controlled} follows easily.
  \end{proof}
\end{lemma}

Deducing regularity from here is rather standard.

\begin{proof}[\hypertarget{proof-eps-reg}{Proof of \cref{thm:eps-reg}}]
  Since the derivative of the map
  \[
    A \in \Sym^{n}(\R) \; \longmapsto \; \frac{A}{\sqrt{\det(\id + A^2)}} \in \Sym^n(\R)
  \]
  at $A = 0$ is the identity, there exists some $0 < \gamma \leq 1$ depending only on $n$ such that this map is a diffeomorphism restricted to $\{ A \in \Sym^n(\R) : |A| \leq \gamma \}$ and hence
  \begin{equation}
    \label{eq:tangent-controls-hessian}
    |A - B| \leq C \left| \frac{A}{\sqrt{\det(\id + A^2)}} - \frac{B}{\sqrt{\det(\id + B^2)}} \right|
  \end{equation}
  whenever $|A|, |B| \leq \gamma$, for a dimensional constant $C > 0$. We construct the $C^{1,1}$ function $f : B_{R/72}^{\pi_0} \to \R$ by applying \cref{lem:lip-approx} on the ball $B_{R/8}$ and with this parameter $\gamma$. As long as $\epsreg \leq \epsed$, we may use \cref{lem:excess-decay} and get that
  \[
    \Exc(T, \Cyl_r^{\pi_0}(x)) \leq C_\idxed \epsreg
  \]
  for every $x \in B_{R/72}$ and every $0 < r < R/9$. Hence, recalling the construction of the set $G_\gamma$ from \eqref{eq:good-set-lip}, if $\epsreg$ is chosen small enough, then $G_\gamma$ is the whole $B_{R/72}$ and it follows that for every $x \in B_{R/72}$, $\spt T \cap (\mathbf{p}^{\pi_0})^{-1}(x) = \{ \Phi^f(x) \}$. Therefore $T \res \Cyl_{R/72}^{\pi_0}$ is supported on the graph of $\Phi^f$. Moreover, it follows from the proof of \cref{lem:lip-approx} that the set $K \subset B_{R/72}$ appearing in \eqref{eq:lip-bound-bad-set} has full measure, which implies that for almost every $x \in B_{R/72}$,
  \[
    \vec{T}(\Phi^f(x)) = \vec{T}^f(\Phi^f(x)) =: \vec{T}^f(x).
  \]
  By \eqref{eq:nearly-horizontal-plane}, if $\epsreg$ is small enough, all the Legendrian planes $\vec{\pi}(x, r) := \vec{\pi}(\Phi^f(x), r)$ with $0 < r \leq \tfrac{R}{36}$ are close enough to $\vec{\pi}_0$ that they can be written as the graph of a symmetric linear map with small norm. Namely,
  \[
    \vec{\pi}(x, r)
    = \frac{(\id, L(x, r))_\# \vec{\pi}_0}{|(\id, L(x, r))_\# \vec{\pi}_0|}
    = \frac{\vec{\pi}_0 + L(x, r)_\# \vec{\pi}_0}{\sqrt{\det(\id + L(x,r)^2)}}
  \]
  for a matrix $L(x, r) \in \Sym^n(\R)$ that satisfies $|L(x, r)| \leq \gamma$. Now we compute using \eqref{eq:tangent-controls-hessian}, \eqref{eq:expansion-pushfwd}, \eqref{eq:jacobian-as-det}, \eqref{eq:tangent-graph} and the coarea formula:
  \begin{align*}
    r^{-n} \int_{B_r(x)} |\Dd^2 f(y) - L(x,2r)|^2 \, \dd y
    &\leq C r^{-n} \bigintss_{B_r(x)} \left| \frac{\Dd^2 f(y)}{\sqrt{\det(\id + (\Dd^2 f(y))^2)}} - \frac{L(x,2r)}{\sqrt{\det(\id + L(x,2r)^2)}} \right|^2 \, \dd y \\
    &\leq C r^{-n} \bigintss_{B_r(x)} \left| \frac{(\Phi^f)_\# \vec{\pi}_0(y)}{\sqrt{\det(\id + (\Dd^2 f(y))^2)}} - \frac{(\id, L(x, 2r))_\# \vec{\pi}_0}{\sqrt{\det(\id + L(x,2r)^2)}} \right|^2 \, \dd y \\
    &= C r^{-n} \int_{B_r(x)} | \vec{T}^f(y) - \vec{\pi}(x,2r)|^2 \, \dd y \\
    &\leq C r^{-n} \int_{\Cyl_r^{\pi_0}(x)} | \vec{T} - \vec{\pi}(x,2r)|^2 \, \dd \|T \res \Cyl_{R/2}^{\pi_0}\| \\
    &\leq C r^{-n} \int_{\Cyl_{2r}^{\pi(x,2r)}(\Phi^f(x))} | \vec{T} - \vec{\pi}(x,2r)|^2 \, \dd \|T \res \Cyl_{R/2}^{\pi_0}\| \\
    &\leq C \Exc(T \res \Cyl_{R/2}^{\pi_0}, \Cyl_{2r}^{\pi(x,2r)}(\Phi^f(x))) \\
    &\leq C \frac{r}{R} \left( \Exc(T, \Cyl_{R}^{\pi_0}) + R\mu\right).
  \end{align*}
  Here we have applied the height bound together with \cref{lem:modulus-from-projections} to change the cylinder of integration, as earlier in this section. Replacing $L(x, 2r)$ by the average of $\Dd^2 f$ over $B_r(x)$ in the left hand side only decreases the integral, therefore Campanato's theorem implies that $f \in C^{2,1/2}(B_{R/72})$ with
  \[
    \left[ \Dd^2 f \right]_{C^{1/2}(B_{R/72})} \leq C R^{-1/2} \left( \Exc(T, \Cyl_{R}^{\pi_0}) + R\mu\right)^{1/2}.
  \]
  This together with \eqref{eq:energy-lip} gives, for every $x \in B_{R/72}$,
  \begin{align*}
    |\Dd^2 f(x)|^2
    &\leq C R^{-n} \int_{B_{R/72}} \left(|\Dd^2 f(x) - \Dd^2 f(y)|^2 + |\Dd^2 f(y)|^2\right) \, \dd y \\
    &\leq C \left( \Exc(T, \Cyl_{R}^{\pi_0}) + R\mu\right) + C \gamma^{-n(n+1)} \Exc(T, \Cyl_R^{\pi_0}) \\
    &\leq C \left( \Exc(T, \Cyl_{R}^{\pi_0}) + R\mu\right)
  \end{align*}
  and the rest of the estimates in \eqref{eq:estimates-eps-reg} now follow by integrating starting at $\Dd f(0) = 0$ and $f(0) = 0$. Finally, now that we know that $\Phi^f$ is $C^{1,1/2}$, the expression \eqref{eq:eps-reg-graph} is a consequence of the constancy theorem.
\end{proof}

\appendix

\section{Metric geometry of the Heisenberg space}
\label{sec:metric-geometry}

Here we recall some facts about the intrinsic metric geometry of $\Heis^n$ with respect to its Carnot--Carath\'eodory distance $d_\mathrm{CC}$. This distance is defined as the infimum of the lengths, computed with the standard subriemannian metric, of smooth horizontal curves joining the two extrema. The infimum is always attained and the resulting distance has an explicit expression, although less convenient than the one induced by the Folland--Kor\'anyi norm \eqref{eq:fk-gauge}, which is bilipschitz equivalent to it. In this section we work with the former distance because we need some equalities with precise constants that we have only been able to find in the literature with $d_\mathrm{CC}$.

A lot of progress has been made in recent years in understanding the space $(\Heis^n, d_\mathrm{CC})$ since the influential book of Gromov \cite{gromov-cc-spaces-from-within}, and a rather satisfactory theory has been built, which encompasses results about Lipschitz extensions, fillings and rectifiable horizontal subsets. In particular, at the time when \cite{schoen-wolfson-jdg} was published, the theory of integral currents in metric spaces had just been born \cite{ambrosio-kirchheim}, and only the isoperimetric inequality for surfaces of Allcock \cite{allcock} was known.

\subsection{Comparison between horizontal currents and metric currents}
The following theorem, relating metric currents in $(\Heis^n, d_\mathrm{CC})$ and horizontal currents in the sense of \cref{defn:horizontal-current}, follows from the work of Ambrosio--Kirchheim \cite{ambrosio-kirchheim, ambrosio-kirchheim-rectifiable} and Williams \cite{williams-currents-carnot}. Note that an analogous result also applies to more general Carnot groups.

We need to endow $\Heis^n$ with a left-invariant Riemannian metric $g_0$ which agrees with the standard subriemannian metric $g_{\Heis^n}$ on the contact distribution; let $d_0$ denote its induced distance. Note that any such distance is locally bi-Lipschitz equivalent to the Euclidean distance coming from the identification $\Heis^n \simeq \R^{2n+1}$ using exponential coordinates. Since the lengths of horizontal curves computed with respect to $g_0$ and $g_{\Heis^n}$ agree, the identity map $I : (\Heis^n, d_\mathrm{CC}) \to (\Heis^n, d_0)$ is $1$-Lipschitz.

Let $\ICmet_k(\Heis^n, d_\mathrm{CC})$ denote the space of $k$-dimensional metric integral currents in the sense of Ambrosio--Kirchheim \cite{ambrosio-kirchheim}, and let $\IC^\hor_k(\Heis^n)$ denote the space of Rumin horizontal currents \cref{defn:horizontal-current}. This is a subset of the space of Federer--Fleming (locally) integral currents $\IC_k(\Heis^n)$, which can be identified with Ambrosio--Kirchheim currents in $\IC_k(\Heis^n, d_0)$.

\begin{thm}
  \label{thm:rumin-metric-currents}
  For any $0 \leq k \leq n$, the pushforward map $I_\# : \ICmet_k(\Heis^n, d_\mathrm{CC}) \to \IC^\hor_k(\Heis^n)$ is well defined and realizes an isomorphism of abelian groups.

  Moreover, given a metric current $T \in \ICmet_k(\Heis^n, d_\mathrm{CC})$, there exists a set $S \subset \Heis^n$ which is $k$-rectifiable with respect to both $d_\mathrm{CC}$ and in the Euclidean sense, an orientation $\vec{S}$ of $\operatorname{Tan}(S, \xi)$ defined for $\Haus^k_{d_0}$-almost every $\xi \in S$, and an integer-valued $\Haus^k_{d_0}$-integrable function $\theta : S \to \Z^+$, such that the following representation formulae hold:
  \begin{equation}
    \label{eq:representation-current}
    (I_\# T)(\alpha) = \int_S \langle \alpha(\xi), \vec{S}(\xi) \rangle \theta(\xi) \, \dd \Haus^k_{d_0}(\xi)
    \qquad \text{for any } \alpha \in \mathcal{D}^k(\Heis^n)
  \end{equation}
  \begin{equation}
    \label{eq:rumin-metric-mass}
    \|T\| = \|I_\# T \|
    = \theta \Haus^k_{d_\mathrm{CC}} \res S 
    = \theta \Haus^k_{d_0} \res S.
  \end{equation}
  Here the density $\theta$ can be computed as
  \begin{equation}
    \label{eq:density-cc-riem}
    \theta(\xi) = \Theta^k_{d_0}(\|I_\# T\|, \xi) = \Theta^k_{d_\mathrm{CC}}(\|T\|, \xi)
  \end{equation}
  for $\Haus^k_{d_0}$-almost every $\xi \in S$.
  \begin{proof}
    Theorem 1.6 of \cite{williams-currents-carnot} gives that, after identifying integral metric currents in $(\Heis^n, d_0)$ with Federer--Fleming integral currents in $(\Heis^n, g_0)$, $I_\#$ is an isomorphism, but just with the inequality
    \[
      \label{eq:williams-masses}
      \frac{k!}{(2n)!} \|T\| \leq \|I_\# T\| \leq \|T\|.
    \]
    To see that the masses coincide, use \cite[Theorem 4.5]{ambrosio-kirchheim} to write 
    \begin{equation}
      \label{eq:current-parametrization}
      T = \sum (f_i)_\# \DBrack{\theta_i},
      \qquad
      \|T\| = \sum \|(f_i)_\# \DBrack{\theta_i}\|
      \qquad \text{and} \qquad
      S = \bigcup_i f_i(K_i),
    \end{equation}
    where $K_i \subset \R^k$ are compact sets, $f_i : \R^k \to \Heis^n$ are bi-Lipschitz with respect to $d_\mathrm{CC}$, $\theta_i \in L^1(\R^k, \Z)$, $f_i(K_i)$ are pairwise disjoint, and $\|T\|$ is concentrated on $S$.
  It follows from the work \cite{scott-d-pauls} (see also \cite{regular-submanifolds-heis, mattila-serapioni-sc})
    that the approximate tangent cones of $S$ exist $\Haus^k_{d_\mathrm{CC}}$-almost everywhere and are all isomorphic to $\R^k$.
    Hence their area factor (in the sense of \cite{ambrosio-kirchheim, ambrosio-kirchheim-rectifiable}) is $1$ and as a result, the discussion in \cite[p.~58]{ambrosio-kirchheim} establishes that
    \begin{equation}
      \label{eq:representation-mass}
      \|T\| = \theta \Haus^k_{d_\mathrm{CC}} \res S,
    \end{equation}
    where $\theta(\xi) = \Theta^k_{d_\mathrm{CC}}(\|T\|, \xi)$ by \cite[Theorem 5.4]{ambrosio-kirchheim-rectifiable} (see also the proof of \cite[Theorem 4.6]{ambrosio-kirchheim}). Comparing the two expressions for the mass of $T$, it follows that for $\|T\|$-almost every $\xi$, $\theta(\xi) = |\theta_i(x)|$ whenever $\xi = f_i(x)$.

    Now an application of the area formula from \cite[Theorem~6.8]{introduction-heisenberg} (plus a standard argument to pass from measures of sets to integrals of functions) on \eqref{eq:representation-mass} together with an application of the standard area formula gives
    \begin{align*}
      \|T\|(\psi)
      &= \int_{S} \psi(\xi) \theta(\xi) \, \dd \Haus^k_{d_\mathrm{CC}}(\xi) \\
      &= \sum_i \int_{K_i} \psi(f_i(x)) \theta(f_i(x)) \Jac(f_i)(x) \, \dd x \\
      &= \int_{S} \psi(\xi) \theta(\xi) \, \dd \Haus^k_{d_0}(\xi)
      = \|I_\# T\|(\psi)
    \end{align*}
    for any $\varphi \in C_c(\Heis^n)$. This shows all the equalities in \eqref{eq:rumin-metric-mass}. As a consequence, we get \eqref{eq:density-cc-riem} again by \cite[Theorem 5.4]{ambrosio-kirchheim-rectifiable}. The formula \eqref{eq:representation-current} is now clear.
  \end{proof}
\end{thm}

\begin{rmk}
  \label{rmk:williams-arbitrary-contact}
  This identification holds in any contact manifold $(M^{2n+1}, \Xi, g)$ with a subriemannian metric $g$ in the contact distribution $\Xi$.
  Since we do not use this anywhere in the paper, we omit the details.
\end{rmk}

\begin{prop}
  \label{prop:blowup-heisenberg}
  Let $T \in \IC^\hor_k(\Heis^n)$. Then for $\|T\|$-almost every $\xi$,
  \begin{equation}
    \label{eq:blowup-heisenberg}
    (\delta_{1/\rho} \circ \ell_{\xi^{-1}})_\# T \weakto \Theta^k_{d_\mathrm{CC}}(\|T\|, \xi) \DBrack{\vec{T}(\xi)} \qquad \text{in } \mathcal{D}_k(\Heis^n) \text{ as } \rho \searrow 0,
  \end{equation}
  where $\vec{T}(\xi)$ is the oriented approximate tangent plane to $T$ at $\xi$ seen as a subgroup of $\Heis^n$.
  \begin{proof}
    This follows from \cite[Theorem~3.14]{mattila-serapioni-sc} by a standard argument: they show that, given a rectifiable set $S$ as in \cref{thm:rumin-metric-currents}, for $\Haus^k$-almost every $\xi \in S$ there exists a unique horizontal $k$-plane $\pi_\xi$ through $\xi$ such that
    \[
      \lim_{r \to 0} \frac{\Haus^k(S \cap \mathscr{B}_r(\xi) \setminus X(\xi, \pi_\xi, s))}{r^k} = 0
    \]
    for every $0 < s < 1$. Here $\mathscr{B}_r(\xi)$ denotes a ball with respect to the Carnot--Carath\'eodory distance and $X(\xi, \pi_\xi, s)$ is a certain cone of opening $s$ centered at $\xi$ around $\pi_\xi$. We omit the precise definition of the cones $X(\xi, \pi_\xi, s)$ since the only fact that we will use is that $\bigcap_{0 < s < 1} X(\xi, \pi_\xi, s) = \pi_\xi$.

    Moreover, the plane $\pi_\xi$ corresponds with the Riemannian tangent plane of $S$ at $\xi$. This follows for example from \cite[Theorem~3.5]{regular-submanifolds-heis} or \cite[Remark~3.8]{mattila-serapioni-sc}. Therefore, if we orient it appropriately and denote $\pi = \ell_{\xi^{-1}}(\pi_\xi)$ we have that $\vec{\pi} = \vec{T}(\xi)$.

    If $\xi$ is in addition a Lebesgue point of $\theta$ with respect to $\Haus^k \res S$, then we evidently have
    \begin{equation}
      \label{eq:convergence-away-cone}
      \lim_{r \to 0} \frac{\|T\|(\mathscr{B}_r(\xi) \setminus X(\xi, \pi_\xi, s))}{r^k} = 0
    \end{equation}
    for every $0 < s < 1$. Define the blow-up maps $\eta_{\xi, \rho}(\zeta) := \delta_{\rho^{-1}}(\ell_{\xi^{-1}}(\zeta))$. Clearly the homogeneity of the Carnot--Carath\'eodory distance and \eqref{eq:rumin-metric-mass} imply that for $\|T\|$-almost every $\xi$,
    \begin{equation}
      \label{eq:control-mass-blowup}
      \|(\eta_{\xi, \rho})_\# T\|(\mathscr{B}_{R})
      = \frac{\|T\|(\mathscr{B}_{\rho R}(\xi))}{\rho^k}
      \xlongrightarrow{\rho \searrow 0} \omega_k R^n \Theta^k_{d_\mathrm{CC}}(\|T\|, \xi) < \infty \qquad \forall R > 0
    \end{equation}
    and \eqref{eq:convergence-away-cone} gives also that
    \begin{equation}
      \label{eq:convergence-away-cone-blowup}
      \lim_{r \to 0} \|(\eta_{\xi, \rho})_\# T\|(\mathscr{B}_R(0) \setminus X(0, \pi, s)) = 0
    \end{equation}
    for every $R > 0$ and every $0 < s < 1$. Now, given any sequence $\rho_i \searrow 0$, \eqref{eq:control-mass-blowup} allows us to take a subsequence $\rho_{i'} \searrow 0$ such that $(\eta_{\xi, \rho_{i'}})_\# T \weakto T_0$ in $\mathcal{D}_k(\Heis^n)$, for an integral current $T_0 \in \IC^\hor_k(\Heis^n)$. Then \eqref{eq:convergence-away-cone-blowup} and the lower semicontinuity of the mass implies that $\|T_0\|(X(0, \pi, s)) = 0$ for every $0 < s < 1$ and hence $T_0$ is supported in $\pi$. Assuming also that $\xi \notin \spt \del T$, we have that $\del T_0 = 0$, so by the constancy theorem $T_0 = Q \DBrack{\vec{\pi}} = Q \DBrack{\vec{T}(\xi)}$ for some integer $Q$.

    Now it is a standard fact that, if $\xi$ is a Lebesgue point for $\vec{T}\theta$ with respect to the mesure $\Haus^k_{d_\mathrm{CC}} \res S$, then $\|(\eta_{\xi,\rho_{i'}})_\# T\| \weakto Q \Haus^k \res \pi$ in the sense of measures, from which the equality $Q = \Theta^k_{d_\mathrm{CC}}(\|T\|, \xi)$ is obvious by \eqref{eq:control-mass-blowup}. To see this, let $f \in C_c(\Heis^n)$ let $\omega$ be the constant $n$-form $\omega = \langle \cdot, \vec{T}(\xi) \rangle$. Then
    \begin{align*}
      Q \int_{\pi} f \, \dd \Haus^k
      &= Q \int_{\pi} f \langle \vec{T}(\xi), \vec{T}(\xi)\rangle \, \dd \Haus^k \\
      &= T_0 (\omega f) 
      = \lim_{i' \to \infty} (\eta_{\xi,\rho_{i'}})_\# T (\omega f) \\
      &= \lim_{i' \to \infty} \rho_{i'}^{-k} \int_S f(\eta_{\xi,\rho_{i'}}(\zeta)) \langle \vec{T}(\xi), \vec{T}(\zeta) \rangle \theta(\zeta) \, \dd \Haus^k_{d_\mathrm{CC}}(\zeta) \\
      &\overset{(\star)}{=} \lim_{i' \to \infty} \rho_{i'}^{-k} \int_S f(\eta_{\xi,\rho_{i'}}(\zeta)) \langle \vec{T}(\xi), \vec{T}(\xi) \rangle\theta(\xi) \, \dd \Haus^k_{d_\mathrm{CC}}(\zeta) \\
      &= \theta(\xi) \lim_{i' \to \infty} \|(\eta_{\xi,\rho_{i'}})_\# T\|(f),
    \end{align*}
    where in $(\star)$ we have used that the following error term vanishes by the Lebesgue point hypothesis:
    \begin{align*}
      &\lim_{i' \to \infty} \rho_{i'}^{-k} \left| \int_S f(\eta_{\xi,\rho_{i'}}(\zeta)) \langle \vec{T}(\xi), \theta(\zeta) \vec{T}(\zeta) - \theta(\xi) \vec{T}(\xi) \rangle \, \dd \Haus^k_{d_\mathrm{CC}}(\zeta) \right| \\
      &\leq \sup |f| \lim_{i' \to \infty} \rho_{i'}^{-k} \int_S |\theta(\zeta) \vec{T}(\zeta) - \theta(\xi) \vec{T}(\xi) | \, \dd \Haus^k_{d_\mathrm{CC}}(\zeta) = 0.
    \end{align*}
    So far we have shown that $(\eta_{\xi,\rho_{i'}})_\# T \weakto \Theta^k_{d_\mathrm{CC}}(\|T\|, \xi) \DBrack{\vec{T}(\xi)}$ for a subsequence $i' \to \infty$, but an easy standard argument extends it to the whole sequence. 
  \end{proof}
\end{prop}

\begin{rmk}
  \label{rmk:recast-metric-currents}
  The existence and regularity theory developed in this paper for the Plateau problem in the Heisenberg group within the context of horizontal Federer--Fleming currents translates directly into metric currents thanks to \cref{thm:rumin-metric-currents}. In particular, for any $S \in \ICmet_{n-1}(\Heis^n, d_\mathrm{CC})$ with $\del S = 0$ there exists $T \in \ICmet_{n}(\Heis^n, d_\mathrm{CC})$ with $\del T = S$ which minimizes the mass among all such currents and such that, for an open set $\mathcal{U}$ such that $\spt T \cap \mathcal{U}$ is dense in $\spt T \setminus \spt S$, $\spt T \cap \mathcal{U}$ is a real analytic Legendrian submanifold.
\end{rmk}

\subsection{Results from metric geometry}

Here we state two results for horizontal Federer--Fleming currents (or Rumin currents) that were first obtained in several works by Basso, Wenger and Young by working intrinsically in $(\Heis^n, d_\mathrm{CC})$.

\begin{thm}[isoperimetric inequality]
  \label{thm:isoperimetric-heis}
  Let $1 \leq k \leq n$, $r > 0$ and $S \in \IC^\hor_{k-1}(\Heis^n)$ with $\spt S \subset \overline{\BHeis}_r(0)$. Then there exists $T \in \IC^\hor_{k}(\Heis^n)$ with $\spt T \subset \overline{\BHeis}_{r+C_\idxiso \Mass(S)^{1/(k-1)}}(0)$ such that $\del T = S$ and
  \begin{equation}
    \Mass(T) \leq C_\idxiso \Mass(S)^{\frac{k}{k-1}},
  \end{equation}
  where $C_\idxiso = C_\idxiso(n, k)$.
  \begin{proof}
    This follows from the discussion in \cite[pp.~5-6]{basso-wenger-young-fillings} for Ambrosio--Kirchheim metric currents in $(\Heis^n, d_{\mathrm{CC}})$ and hence for horizontal currents thanks to \cref{thm:rumin-metric-currents}. The bound in the support is given in \cite[Lemma 3.4]{wenger-isoperimetric-cone}. 

    As explained in the cited paper, the result that we need (for compactly supported currents) actually follows from the earlier work of Young \cite{young-filling-riemannian} and Wenger \cite{wenger-isoperimetric-asymp}.
  \end{proof}
\end{thm}

\begin{thm}[flat convergence]
  \label{thm:flat-convergence}
  Let $1 \leq k \leq n$ and $R_j \in \IC^\hor_{k-1}(\Heis^n)$ be a sequence of horizontal currents with $\del R_j = 0$ and $\spt R_j \subset K$ for some compact set $K \subset \Heis^n$. Suppose that $R_j$ have uniformly bounded masses and converge weakly (in the Federer--Fleming sense) to a current $R$. Then there exist currents $S_j \in \IC^\hor_{k}(\Heis^n)$ and positive real numbers $s_j \to 0$ such that
  \[
    \del S_j = R_j - R,
    \qquad
    \Mass(S_j) \to 0
    \qquad \text{and} \qquad
    \spt S_j \subset \BHeis_{s_j}(K).
  \]
  \begin{proof}
    We will deduce this from the main theorem in \cite{wenger-flat-convergence}. It is clear that $(\Heis^n, d_\mathrm{CC})$ is complete, quasiconvex (in fact it is a geodesic space), and again by \cite[pp.~5-6]{basso-wenger-young-fillings}, $(\Heis^n, d_\mathrm{CC})$ enjoys coning inequalities for $\ICmet_{k'}(\Heis^n, d_\mathrm{CC})$ for each $0 \leq k' \leq n-1$.

    Let $\tilde{R}, \tilde{R}_{j} \in \ICmet_{k-1}(\Heis^n, d_\mathrm{CC})$ be the metric currents corresponding to $R, R_j$ via \cref{thm:rumin-metric-currents}. If we can show that $\tilde{R}_j \weakto \tilde{R}$, then \cite[Theorem~1.4]{wenger-flat-convergence} will give us metric currents $\tilde{S}_j \in \ICmet_{k}(\Heis^n, d_\mathrm{CC})$ with $\del \tilde{S}_j = \tilde{R}_j - \tilde{R}$ and $\Mass(\tilde{S}_j) \to 0$. These facts translate automatically into the corresponding ones for $S_j := I_\# \tilde{S}_j$ and $R_j = I_\# \tilde{R}_j$, and the control on the support comes once more from \cite[Lemma~3.4]{wenger-isoperimetric-cone}.

    Thus we just need to show that, given $f_1, f_2, \ldots, f_k \in \Lip(\Heis^n, d_\mathrm{CC})$ with $f_1$ bounded, 
    \begin{equation}
      \label{eq:weak-convergence-ak}
      \tilde{R}_j(f_1 \, \dd f_2 \wedge \cdots \wedge \dd f_k)
      \xrightarrow{j \to \infty} \tilde{R}(f_1 \, \dd f_2 \wedge \cdots \wedge \dd f_k).
    \end{equation}
    Following \cite[Lemma~7.4]{williams-currents-carnot}, we can approximate the functions $f_i$ by smooth functions $f_i^\epsilon \in C^\infty_c(\Heis^n)$ by convolution (the support of $f_i^\epsilon$ may be made compact by cutting off the functions away from $K$). In particular,
    \[
      \Lip_{d_\mathrm{CC}}(f_i^\epsilon) \leq C
      \qquad \text{and} \qquad
      f_i^\epsilon \xrightarrow{\epsilon \to 0} f_i \text{ uniformly in a neighborhood of } K.
    \]
    Now we write
    \begin{align*}
      (\tilde{R}_j-\tilde{R})(f_1 \, \dd f_2 \wedge \cdots \wedge \dd f_k)
      &= (\tilde{R}_j-\tilde{R})((f_1 - f_1^\epsilon) \, \dd f_2 \wedge \cdots \wedge \dd f_k) \\
      &\quad + (\tilde{R}_j-\tilde{R})(f_1^\epsilon \, \dd (f_2 - f_2^\epsilon) \wedge \cdots \wedge \dd f_k) + \cdots \\
      &\quad + (\tilde{R}_j-\tilde{R})(f_1^\epsilon \, \dd f_2^\epsilon \wedge \cdots \wedge \dd (f_k - f_k^\epsilon)) \\
      &\quad + (\tilde{R}_j-\tilde{R})(f_1^\epsilon \, \dd f_2^\epsilon \wedge \cdots \wedge \dd f_k^\epsilon)
    \end{align*}
    and observe that, since $\del (\tilde{R}_j - \tilde{R}) = 0$,
    \begin{align*}
      &(\tilde{R}_j-\tilde{R})(f_1^\epsilon \, \dd f_2^\epsilon \wedge \cdots \wedge \dd (f_i - f_i^\epsilon) \wedge \cdots \wedge \dd f_k) \\
      &= (\tilde{R}_j-\tilde{R})(\dd f_2^\epsilon \wedge \cdots \wedge \dd (f_1^\epsilon (f_i - f_i^\epsilon)) \wedge \cdots \wedge \dd f_k)
      - (\tilde{R}_j-\tilde{R})((f_i - f_i^\epsilon) \, \dd f_2^\epsilon \wedge \cdots \wedge \dd f_1^\epsilon \wedge \cdots \wedge \dd f_k) \\
      &= - (\tilde{R}_j-\tilde{R})((f_i - f_i^\epsilon) \, \dd f_2^\epsilon \wedge \cdots \wedge \dd f_1^\epsilon \wedge \cdots \wedge \dd f_k).
    \end{align*}
    Therefore
    \begin{align*}
      |(\tilde{R}_j-\tilde{R})(f_1 \, \dd f_2 \wedge \cdots \wedge \dd f_k)|
      &\leq |(\tilde{R}_j-\tilde{R})((f_1 - f_1^\epsilon) \, \dd f_2 \wedge \cdots \wedge \dd f_k)| \\
      &\quad + |(\tilde{R}_j-\tilde{R})((f_2 - f_2^\epsilon) \, \dd f_1^\epsilon \wedge \cdots \wedge \dd f_k)| + \cdots \\
      &\quad + |(\tilde{R}_j-\tilde{R})((f_k - f_k^\epsilon) \, \dd f_2^\epsilon \wedge \cdots \wedge \dd f_1^\epsilon)| \\
      &\quad + |(\tilde{R}_j-\tilde{R})(f_1^\epsilon \, \dd f_2^\epsilon \wedge \cdots \wedge \dd f_k^\epsilon)| \\
      &\leq C (\Mass(\tilde{R}_j) + \Mass(\tilde{R})) \left( \sup_K |f_1 - f_1^\epsilon| + \cdots + \sup_K |f_k - f_k^\epsilon|\right) \\
      &\quad + |(\tilde{R}_j-\tilde{R})(f_1^\epsilon \, \dd f_2^\epsilon \wedge \cdots \wedge \dd f_k^\epsilon)|.
    \end{align*}
    Finally, since $f_1^\epsilon \, \dd f_2^\epsilon \wedge \cdots \wedge \dd f_k^\epsilon \in \mathcal{D}^{k-1}(\Heis^n)$ is a smooth compactly supported $(k-1)$-form, we may replace the last term by $|(R_j-R)(f_1^\epsilon \, \dd f_2^\epsilon \wedge \cdots \wedge \dd f_k^\epsilon)|$, which converges to zero when we send $j \to \infty$ by our hypothesis. Then \eqref{eq:weak-convergence-ak} follows after letting $\epsilon \searrow 0$.
  \end{proof}
\end{thm}

\section{Absence of monotonicity formula in higher dimensions}
\label{sec:no-monotonicity}

In this section we show that the phenomenon of monotonicity of the area density of a \emph{Hamiltonian-stationary} smooth $n$-dimensional Legendrian submanifold of $\Heis^n$ is special of dimension $n=2$. In fact, the surface that we produce does not have a density lower bound on large scales and thus cannot be globally area-minimizing as those considered in this paper.

Let $\T^n := \R^n / (2\pi \Z)^n$ and $\phi : \T^n \into \C^n$ be the standard parametrization of the generalized Clifford torus $\S^1 \times \cdots \times \S^1 \subset \C^n \simeq \R^{n} \times \R^n$ (translated by $(-1, \ldots, -1)$ for convenience), which is an isometry:
\[
  \phi(t_1, \ldots, t_n) = (\cos t_1 - 1, \ldots, \cos t_n - 1, \sin t_1, \ldots, \sin t_n).
\]
It is well known that this surface is Hamiltonian-stationary (see \cite{oh-hamiltonian-stationary}) but not exact. However, it has an exact covering which lifts to $\Heis^n$ as
\begin{gather*}
  \tilde{\phi} : \mathcal{C} := \R^n / \langle 2\pi(e_1-e_2), \ldots 2\pi (e_1 - e_n)\rangle_\Z \longrightarrow \Heis^n \\
                 (t_1, \ldots, t_n)  \longmapsto (\cos t_1 - 1, \ldots, \cos t_n - 1, \sin t_1, \ldots, \sin t_n, \varphi(t_1, \ldots, t_n))
\end{gather*}
where $\varphi(t_1, \ldots, t_n) = \tfrac{1}{2} (t_1 + \cdots + t_n - \sin t_1 - \cdots - \sin t_n)$ is well defined and satisfies
\[
  \dd \varphi = \frac{1}{2} \left( \dd t_1 + \cdots + \dd t_n - \cos t_1 \dd t_1 - \cdots - \cos t_n \dd t_n \right) = \phi^* \lambda
\]
for the Liouville form $\lambda = \tfrac{1}{2} (\Vec{x} \cdot \dd \Vec{y} - \Vec{y} \cdot \dd \Vec{x})$. It is clear that $\tilde{\phi}$ is a Legendrian embedding, which is still an isometry and Hamiltonian-stationary because these two properties are local and preserved by Legendrian lifts.

For $r$ small, $\tilde{\phi}^{-1}(\BHeis_r) \subset \{ t_1^2 + \cdots + t_n^2 < r^2 \}\subset \tilde{\phi}^{-1}(\BHeis_{r+o(r)})$, hence $\vol(\tilde{\phi}(\mathcal{C}) \cap \BHeis_r) = \omega_n r^n + o(r^n)$.

To understand the large scale behavior, first observe that since $\phi$ is bounded, $\tilde{\phi}^{-1}(\BHeis_r) \subset \{ 4 |\varphi| < r^2 \} \subset \tilde{\phi}^{-1}(\BHeis_{r+o(r)})$ and as a result $\tilde{\phi}^{-1}(\BHeis_{r-o(r)}) \subset \left\{ |t_1 + \cdots + t_n| < \tfrac{1}{2} r^2 \right\} \subset \tilde{\phi}^{-1}(\BHeis_{r+o(r)})$ for $r$ large. To compute the volume of this region, we use the fundamental domain $\mathcal{C}_0 := \R \times [0, 2\pi) \times \cdots \times [0, 2\pi)$ of $\mathcal{C}$. Then
\[
  \vol(\tilde{\phi}(\mathcal{C}) \cap \BHeis_r)
  \sim
  \vol \left(\left\{ (t_1, t_2, \ldots, t_n) \in \mathcal{C}_0 : |t_1 + \cdots + t_n| < \frac{1}{2} r^2 \right\} \right)
  \sim (2\pi)^{n-1} r^2,
\]
which shows that the quotient $(\omega_n r^n)^{-1} \vol(\tilde{\phi}(\mathcal{C}) \cap \BHeis_r)$ actually tends to zero for large $r$ when $n > 2$.

\printbibliography
\end{document}